\documentclass[12pt]{amsart} 
\usepackage[normalem]{ulem}  
\usepackage{lmodern}
\usepackage{amsmath}
\usepackage{amssymb}
\usepackage{dsfont}  % to do a doublestroke 1
\usepackage{hyperref}
\usepackage{graphicx}
\usepackage[all]{xy}
\usepackage{latexsym}
\usepackage{epstopdf}
\usepackage{pgf}
\usepackage{tikz}
\usepackage{tikz-cd}
\usepackage{pgf}
\usepackage{float}
\usepackage{subcaption} 
\usepackage{xcolor}
\usepackage{algorithm2e}
\usepackage{color}
\usepackage{ytableau}
\usepackage{adjustbox}
\usepackage{mathrsfs}
\usepackage{enumerate}   
\usepackage{mathbbol}
\usepackage{quiver}

\usetikzlibrary{matrix,arrows,backgrounds,shapes.misc,shapes.geometric,patterns,calc,positioning}
\usetikzlibrary{calc,shapes}
\usetikzlibrary{decorations.pathmorphing}
\usetikzlibrary{positioning,arrows}
\usetikzlibrary{decorations.markings}

\usepackage[left=3cm,right=3cm,top=2.8cm,bottom=2.8cm]{geometry}
\newtheorem{theorem}{Theorem}[section]

\newtheorem{lemma}[theorem]{Lemma}
\newtheorem{corollary}[theorem]{Corollary}

\newtheorem{proposition}[theorem]{Proposition} 
\newtheorem{conjecture}[theorem]{Conjecture}

\theoremstyle{remark}
\newtheorem{remark}[theorem]{Remark}

\theoremstyle{definition}
\newtheorem{definition}[theorem]{Definition} 
\newtheorem{example}[theorem]{Example} 

\numberwithin{equation}{section}

\newcommand{\rad}{\operatorname{rad}\nolimits}

\newcommand{\Hom}{\operatorname{Hom}\nolimits}

\newcommand{\Ext}{\operatorname{Ext}\nolimits}
\newcommand{\End}{\operatorname{End}\nolimits}

\newcommand{\add}{\operatorname{add}\nolimits}

\newcommand{\CC}{\mathbb{C}}
\newcommand{\PP}{\mathbb{P}}
\newcommand{\QQ}{\mathbb{Q}}
\newcommand{\ZZ}{\mathbb{Z}}

\newcommand{\Gr}{\mathrm{Gr}}

\newcommand{\SSYT}{{\rm SSYT}}

\newcommand{\ch}{{\rm ch}}
\newcommand{\sgn}{{\rm sgn}}
\newcommand{\TT}{\mathbb{T}}
\newcommand{\bT}{\mathbf{T}}

\DeclareMathOperator{\diag}{diag}

\newcommand{\thfrac}[3]{ 	
\begin{aligned} #1\\ \hline \\[-3\jot] #2 \\ \hline \\[-3\jot] #3 \end{aligned}
}
\newcommand{\ffrac}[4]{ 	
\begin{aligned} #1\\ \hline \\[-3\jot] #2 \\ \hline \\[-3\jot] #3 \\ \hline \\[-3\jot] #4 \end{aligned}
}

\newcommand\scalemath[2]{\scalebox{#1}{\mbox{\ensuremath{\displaystyle #2}}}}

\begin{document}

%\title{Additive and monoidal categorifications of cluster algebras}
%\title{Rigid modules and real modules}
\title[Additive and monoidal categorifications]{A correspondence between additive and monoidal categorifications with application to Grassmannian cluster categories}
\author{Karin Baur}
\address{School of Mathematics, University of Leeds, Leeds, LS2 9JT, United Kingdom}
\address{Faculty for Mathematics, Ruhr University Bochum, Germany}
\email{karin.baur@rub.de}
\author{Changjian Fu}
\address{Department of Mathematics, Sichuan University, Chengdu 610064, P.R.China}
\email{changjianfu@scu.edu.cn}
\author{Jian-Rong Li}
\address{Faculty of Mathematics, University of Vienna, Oskar-Morgenstern Platz 1, 1090 Vienna, Austria}
\email{lijr07@gmail.com}
\date{}

\begin{abstract}
Building on work of Derksen--Fei and Plamondon, we formulate a conjectural correspondence between additive and monoidal categorifications of cluster algebras, which reveals a new connection between the additive reachability conjecture and the multiplicative reachability conjecture. Evidence for this conjecture is provided by results on Grassmannian cluster algebras and categories in the tame types. Moreover, we give a construction of the generic kernels introduced by Hernandez and Leclerc for type $\mathbb{A}$ via the Grassmannian cluster categories. As an application of the correspondence, we construct rigid and non-rigid indecomposable modules in Grassmannian cluster categories. 
\end{abstract}

\maketitle

\setcounter{tocdepth}{1}
\tableofcontents 

\section{Introduction}

\subsection{Cluster algebras} Fomin and Zelevinsky introduced cluster algebras around 2000 \cite{FZ02} to give a combinatorial framework for the study of canonical
bases in quantum groups \cite{Kas91, Lus90} and total positivity in algebraic groups \cite{Lus94}. A cluster algebra is a commutative algebra defined with some initial data called an initial seed and using a procedure called mutation. Cluster algebras have connections and applications to various areas in mathematics and physics, for example, finite-dimensional algebras \cite{BBM06, Kel10}, KLR algebras \cite{KKKO}, Lie theory \cite{BFZ05, GLS13}, %, HL10}, 
mirror symmetry \cite{GHKK18}, Poisson geometry \cite{GSV03}, quantum affine algebras \cite{HL10}, scattering amplitudes and amplituhedron in physics \cite{ABCGPT16, ALS21, PSW22, DFGK20, HP21}, stability conditions in the sense of Bridgeland \cite{Bri07}, symplectic geometry \cite{CG22}, and Teichm\"{u}ller theory \cite{FG06}. 

In the development of cluster algebras, categorification has played an important role. Here,
categorification is a process of finding category-theoretic analogs of set-theoretic concepts by replacing sets with categories, functions with functors, and equations between functions by natural isomorphisms between functors. Categofications provide a deeper understanding of the original objects, see for example \cite{LS22}. There are two types of categorifications of cluster algebras: one is called additive categorification \cite{BMRRT06,CCS06} and the other is called monoidal categorification \cite{HL10}. 

\subsection{Additive categorification} 
In \cite{BMRRT06}, Buan, Marsh, Reineke, Reiten, and Todorov introduced cluster categories as a quotient of the bounded derived category of the module category of a finite-dimensional hereditary algebra. Independently, Caldero, Chapoton, and Schiffler introduced cluster categories in the type $A_n$ case in~\cite{CCS06}. Their categories are equivalent to the ones of \cite{BMRRT06} for type $A$. In the simply-laced Dynkin types, these categories categorify the corresponding cluster algebras. 

Categorifications of (acyclic) cluster algebras %associated with an acyclic quiver 
%using the cluster category associated with 
via the path algebra of the underlying quiver
were studied in \cite{BMR07, BMR08, BMRT07, CC06, CK06, CK08, MRZ03}, for example. Gei{\ss}, Leclerc, and Schr\"{o}er \cite{GLS06, GLS07, GLS08} categorified cluster algebras using subcategories of the category of modules over a preprojective algebra associated to an acyclic quiver. Buan, Iyama, Reiten, and Scott \cite{BIRS09}   
introduced the notion of cluster structure and constructed 2-Calabi–Yau categories with cluster structures related to preprojective algebras of non-Dynkin quivers associated with elements in the Coxeter group. 

Denote by $\Gr(k,n)$ the Grassmannian, of $k$-planes in $\mathbb{C}^n$ and $\CC[\Gr(k,n)]$ its homogeneous coordinate ring. It was shown by Scott \cite{Sco} that the ring $\CC[\Gr(k,n)]$ has a cluster algebra structure. In \cite{JKS16}, Jensen, King and Su
introduced the Grassmannian cluster category ${\rm CM}(B_{k,n})$ of maximal Cohen--Macaulay modules over an algebra $B_{k,n}$ which provided an additive categorification of $\CC[\Gr(k,n)]$.
It can be viewed as a completion of the category Sub$Q_k$ introduced by Geiss, Leclerc and Schr\"oer \cite{GLS08} by an additional projective-injective object. 

%%%%%%%%%
%
\subsection{Monoidal categorification} Hernandez and Leclerc introduced monoidal categorifications of cluster algebras in \cite{HL10}. For details we refer to Section~\ref{ss:monoidal-cat} below.

Since their appearance, monoidal categorifications of cluster algebras have been studied intensively, see e.g.~\cite{DS23, KKKO, KKOP22, HL10, HL15, HL16, Nak11, Qin17}. 

In this article, we are interested in the categories which first appeared in~\cite{HL10}. 
For every complex simple Lie algebra $\mathfrak{g}$ and every $\ell \ge 0$, Hernandez and Leclerc introduced a 
subcategory $\mathscr{C}_{\ell}^{\mathfrak{g}}$ of the category $\mathscr{C}$ of finite dimensional $U_q(\widehat{\mathfrak{g}})$-modules. It provides a monoidal categorification of a certain cluster algebra $\mathcal{A}_{\ell}^{\mathfrak{g}}$. Roughly speaking, the $q$-character $\chi_q$ of a simple module in $\mathscr{C}_{\ell}^{\mathfrak{g}}$ gives rise to a dual canonical basis element of the cluster algebra $\mathcal{A}_{\ell}^{\mathfrak{g}}$, see Section~\ref{ss:monoidal-cat} below.  In the case of $\mathfrak{g}=\mathfrak{sl}_k$, the cluster algebra $\mathcal{A}_{\ell}^{\mathfrak{sl}_k}$ is isomorphic to a certain quotient $\CC[\Gr(k,n,\sim)]$ of $\CC[\Gr(k,n)]$, for $n=\ell+k+1$, see Section 13.9 in \cite{HL10}.  This fact provides a link with the additive categorification associated with the Grassmannian which we will exploit. 

\subsection{Correspondence between additive and monoidal categorifications}
A natural question is to see how additive categorifications and monoidal categorifications of cluster algebras are related. In \cite{KKOP21b}, Kashiwara-Kim-Oh-Park constructed cluster algebra structures on the Grothendieck rings of certain monoidal subcategories of the category of finite-dimensional representations of a quantum loop algebra, generalizing Hernandez-Leclerc's pioneering work in \cite{HL10}. In \cite{Con24}, by studying the relation between the representation theory of quantum groups (monoidal categorification side) and that of quivers with relations (additive categorification side), Contu solved an open problem proposed by Kashiwara-Kim-Oh-Park in \cite{KKOP21b}: finding explicit quivers for the seeds of the cluster algebras in \cite{KKOP21b}. In \cite{Fujita2024}, Fujita studied the singularities of normalized R-matrices between arbitrary simple modules over the quantum loop algebra of type ADE in Hernandez-Leclerc's level-one subcategory $\mathscr{C}_1$ \cite{HL10}. He showed that the pole orders of these R-matrices coincide with the dimensions of $E$-invariants between the corresponding decorated representations of Dynkin quivers. His result can be seen as a correspondence of numerical characteristics between additive and monoidal categorifications of cluster algebras of finite ADE type.

The goal of this paper is to study the correspondence between additive and monoidal categorifications. In particular, we study the relations between the category %relation between additive categorifications and monoidal categorifications of cluster algebras. In particular, we study relations between the category 
${\rm CM}(B_{k,n})$ and the category $\mathscr{C}_{\ell}^{\mathfrak{sl}_k}$. 

Following \cite{P13}, for every simple $U_q(\hat{\mathfrak{g}})$-module $L(M)$ with dominant monomial $M$, we associate an open dense subset $\mathcal{O}_{\mathbf{g}_{M}}$ of an affine space  related to the $\mathbf{g}$-vector $\mathbf{g}_{M}\in \mathbb{Z}^n$ of $L(M)$. Applying a construction of Derksen--Fei \cite{DFGK20}, we get the so-called symmetrized  $E$-invariant function $\mathbb{E}(-,-)$ on $\mathcal{O}_\mathbf{g}\times \mathcal{O}_{\mathbf{h}}$ and its `degenerate' version $\mathfrak{e}(-)$ 
on $\mathcal{O}_\mathbf{g}$, for integer vectors $\mathbf{g},\mathbf{h}\in \mathbb{Z}^n$. Denote by $\mathfrak{e}(\mathbf{g},\mathbf{h})$ the generic value of $\mathbb{E}(-,-)$ on $\mathcal{O}_\mathbf{g}\times \mathcal{O}_{\mathbf{h}}$ and by $\mathfrak{e}(\mathbf{g})$ the generic value of $\mathfrak{e}(-)$  on $\mathcal{O}_\mathbf{g}$, see Section~\ref{ss:function-e} for details. 

We expect that $\mathfrak{e}(-,-)$ and $\mathfrak{e}(-)$ provide links between the additive and monoidal categorifications: we conjecture that $L(M)$ is real if and only if $\mathfrak{e}(\mathbf{g}_M)=0$ (Conjecture \ref{c:monoidal-additive-correspondence}). Furthermore, for real dominant monomials $M$ and $N$, we conjecture that $\chi_q(L(M))\chi_q(L(N))=\chi_q(L(MN))$ if and only if $\mathfrak{e}(\mathbf{g}_M,\mathbf{g}_N)=0$ (cf. Conjecture \ref{c:monoidal-additive-corr-compatible}).
Under the assumption that Conjecture \ref{c:monoidal-additive-correspondence} holds, we show that the additive reachability conjecture in additive categorification is equivalent to the multiplicative reachability conjecture in monoidal categorification (cf. Proposition \ref{p:equivalence-add-mult}). We provide evidence for Conjecture \ref{c:monoidal-additive-correspondence} by proving it in the case where the underlying cluster algebra is of finite type, Theorem~\ref{thm:conjecture-true-finite-type}. The value of $\mathfrak{e}(-,-)$ also gives a characterization for when two given cluster variables form an exchange pair (cf. Proposition \ref{p:exchange-pair}). 
In particular, we establish the additive reachability conjecture for Grassmannian cluster categories of tame type in Section \ref{sec:Evidences in Grassmannian cluster categories of tame type}, see Theorem~\ref{thm:CMGr39-48-connected}.  {This gives further 
evidence for Conjecture \ref{c:monoidal-additive-correspondence}.}
%In particular, we establish the additive reachability conjecture for Grassmannian cluster categories of tame type, see Theorem~\ref{thm:CMGr39-48-connected}.
%\ref{thm:CMGr39 amd CMGr48 are connected}.

%%
%
\subsection{From real prime modules to rigid indecomposable modules}

We show that there are one to one correspondences between the isomorphism classes of reachable rigid indecomposable modules in ${\rm CM}(B_{k,n})$, the set of reachable prime real semistandard Young tableaux in $\SSYT(k, [n])$, and the set of isomorphism classes of reachable prime real modules in $\mathscr{C}_{\ell}^{\mathfrak{sl}_k}$, where $n=k+\ell+1$, see Theorem~\ref{thm: cluster variables in Grkn and CMBkn}(1) and (2). As an application, given a reachable real prime tableau or a reachable real prime module, we are able to compute its ${\bf g}$-vector. From that, we can find
the corresponding rigid indecomposable module in ${\rm CM}(B_{k,n})$ explicitly, see Section \ref{sec:correspondence between real modules and rigid modules}.

\subsection{Hernandez and Leclerc's generic kernels}

Let $\Gamma^-$ be the quiver of type $\mathbb{A}_{k-1}$ introduced by Hernandez and Leclerc in \cite{HL16}, see Section~\ref{subsec:Application to generic kernels of Hernandez-Leclerc type A}. Denote by $V^-$ its vertex set. Denote by $W$ the formal sum of all oriented $3$-cycles (up to cyclic permutations) and $A=J(\Gamma^-,W)$ the Jacobian algebra associated with $(\Gamma^-,W)$. Let $I(i,m)$ be the indecomposable injective $A$-module associated with vertex $(i,m)$ of $\Gamma^-$. Hernandez and Leclerc proved in~\cite{HL16} that, for any $v \in \ZZ_{\ge 1}$, there is a morphism $f\in \Hom_A(I(i,m),I(i,m-2v))$ such that $\ker f$ is finite dimensional. Such a morphism $f$ is referred to as a \textit{generic homomorphism} and $K_{v,m}^{(i)}:=\ker f$ is referred as a \textit{generic kernel}.

For any positive integer $\ell$, let $\Gamma^{-}_{-2\ell}$ be the (finite) full subquiver of $\Gamma^-$ consisting of the vertex set $\{(i,m)\in V^-~|~m\geq -2\ell\}$.  Denote by $A_{-2\ell}$ the associated Jacobian algebra. Set $n=k+\ell+1$. Let $\underline{\rm CM}(B_{k,n})$ be the stable category of the Frobenius category ${\rm CM}(B_{k,n})$, and $\tau$ the Auslander-Reiten translaton of $\underline{\rm CM}(B_{k,n})$.
We construct an explicit cluster-tilting object $L=\bigoplus\limits_{(i,m)\in \Gamma^-_{-2\ell}}L_{i,m}$ of $\underline{\rm CM}(B_{k,n})$ such that $A_{-2\ell}\cong \End_{\underline{\rm CM}(B_{k,n})}(L)$.
For every vertex $(i,m)$ of $\Gamma_{-2\ell}^-$ and every positive integer $v \le \frac{m+2\ell+(-1)^{i+1}}{2}$, we introduce an indecomposable rigid module $M_{v,m}^{(i)}$ in $\operatorname{CM}(B_{k,n})$.
% {\color{blue}We prove that for every vertex $(i,m)$ of $\Gamma_{-2\ell}^-$ and every positive integer $v\le \frac{m+2\ell+(-1)^{i+1}}{2}$ there is a unique indecomposable rigid module $M_{v,m}^{(i)}$ in ${\rm CM}(B_{k,n})$ fitting into a triangle with $L_{i,m}\to L_{i,m-2v}$ in  $\underline{\rm CM}(B_{k,n})$ 
% %
% %\begin{align*}
% %\tau T_{i,m} \to \tau T_{i,m-2v} \to \tau M_{v,m}^{(i)}\to \tau^2T_{i,m},
% %\end{align*}
% %
% \begin{align*}
% L_{i,m} \to L_{i,m-2v} \to M_{v,m}^{(i)}\to \tau L_{i,m},
% \end{align*}
% where $\tau$ is the Auslander-Reiten translation of $\underline{\rm CM}(B_{k,n})$.} 
Both $M_{v,m}^{(i)}$ and $\tau M_{v,m}^{(i)}$ are rank one modules and we provide explicit descriptions of them, 
%we give explicit %expression of $M_{v,m}^{(i)}$ and $\tau M_{v,m}^{(i)}$, 
see Theorem~\ref{thm:generic-kernel-grass}. 
For $(i,m)\in V^-$ and $v\geq 1$, in order to determine $K_{v,m}^{(i)}$, it suffices to take $\ell\geq \frac{2v+2-m}{2}$. By applying Theorem~\ref{thm:generic-kernel-grass}, we prove that $K_{v,m}^{(i)}\cong \Hom_{\underline{\rm CM}(B_{k,n})}(L,\tau M_{v,m}^{(i)})$, see Theorem \ref{thm:realization-generic-kernel}.  In particular, this  yields an explicit and computable construction of Hernandez and Leclerc's generic kernels in the type $\mathbb{A}$ case. As an application, this  gives a criterion for checking whether two cluster variables corresponding to Kirillov-Reshetikhin modules belong to a common cluster, see Remark~\ref{rem:criterion-compatibility}.

\subsection{New indecomposable modules in Grassmannian cluster categories}

We use mutation of tableaux from \cite[Section 4]{CDFL20} to obtain additional cluster variables in $\CC[\Gr(k,n)]$. 
From these cluster variables, using the results of Section~\ref{sec:correspondence between real modules and rigid modules}, we compute explicit rigid indecomposable modules in the two tame types ${\rm CM}(B_{3,9})$ and in ${\rm CM}(B_{4,8})$ up to rank $4$, see Sections~\ref{subsec:rigid indecomposable modules Gr39} 
and~\ref{subsec:rigid indecomposable modules Gr48}. We conjecture
that these are all rigid indecomposable modules up to rank $4$ in ${\rm CM}(B_{3,9})$ and in ${\rm CM}(B_{4,8})$. 

We then consider non-real 
$U_q(\widehat{\mathfrak{sl}_3})$-modules and their semistandard Young tableaux to find non-rigid objects of these two  Grassmannian cluster categories: 

The $U_q(\widehat{\mathfrak{sl}_3})$-module corresponding to the semistandard Young tableau $\ytableausetup{centertableaux}
\scalemath{0.6}{\begin{ytableau}
1 & 2 & 3 \\
4 & 5 & 6 \\
7 & 8 & 9
\end{ytableau}}$ is non-real.  
We associate a module in 
${\rm CM}(B_{3,9})$ to it and prove that it is non-rigid, see Section~\ref{subsec:non-rigid-module-example-Gr39}. 
In Section~\ref{subsec:non-rigid-module-example-Gr48} we establish a non-rigid module 
of ${\rm CM}(B_{4,8})$. It 
corresponds to the non-real tableau (see Definition~\ref{def:real-tableau_vector}) $\ytableausetup{centertableaux}
\scalemath{0.6}{\begin{ytableau}
1 & 2  \\
3 & 4  \\
5 & 6  \\
7 & 8
\end{ytableau}}$.
We conjecture that these two non-rigid modules are indecomposable. Finally, in Sections~\ref{subsec:non-rigid modules Gr39} and~\ref{subsec:non-rigid modules Gr48}, we apply braid group actions to construct new indecomposable non-rigid modules from the indecomposable non-rigid modules of Section \ref{subsec:non-rigid-module-example-Gr39} and \ref{subsec:non-rigid-module-example-Gr48}. 

% Furthermore, we apply braid group actions to construct indecomposable non-rigid modules in ${\rm CM}(B_{k,n})$ from known indecomposable non-rigid modules, see Sections \ref{subsec:Quasi-homomorphisms of cluster algebras and braid group actions}, \ref{subsec:non-rigid-module-example-Gr39}, and \ref{subsec:non-rigid-module-example-Gr48}. 

%%
%
\subsection{Organization of the paper} 
In Section \ref{sec:preliminaries}, we recall the necessary background. 
In Section~\ref{sec:correspondence between additive and monoid categorifications}, we study the link between additive and monoidal categorifications. In Section \ref{sec:Evidences in Grassmannian cluster categories of tame type}, we establish the additive reachability conjecture for Grassmannian cluster categories of tame type and demonstrate evidence for Conjecture \ref{c:monoidal-additive-correspondence} characterising real simple $U_q(\hat{\mathfrak{g}})$-modules. In Section \ref{sec:correspondence between real modules and rigid modules}, we give an explicit method to construct reachable rigid indecomposable modules from reachable real prime modules. In Section \ref{sec:realization of generic kernels of Hernandez-Leclerc}, we give a construction of the generic kernels introduced by Hernandez and Leclerc for type $\mathbb{A}$ via Grassmannian cluster categories. This gives us a criterion for checking whether two cluster variables corresponding to Kirillov-Reshetikhin modules are compatible. In Section \ref{sec:application to construction of indecomposable modules}, we construct rigid indecomposable modules and (conjectural) indecomposable non-rigid modules in Grassmannian cluster categories. 

\subsection*{Acknowledgements}
The authors would like to thank Ellis Caird and Alastair King for helpful discussions. We are also grateful to the reviewer for their insightful comments and constructive suggestions, which have significantly improved the manuscript's quality and clarity. JRL would like to thank Tianyuan Mathematical Center in Southwest China and the Department of Mathematics at Sichuan University for  hospitality where part of the work has been done. KB was supported by the EPSRC Programme Grant EP/W007509/1 and by the Royal Society Wolfson Award RSWF/R1/180004. CJF is supported by NSF of China (Nos. 11971326, 12471037, 12571040) and Fundamental Research Funds for the Central Universities. JRL is supported by the Austrian Science Fund (FWF): P-34602, Grant DOI: 10.55776/P34602, and PAT 9039323, Grant-DOI 10.55776/PAT9039323, and NSF of China (No. 12471023).

\section{Preliminaries} \label{sec:preliminaries}

In this section, we recall the results we will use in the later sections. 

%results on cluster algebras \cite{BFZ05, FZ02, FZ07} %, W14}, 
%additive categorification of cluster algebras \cite{BIRS09,FK10},  quantum affine algebras \cite{CP94, FR98}, monoidal categorifictions of cluster algebras \cite{HL10} and on Grassmannian cluster algebras \cite{Sco, CDFL20} which will be needed in the paper. 

%%
\subsection{Notation}\label{sec:notation}
For integers $a \le b$, we use the notation $[a,b]$ for $\{i: a \le i \le b\}$ and $[a]$ for $\{i: 1 \le i \le a\}$ (for positive $a$). For any integer $a\in \mathbb{Z}$, we write $[a]_+$ for  
$\max(a,0)$ 
%$[a]_+=\max(a,0)$ 
and \begin{align*}
\sgn(a)=\begin{cases}
    1, & a>0, \\
    0, & a=0, \\
    -1, & a<0.
\end{cases}    
\end{align*} 
For any integer vector $\mathbf{g}=(g_1,\dots, g_n)\in \mathbb{Z}^n$, we write 
$\mathbf{g}_+:=([g_1]_+,\dots, [g_n]_+)\in \mathbb{Z}^n$ and $\mathbf{g}_-:=([-g_1]_+,\dots, [-g_n]_+)\in \mathbb{Z}^n$.

Let $\mathcal{T}$ be an additive category. For $M\in \mathcal{T}$, we denote by $\add M$ the subcategory of $\mathcal{T}$ consisting of objects which are finite direct sums  of  direct summands of $M$.
Let  $T=\oplus_{i=1}^nT_i\in \mathcal{T}$  with indecomposable direct summands $T_1,\dots,T_n$. We say that $T$ is {\em basic} if $T_i\not\cong T_j$ for every $i\ne j$. Let $T=\oplus_{i=1}^nT_i$ be basic and let $\alpha=(a_1,\dots, a_n)$ be an integer vector with $a_i\ge 0$ for all $i$. Then we define 
\[
T^\alpha=T_1^{a_1}\oplus\cdots\oplus T_n^{a_n}, 
\]
where $T_i^{a_i}$ is the direct sum of $a_i$ copies of $T_i$.

\subsection{Cluster algebras} \label{subsec:cluster algebras}

Let $Q=(Q_0,Q_1,s,t)$ be a finite quiver (i.e. directed graph) with vertex set $Q_0$, arrow set $Q_1$ and with maps 
$s,t: Q_1 \to Q_0$ taking an arrow to its source and target, respectively. We assume that $Q$ has no loops or $2$-cycles. 
(Note that this is sometimes 
called a cluster quiver). 
We identify $Q_0$ with $[m] = \{1,\dots,m\}$. 
As part of the data of $Q$, one further declares vertices $1,\dots,n$ as mutable and vertices $n+1,\dots,m$ as frozen, for some $n\le m$: 
For $k \in [n]$, the mutated quiver $\mu_k(Q)$ is a quiver on the same vertex set as $Q$, 
% and value of $n$,
and with arrows obtained as follows:
\begin{enumerate}
\item[(i)] for each sub-quiver $i \to k \to j$, add a new arrow $i \to j$,

\item[(ii)] reverse the orientation of every arrow with target or source equal to $k$,

\item[(iii)] remove the arrows in a maximal set of pairwise disjoint $2$-cycles. 
\end{enumerate}

%Let $m \ge n$ be positive integers and
Let $\mathscr{F}$ be the 
field of rational functions in $m$ independent variables  over $\QQ$. 
 
A {\em labeled seed} in $\mathscr{F}$  is a pair $({\bf x}, Q)$ where  ${\bf x} = (x_1, \ldots, x_m)$ form a free generating set for $\mathscr{F}$, and $Q$ is a quiver on vertices $1, \ldots, m$, whose vertices $1, \ldots, n$ are mutable, and whose vertices $n+1, \ldots, m$ are frozen. The set ${\bf x}$ is called the {\em (labeled)  cluster} of the (labeled) seed $({\bf x}, Q)$, {and the quiver $Q$ is referred to as the {\em exchange quiver}.} The variables $x_1, \ldots, x_m$ are called {\em cluster variables}. Specifically, $x_1, \ldots, x_n$ are the {\em mutable variables}, while $x_{n+1}, \ldots, x_m$ are referred to as {\em frozen variables} (or {\em coefficients}).

For a seed $({\bf x}, Q)$ and $k \in [n]$, the mutated seed $\mu_k({\bf x}, Q)$ is $({\bf x}', \mu_k(Q))$, where ${\bf x}' = (x_1', \ldots, x_m')$ with $x_j'=x_j$ for $j\ne k$ and where $x_k' \in \mathscr{F}$ is determined by
\begin{align*}
x_k' x_k = \prod_{\alpha \in Q_1, s(\alpha)=k} x_{t(\alpha)} + \prod_{\alpha \in Q_1, t(\alpha)=k} x_{s(\alpha)}.
\end{align*} 

Consider the $n$-regular tree $\TT_n$ whose edges are labeled by $1, \ldots, n$, so that the $n$ edges emanating from each vertex receive different labels (note that for $n\ge 2$, the tree $\TT_n$ is an infinite graph). We will refer to $\TT_n$ as a {\em labeled tree}. A {\em cluster pattern} (of $({\bf x},Q)$) is an assignment of a labeled seed $\Sigma_t = ({\bf x}_t, Q_t)$ (which can be obtained from $({\bf x},Q)$ by iterated mutations) to every vertex $t \in \TT_n$ such that the seeds assigned to the end-points of any edge $\xymatrix{t\ar@{-}[r]^k &t'}$ are related to each other by the mutation $\mu_k$.  It is clear that a cluster pattern is uniquely determined by assigning $(\mathbf{x}, Q)$ to a fixed vertex $t_0\in \mathbb{T}_n$. 
We refer to $t_0$ {as} the {\em root vertex} and {to} $(\mathbf{x},Q)$ {as} the {\em initial seed} of this cluster pattern. For a given cluster pattern, we write 
\[
{\bf x}_t=(x_{1;t},\dots, x_{m;t})\ \text{and}\  {\bf x}_{t_0}=(x_1,\dots, x_m).
\]
to denote the cluster at $t\in \TT_n$ and the initial cluster, respectively.
Clearly, for every vertex $t\in \TT_n$, we have $x_{j;t}=x_j$ whenever $n<j\leq m$.

Given a cluster pattern of $({\bf x},Q)$, let $\mathcal{X} = \cup_{t \in \TT_n} {\bf x}_t = \{ x_{i,t}: t \in \TT_n, i \in [m] \}$. 
The elements $x_{i,t} \in \mathcal{X}$ are cluster variables. Two cluster variables are {\em compatible} if they are in a common cluster. A {\em cluster monomial} is a monomial in the cluster variables of a single cluster.
The {\em cluster algebra} $\mathcal{A}:=\mathcal{A}(Q)$ associated with the cluster pattern is the $\ZZ$-subalgebra $\ZZ[\mathcal{X}]$ of $\mathscr{F}$ generated by all cluster variables (appearing in the cluster pattern). The number $n$ of mutable variables in any cluster is called the {\em rank} of the cluster algebra $\mathcal{A}$. In this generality, $\mathcal{A}$ is called a cluster algebra defined by a quiver, or a skew-symmetric cluster algebra of geometric type. 
If $n=m$, i.e. if there are no frozen variables, we say that $\mathcal{A}$ is a {\em cluster algebra with trivial coefficients}.

\begin{remark}\label{rem:Q-gives-B}
Let $Q$ be a finite quiver with no loops or $2$-cycles whose vertices are labeled $1, \ldots, m$. Then we may encode $Q$ by an $m \times m$ skew-symmetric matrix $\widetilde{B}=\widetilde{B}(Q)=(b_{ij})_{i, j}$, where $b_{ij} = -b_{ji} = \ell$ whenever there are $\ell$ arrows from vertex $i$ to vertex $j$. The arrows between frozen vertices do not affect seed mutation. Therefore, one often omits the data corresponding to such arrows, i.e. one omits the last $m-n$ columns. The resulting matrix $\widetilde{B}$ is an $m \times n$ matrix. For $t\in \TT_n$ and $\Sigma_t=({\bf x}_t, Q_t)$ the corresponding seed, we also write this as $\Sigma_t=({\bf x}_t, \widetilde{B}_t)$, where $\widetilde{B}_t$ is the $m \times n$ matrix corresponding to $Q_t$, omitting the last $m-n$ columns. We also denote the cluster algebra $\mathcal{A}(Q)$  by $\mathcal{A}(\widetilde{B})$.
\end{remark}

We recall the upper cluster algebra introduced in~\cite{BFZ05}. Consider a cluster pattern 
$\TT_n$ and $t\in \TT_n$. We denote by $\ZZ[{\bf x}_t^{\pm}]=\ZZ[x_{n+1},\ldots, x_m][x_{1;t}^{\pm 1},\dots, x_{n;t}^{\pm 1}]$ the ring of Laurent polynomials in the cluster variables $x_{1;t},\dots, x_{n;t}$ with coefficients in $\ZZ[x_{n+1},\ldots, x_m]$ and set
\begin{align*}
    \mathcal{U}(\Sigma_t) := \ZZ[{\bf x}_t^{\pm}] \cap \ZZ[{\bf x}_{t_1}^{\pm}] \cap \cdots \cap \ZZ[{\bf x}_{t_n}^{\pm}],
\end{align*}
where $t_1,\dots, t_n\in \TT_n$ are the $n$ vertices adjacent to $t$. The {\em upper cluster algebra } 
$\overline{\mathcal{A}}:=\overline{\mathcal{A}}(\tilde{B}):=\overline{\mathcal{A}}(Q)\subset \mathcal{F}$ associated with the given cluster pattern is defined by
\begin{align*}
    \overline{\mathcal{A}}(Q) = \bigcap\limits_{t\in \TT_n} \mathcal{U}(\Sigma_t).
\end{align*}

\begin{remark}

  Alternatively, one can choose a different coefficient ring to define cluster algebras and upper cluster algebras.
  Let $\PP = {\rm Trop}(x_{n+1}, \ldots, x_m)$ be the tropical semifield generated by the frozen variables $x_{n+1}, \ldots, x_m$, i.e. the abelian group freely generated by $x_{n+1}, \dots, x_m$ with addition $\oplus$ defined by
\begin{align*}
    \prod_{i=n+1}^m x_j^{a_j} \oplus \prod_{j=n+1}^m x_j^{b_j} = \prod_{j=n+1}^m {x}_j^{\min(a_j, b_j)}.
\end{align*}
Denote by $\ZZ\PP$ the group ring of $\PP$ with integer coefficients.   One may define the cluster algebra $\mathcal{A}_{\mathrm{loc}}$ associated with a cluster pattern as the $\mathbb{Z}\mathbb{P}$-subalgebra of $\mathcal{F}$ generated by all cluster variables. In terms of our convention, this is precisely the localization $\mathcal{A}_{\rm loc}=\mathcal{A}[x_{n+1}^{-1},\ldots,x_m^{-1}]$.
\end{remark}

For the initial seed $\Sigma_{t_0}=({\bf x}, \widetilde{B}=(b_{ij})_{m\times n})$ of a cluster algebra $\mathcal{A}$, we set  
\begin{align} \label{eq:definition of y hat}
\widehat{y}_{j} = \prod_{i=1}^m x_{i}^{b_{ij}}\ \text{and}\ y_j=\prod_{i=n+1}^mx_i^{b_{ij}}, \quad 1 \le j \le n.
\end{align}  
For integer vectors $\alpha=(a_1,\dots, a_n)\in \mathbb{Z}^n$ and $\beta=(b_1,\dots, b_m)\in \mathbb{Z}^m$, we also use the notation 
\begin{align*}
    {\bf x}^\alpha:=\prod_{i=1}^nx_i^{a_i}\ \text{ and } \ {\bf x}^\beta:=\prod_{i=1}^mx_i^{b_i}.
\end{align*}

Let $M$ be a cluster monomial in the cluster variables of the seed $\mathbf{x}_t$ of $\mathcal{A}$. According to \cite[Corollary 6.3]{FZ07}, there is an integer vector $\mathbf{g}_M\in \mathbb{Z}^n$ and a polynomial $F_M(z_1,\ldots,z_n)\in \mathbb{Z}[z_1,\ldots,z_n]$ such that
\begin{align}\label{for:separation-formula}
    M={\bf x}^{\mathbf{g}_M}\frac{F_M(\hat{y}_1,\ldots, \hat{y}_n)}{F_M|_{\mathbb{P}}(y_1,\ldots, y_n)},
\end{align}
The vector $\mathbf{g}_M$ is called the {\em $\mathbf{g}$-vector} of $M$ and $F_M$ the {\em $F$-polynomial} of $M$. The formula (\ref{for:separation-formula}) is called the {\em 
 separation formula} of $M$ (cf. \cite[Section 6]{FZ07}). It is known that the $F$-polynomial $F_M(z_1,\ldots, z_n)$ has a constant of $1$ (cf. \cite{GHKK18}). Consequently, the separation formula \eqref{for:separation-formula} takes the form
 \begin{align*}
     M=\mathbf{x}^{\hat{\mathbf g}_M}F_M(\hat{y}_1,\ldots, \hat{y}_n)
 \end{align*}
 for a unique integer vector $\hat{\mathbf{g}}_M\in \ZZ^m$. Let $\operatorname{pr}_n: \ZZ^m \to \ZZ^n$ be the projection onto the first $n$ components. Then $\operatorname{pr}_n(\hat{\mathbf{g}}_M) = \mathbf{g}_M$, and the final $m-n$ entries of $\hat{\mathbf{g}}_M$ are non-negative. The vector $\hat{\mathbf{g}}_M$ is called the {\em extended $\mathbf{g}$-vector} of $M$.

\begin{lemma}\label{l:prod-cluster-monomial}
Let $t,t_1,t_2$ be vertices of $\TT_n$, let $M_1$ be a cluster monomial in $\mathbf{x}_{t_1}$ and $M_2$ a cluster monomial in $\mathbf{x}_{t_2}$. 
If the product $M_1M_2$ is a cluster monomial in $\mathbf{x}_t$, then both $M_1$ and $M_2$ are cluster monomials in $\mathbf{x}_t$.
\end{lemma}

\begin{proof}
We may take $\mathbf{x}_t:=\mathbf{x}:=(x_1,\dots, x_n)$ as the initial seed and by specializing the coefficients to $1$, we may assume that $\mathcal{A}$ is a cluster algebra with trivial coefficients. 
For $1\le j\le n$ let $\hat{y}_j$ be defined as in Equation (\ref{eq:definition of y hat}), for the cluster $\mathbf{x}$.
By the separation formula (\ref{for:separation-formula})(cf.\cite[Corollary 6.3]{FZ07}), we can write 
\[
M_1=\mathbf{x}^{\mathbf{g}_{M_1}}F_{M_1}(\hat{y}_1,\dots, \hat{y}_n) \hskip.3cm 
\mbox{ and } \hskip.3cm
M_2=\mathbf{x}^{\mathbf{g}_{M_2}}F_{M_2}(\hat{y}_1,\dots, \hat{y}_n),
\]
where $\mathbf{g}_{M_1}$ and $\mathbf{g}_{M_2}$ are the $\mathbf{g}$-vectors of $M_1$ and $M_2$  with respect to the cluster $\mathbf{x}$, while 
$F_{M_i}(z_1,\dots, z_n)$ is the $F$-polynomials of $M_i$ 
% are $F$-polynomials of $M_1$ and $M_2$ respectively, 
with non-negative integer coefficients and constant term $1$ (cf.  \cite[Theorem 1.1]{LR15}  or \cite[Theorem 0.3 (6)]{GHKK18}).
It follows that
\begin{align}\label{eq:product-F-poly}
F_{M_1}(\hat{y}_1,\dots, \hat{y}_n)F_{M_2}(\hat{y}_1,\dots, \hat{y}_n)=\mathbf{x}^{\mathbf{h}}
\end{align}
for some integer vector $\mathbf{h}\in \mathbb{Z}^n$.
Assume that at least one of $M_1$ and $M_2$, say $M_1$, is not  a cluster monomial of $\mathbf{x}$. Then $F_{M_1}$ is not a constant polynomial. It follows that the left hand side of (\ref{eq:product-F-poly}) can not be a Laurent monomial, a contradiction.
\end{proof}

% We recall the definition of upper cluster algebras, \cite{BFZ05}.
% Let $\Sigma=({\bf x}, \widetilde{B})$ be a seed of 
% $\mathcal{A}$. Then we write ${\bf x}_j = ({\bf x} \setminus \{x_j\}) \cup \{x_j'\}$ for the cluster obtained through mutation at $j$, $1\le j\le n$. 
% Denote by $\ZZ\PP[{\bf x}^{\pm }] = \ZZ\PP[x_1^{\pm 1}, \ldots, x_n^{\pm 1}]$ the ring of Laurent polynomials in the cluster variables $x_1,\dots, x_n$ with coefficients in $\ZZ\PP$. Let 
% \begin{align*}
%     \mathcal{U}(\Sigma) := \ZZ\PP[{\bf x}^{\pm}] \cap \ZZ\PP[{\bf x}^{\pm}_1] \cap \cdots \cap \ZZ\PP[{\bf x}^{\pm}_n].
% \end{align*}
% Fix an initial seed $\Sigma_0$. The {\em upper cluster algebra} $\overline{\mathcal{A}} = \overline{\mathcal{A}}(\Sigma_0) \subset \mathcal{F}$ is defined by 
% \begin{align*}
%     \overline{\mathcal{A}}(\Sigma_0) = \bigcap\limits_{\Sigma \sim \Sigma_0} \mathcal{U}(\Sigma),
% \end{align*}
% where $\Sigma \sim \Sigma_0$ means that $\Sigma$ can be obtained from $\Sigma_0$ by a sequence of mutations. 

%%
%
\subsection{Additive categorifications of cluster algebras}\label{ss:additive-cat}
Let $\mathbb{k}$ be an algebraically closed field and $\mathcal{F}$ a $\mathbb{k}$-linear  Krull--Schmidt Frobenius category which is stably $2$-Calabi-Yau, that is, the stable category $\mathcal{C}:=\underline{\mathcal{F}}$ (whose objects are the same as the objects of $\mathcal{F}$ and whose morphisms are considered up to morphisms factoring through projectives) is a $\Hom$-finite $2$-Calabi-Yau triangulated category. 
We denote the (pairwise non-isomorphic) indecomposable projective objects of $\mathcal{F}$ by 
$P_1,\dots, P_{m-n}$.

% {\color{red}We refer to \cite{NakaokaPalu19, LiuNakaoka19} for fundamentals of (Frobenius) extriangulated categories. For our purpose, it suffices to consider the following two typical examples:
% \begin{itemize}
%     \item[(1)] Stably $2$-Calabi-Yau categories, i.e., Frobenius categories whose stable categories are $2$-Calabi-Yau \cite{BIRS09};
%     \item[(2)] Higgs categories associated with ice quivers potentials \cite{Wu23,CKW23}.
% \end{itemize}
% The statements of this section are well-established for stably $2$-Calabi-Yau categories, but their proofs are valid for Frobenius extriangulated categories.
% }

Let $\Sigma$ be the suspension functor of $\mathcal{C}$. The following is well-known(cf. \cite{Happel87}).
\begin{lemma}
For any $L,N\in \mathcal{F}$, we have $\Ext^i_{\mathcal{F}}(L,N)\cong \Ext^i_{\mathcal{C}}(L,N)\cong \Hom_\mathcal{C}(L,\Sigma^i N)$ for any $i\geq 1$. 
\end{lemma}

An object $T$ in $\mathcal{F}$ (resp. in $\mathcal{C}$)  is a {\it cluster-tilting object} provided that $\Ext^1_{\mathcal{F}}(T,T)=0$ (resp. $\Ext^1_{\mathcal{C}}(T,T)=0$ ) and $\Ext^1_{\mathcal{F}}(T,X)=0$ (resp. $\Ext^1_{\mathcal{C}}(T,X)=0$) implies that $X\in \add T$. 
Two basic cluster-tilting objects are {\it mutations} of each other if they only differ in one indecomposable direct summand, we refer to \cite{BIRS09,WWZ2024} for details. We say that a rigid object $M$ is {\it reachable from $T$} if there is a basic cluster-tilting object $T'$ with $M\in \add T'$ and such that $T'$ can be obtained from $T$ by (a sequence of) mutations. 
Denote by $\text{ct}(\mathcal{F})$ the set of isomorphism classes of basic cluster-tilting objects of $\mathcal{F}$ and $\text{ct}(\mathcal{C})$ the set of isomorphism classes of basis cluster-tilting objects of $\mathcal{C}$. The {\em exchange graph} $\mathcal{E}(\mathcal{F})$ (resp. $\mathcal{E}(\mathcal{C})$) of $\mathcal{F}$ (resp. of $\mathcal{C}$) has vertex set indexed by $\text{ct}(\mathcal{F})$ (resp. $\text{ct}(\mathcal{C})$), and two basic cluster-tilting objects are connected by an edge if and only if they are linked by one mutation. 
% {\color{red}The following has been proved by \cite[Lemma II 1.3]{BIRS09} for stably $2$-Calabi-Yau categories, which can be generalized to the general setting in an obvious way.}
\begin{lemma}\cite[Lemma II 1.3]{BIRS09}\label{l:isom-exchange-graph}
There is a bijection  $\widehat{}:\operatorname{ct}(\mathcal{C})\to \operatorname{ct}(\mathcal{F})$ given by $M\mapsto M\oplus P_1\oplus\cdots\oplus P_{m-n}$, which is compatible with mutations. As a consequence, $\widehat{}$ induces an isomorphism between $\mathcal{E}(\mathcal{C})$ and $\mathcal{E}(\mathcal{F})$.
\end{lemma}

\begin{definition}\label{def:Grothendieck}
Let $\mathcal{M}$ be an additive category. \\

(1) The {\em split Grothendieck group $K_0^{sp}(\mathcal{M})$} is the abelian group generated by the (isomorphism classes of) objects of $\mathcal{M}$, quotiented out by the relations 
$[M_1\oplus M_2]-[M_1]-[M_2]=0$. 

(2) Assume that $\mathcal{M}$ is an abelian category. The {\em Grothendieck group 
$K_0(\mathcal{M})$} is the abelian group generated by the (isomorphism classes of) objects of $\mathcal{M}$, quotiented out by the relations 
$[A]+[B]=[C]$ for every short exact sequence $0\to A\to C\to B\to 0$ in $\mathcal{M}$. 
If $\mathcal{M}$ also has a 
monoidal structure, the Grothendieck $K_0(\mathcal{M})$ 
has a ring structure given by $[A]\cdot[B]=[A\otimes B]$.

%for every short exact sequence $0\to A\to C\to B\to 0$ in $\mathcal{M}$.
\end{definition}

Let $T$ be a basic cluster-tilting object of $\mathcal{C}$ and  $\widehat{T}$ the associated basic cluster-tilting object of $\mathcal{F}$.  
%Denote by $K_0^{sp}(\add \widehat{T})$ the split Grothendieck group of the category $\add \widehat{T}$. 
%
% (((, {\it i.e.} the quotient of the free abelian group on the set of isomorphism class $[N]$ of $N\in \add \widehat{T}$, modulo the subgroup generated by all elements $[N_1\oplus N_2]-[N_1]-[N_2]$. )))
For any object $M\in \mathcal{F}$, there is a short exact sequence $0\to L_1\to L_0\to M\to 0$, where $L_0,L_1\in \add \widehat{T}$. 
We call 
\[
\operatorname{ind}_{\widehat{T}}(M)=[L_0]-[L_1]\in K_0^{sp}(\add \widehat{T})
\]
the {\em index} of $M$ with respect to $\widehat{T}$.
Note that if $T=\bigoplus_{i=1}^nT_i$ with indecomposable direct summands $T_1,\dots, T_n$, then we may identify $K_0^{sp}(\add \widehat{T})$ with $\mathbb{Z}^{m}$ by sending $[T_i]$ to $e_i$ for $1\leq i\leq n$ and $[P_j]$ to $e_{j+n}$ for $1\leq j\leq m-n$, where $e_1,\dots, e_{m}$ is the standard basis of $\mathbb{Z}^{m}$.

The quiver $Q_{\widehat{T}}$ of $\widehat{T}$ has vertex set indexed by the indecomposable direct summands of $\widehat{T}$ and the number of arrows $T_i\to T_j$ in $Q_{\widehat{T}}$ between two indecomposable direct summands $T_i$ and $T_j$ is given by the dimension of the space of irreducible maps $\rad(T_i,T_j)/\rad^2(T_i,T_j)$, where $\rad(-,-)$ is the radical of the subcategory $\add \widehat{T}$. We denote by $Q_T$ the quiver obtained from $Q_{\widehat{T}}$ by deleting  all arrows between vertices corresponding to projective-injectives. The quiver $Q_T$ is called the {\em extended quiver of $T$}. We always assume that the extended quiver of any basic cluster-tilting object is finite.
The category $\mathcal{F}$ is said to have a {\em cluster structure}  if the extended quiver of every basic cluster-tilting object of $\mathcal{C}$ has no loops nor $2$-cycles (cf. \cite[Section II.1]{BIRS09}).

Now assume that $\mathcal{F}$ has a cluster structure and that $T$ is a basic cluster-tilting object of $\mathcal{C}$ such that $\End_{\mathcal{F}}(\widehat{T})$ is Noetherian.
By regarding the vertices in $Q_T$ which correspond to projective-injectives of $\mathcal{F}$ as frozen, we obtain a cluster algebra $\mathcal{A}(Q_T)$ as in Section~\ref{subsec:cluster algebras}. It follows that $(\mathcal{F},\widehat{T})$ is an additive categorification of $\mathcal{A}(Q_T)$ (cf. \cite[Defintion 5.1]{FK10}). 
We write $T=\bigoplus_{i=1}^{n}T_i$ and 
$\widehat{T}=\bigoplus_{i=1}^{m}T_i$, where 
$T_{j+n}:=P_j$ for $1\leq j\leq m-n$. Let $B_{\widehat{T}}$ be the associated  integral $m\times m$-matrix, i.e. 
$B_{\widehat{T}}:=(b_{ij})$ 
where $b_{ij}=\#\{T_i\to T_j\}-\#\{T_j\to T_i\}$. 
Denote by $B_T$ the submatrix of $B_{\widehat{T}}$ formed by the first $n$ columns, i.e. the $m\times n$ matrix corresponding to $Q_T$ as in Remark~\ref{rem:Q-gives-B}.

For any $M\in \mathcal{F}$, $\Ext_\mathcal{F}^1(\widehat{T}, M)\cong \Ext_\mathcal{F}^1(T, M) $ is a right $\End_\mathcal{C}(T)$-module and we define the Caldero--Chapoton \cite{FK10}  map $\mathbf{X}_M^T$ of $M$  as follows 

\[
\mathbf{X}_M^{\widehat{T}}:=\mathbf{x}^{\operatorname{ind}_{\widehat{T}}(M)}\sum_{e}\chi(\text{Gr}_e(\Ext^1_\mathcal{F}(T,M))\mathbf{x}^{B_T e}\in \mathbb{Z}[x_1^{\pm},\dots, x_{m}^{\pm}], 
\]
where 
\begin{itemize}
% \item for an integer vector $\alpha=(a_1,\dots, a_{m+n})$, we denote by $\mathbf{x}^\alpha=x_1^{a_1}\cdots x_{m+n}^{a_{m+n}}$;
\item $\text{Gr}_e(\Ext^1_\mathcal{F}(T,M))$ is the quiver Grassmanian of $\Ext^1_\mathcal{F}(T,M)$ consisting of sub $\End_\mathcal{C}(T)$-modules with dimension vector $e$ and where $e$ is understood as a column vector;
\item $\chi(\text{Gr}_e(\Ext^1_\mathcal{F}(T,M))$ is the Euler characteristic of $\text{Gr}_e(\Ext^1_\mathcal{F}(T,M))$.
\end{itemize}  
% We remark that the above formula for $\mathbf{X}_M^{\hat T}$ is equivalent to the ones in \cite[Definition 4.2]{WWZ2024} and \cite[Definition 5.3]{FMP23} by the arguments in the proof of \cite[Proposition 5.4]{FMP23}.
The following result justifies the name of additive categorification, which follows the proof of \cite[Theorem 5.4]{FK10}.

\begin{theorem}\label{t:additve-cat}
Keep the assumptions and notation as above. 
\begin{itemize}
\item[(1)] The map $M\mapsto \mathbf{X}_M^{\hat T}$ induces a bijection from the set of isomorphism classes of indecomposable rigid  objects of $\mathcal{F}$ reachable from $T$ onto
the set of cluster variables of $\mathcal{A}(Q_T)$. Under this bijection, the cluster-tilting
objects which are reachable from $T$ correspond to the clusters of $\mathcal{A}(Q_T)$.
\item[(2)] {If $M$ is a rigid object reachable from $T$, then $\operatorname{ind}_{\widehat{T}}(M)$  is the extended $\mathbf{g}$-vector of the cluster monomial $\mathbf{X}_M^{\hat T}$.}
\item[(3)] If $M$ is an indecomposable rigid nonprojective object reachable from $T$, then the truncation of $\operatorname{ind}_{\widehat{T}}(M)$ with respect to the first $n$ coordinates coincides with the $\mathbf{g}$-vector of the cluster variable $\mathbf{X}_M^{\hat T}$.
    \end{itemize}
\end{theorem}
As a direct consequence of  Lemma \ref{l:prod-cluster-monomial} and Theorem \ref{t:additve-cat}, we obtain the categorical counterpart of Lemma \ref{l:prod-cluster-monomial}.
\begin{corollary}\label{c:prod-cluster-monomial}
    Let $L,M,N$ be rigid objects of $\mathcal{C}$ which are reachable from $T$. If $\mathbf{X}_L^{\hat T}\mathbf{X}_N^{\hat T}=\mathbf{X}_M^{\hat T}$, then $L, N\in \add M$.
\end{corollary}
It is clear that the connectivity of $\mathcal{E}(\mathcal{C})$ is equivalent to the property that every rigid object is reachable from $T$. 
The following folklore conjecture is known as the {\it additive reachability conjecture} (cf. \cite[Remark 5.9]{Qin20}).
\begin{conjecture}\label{c:addtive-reach-conj}
Keep the assumptions and  notation as above.
If there is a path in $\mathcal{E}(\mathcal{C})$ connecting $T$ and $\Sigma T$, then $\mathcal{E}(\mathcal{C})$ has exactly one connected component.
\end{conjecture}

\begin{remark}
The existence of a path linking $T$ and $\Sigma T$ in $\mathcal{E}(\mathcal{C})$ is equivalent to the existence of a green-to-red sequence \cite{KelD20} for the Gabriel quiver of $\End_{\mathcal{C}}(T)$ (equivalently, for the principal part of $B_T$), which in turn is equivalent to the property that the associated cluster algebra is injective-reachable in the sense of Qin \cite{Qin17}.
It was conjectured by \cite[Conjecture 2]{Mi18} that if a skew-symmetrizable integer matrix $B$ admits a green-to-red sequence, then the associated cluster algebra $\mathcal{A}(B)$ coincides with its upper cluster algebra $\overline{\mathcal{A}}(B)$ under a particular choice of coefficients.
\end{remark}

\begin{remark}\label{rem:existence-additive-cat}
All results in this section can be generalized to the setting where $\mathcal{F}$ is a $\mathbb{k}$-linear  Krull--Schmidt Frobenius extriangulated category which is stably $2$-Calabi-Yau. For example, \cite{WWZ2024,FMP23} established the Caldero--Chapoton map and its multiplication formulas in this general setting. We refer to \cite{NakaokaPalu19, LiuNakaoka19} for fundamentals of (Frobenius) extriangulated categories. 

Every injective-reachable cluster algebra of skew-symmetric type admits an additive categorification by Frobenius extriangulated category whose stable category is $2$-Calabi-Yau. Indeed, let $\mathcal{A}=\mathcal{A}(Q)$ be an injective-reachable cluster algebra associated to a finite quiver $Q$ as in Section \ref{subsec:cluster algebras}. Denote by $Q^{mu}$ the full subquiver of $Q$ consisting of mutable vertices. Let $W$ be a non-degenerate potential on $Q$ and $W^{mu}$ the restriction of $W$ on $Q^{mu}$, cf. \cite{DWZ08}. Since $\mathcal{A}$ is injective-reachable, we know that the Jacobian algebra $J(Q^{mu},W^{mu})$ is finite dimensional, see \cite[Proposition 8.1]{BDM14}. According to \cite[Theorems 1.1 \& 1.2]{CKW23}, the Higgs category $\mathcal{H}$ associated with $(Q,W)$ is a Frobenius extriangulated category which is stably $2$-Calabi-Yau, which provides an additive categorification of $\mathcal{A}$.
\end{remark}

\subsection{Quantum affine algebras} \label{subsec:qu_aff_algebras-qcharacters}

Let $\mathfrak{g}$ be a simple finite-dimensional Lie algebra and $I$ the set of vertices of the Dynkin diagram of $\mathfrak{g}$. Denote by $\{ \omega_i: i \in I \}$, $\{\alpha_i : i \in I\}$, $\{\alpha_i^{\vee} : i \in I\}$ the set of fundamental weights, the set of simple roots, the set of simple coroots, respectively. Denote by $P$ the integral weight lattice and
$P^+$ the set of dominant weights. The Cartan matrix is $C = (\alpha_j(\alpha_i^{\vee}))_{i,j \in I}$. Let $D = \diag(d_i: i \in I)$, where $d_i$'s are minimal positive integers such that $DC$ is symmetric.  

The quantum affine algebra $U_q(\widehat{\mathfrak{g}})$ is a Hopf algebra that is a $q$-deformation of the universal enveloping algebra of the affine Lie algebra $\widehat{\mathfrak{g}}$ of $\mathfrak{g}$ \cite{Dri85, Jim85}. In this paper, we take $q$ to be a non-zero complex number which is not a root of unity. 

Denote by $\mathcal{P}$ the free abelian group generated by formal variables $Y_{i, a}^{\pm 1}$, $i \in I$, $a \in \CC^*$ {and} denote by $\mathcal{P}^+$ the submonoid of $\mathcal{P}$ generated by $Y_{i, a}$, $i \in I$, $a \in \CC^*$. 
The elements in $\mathcal{P}^+$ are called {\it dominant monomials}. 
Chari and Pressley \cite[\S 3]{CP95a} proved that any finite-dimensional simple object in the category $\mathscr{C}$ of finite dimensional $U_q(\widehat{\mathfrak{g}})$-modules is a highest $l$-weight module (denoted by $L(M)$) with a highest $l$-weight $M \in \mathcal{P}^+$. 
The so-called $l$-weights (or loop weights) are generalisations of weights to the setting of quantum affine algebras. The $l$-weights can be identified with elements in $\mathcal{P}$ (see below in this section for more details). 

Frenkel and Reshetikhin \cite{FR98} introduced the theory of $q$-characters which is a powerful tool for studying representations of quantum affine algebras. 
The $q$-character map is an injective ring morphism $\chi_q$ from the Grothendieck ring of $\mathscr{C}$ to $\mathbb{Z}\mathcal{P} = \mathbb{Z}[Y_{i, a}^{\pm 1}]_{i\in I, a\in \mathbb{C}^*}$. 
For a $U_q(\widehat{\mathfrak{g}})$-module $V$, $\chi_q(V)$ encodes the decomposition of $V$ into common generalized eigenspaces for the action of a large commutative subalgebra of $U_q(\widehat{\mathfrak{g}})$ (the loop-Cartan subalgebra). These generalized eigenspaces are called {\it $l$-weight spaces} and generalized eigenvalues are called {\it $l$-weights}. One can identify $l$-weights with monomials in $\mathcal{P}$ \cite{FR98}. Then
the $q$-character of a $U_q(\widehat{\mathfrak{g}})$-module $V$ is given by (see \cite{FR98})
\begin{align*}
\chi_q(V) = \sum_{  m \in \mathcal{P}} \dim(V_{m}) m \in \mathbb{Z}\mathcal{P},
\end{align*}
where $V_{m}$ is the $l$-weight space with $l$-weight $m$.

\subsection{Monoidal categorifications of cluster algebras}\label{ss:monoidal-cat}

Hernandez and Leclerc \cite{HL10} introduced the concept of monoidal categorification of cluster algebras. Let $(\mathscr{M}, \otimes)$ be a monoidal category. A simple object $S$ of $\mathscr{M}$ is called {\it real} if $S \otimes S$ is simple \cite{Lec03}. A simple object $S$ is called {\it prime} if there exists no non-trivial factorization $S \cong S' \otimes S''$ \cite{CP97}. An abelian monoidal category $\mathscr{M}$ is said to be a {\it monoidal categorification} of the cluster algebra $\mathcal{A}$ if the Grothendieck ring 
$K_0(\mathscr{M})$ of $\mathscr{M}$ is isomorphic to $\mathcal{A}$ and if  
\begin{enumerate}[(i)]
% \item the cluster monomials of $\mathcal{A}$ are the classes of all real simple objects of $\mathscr{M}$,
% \item the cluster variables of $\mathcal{A}$ (including the frozen variables) are the classes of all real prime simple objects of $\mathscr{M}$. 
\item any cluster monomial of $\mathcal{A}$ corresponds to the class of a real simple object of $\mathscr{M}$,
\item any cluster variable of $\mathcal{A}$ corresponds to the class of a real simple prime object of $\mathscr{M}$.
\end{enumerate}
The above definition is weaker than the original definition of monoidal categorification in \cite[Definition 2.1]{HL10}, see Section 5 of \cite{kashiwara2018crystal}.

Recall that $\mathscr{C}$ is the category of finite-dimensional $U_q(\widehat{\mathfrak{g}})$-modules (Section~\ref{subsec:qu_aff_algebras-qcharacters}). 
In \cite{HL10},  Hernandez and Leclerc introduced a full subcategory $\mathscr{C}_{\ell} := \mathscr{C}_{\ell}^{\mathfrak{g}}$ of $\mathscr{C}$ for every $\ell \in \mathbb{Z}_{\geq 0}$ (cf. also \cite{HL16}). We recall this now. 
We fix $a \in \CC^*$ and write 
$Y_{i,s}$ for $Y_{i,aq^s}$, $i \in I$, $s \in \ZZ$.

\begin{definition}[{\cite[Section 2.2]{HL15} and \cite[Definition 4.1]{FHOO23}}] \label{def:height} 
Let $I$ be the set of vertices {of a} simply-laced Dynkin diagram. Fix an orientation of the Dynkin diagram. A {\it height function} for this orientation is a function $\xi:I\to \ZZ$ such that whenever there is an arrow $i\to j$ in the oriented Dynkin diagram, we have 
$\xi(j)-\xi(i)=1$. 
\end{definition} 
See Remark~\ref{rem:height-functions} below for two examples of height functions for different 
orientations of a type $\mathbb{A}$ Dynkin diagram. 

For $\ell \in \ZZ_{\ge 0}$, denote by $\mathcal{P}_\ell$ the subgroup of $\mathcal{P}$ generated by $Y_{i,\xi(i)-2r}^{\pm 1}$, $i \in I$, $r \in [0, d \ell]$, where $d$ is the maximal diagonal element in the diagonal matrix $D$ from Section~\ref{subsec:qu_aff_algebras-qcharacters}, and $\xi: I \to \ZZ$ is 
a height function. 

Denote by $\mathcal{P}^+_\ell$ the submonoid of $\mathcal{P}^+$ generated by $Y_{i,\xi(i)-2r}$, $i \in I$, $r \in [0, d \ell]$. 
The subcategory $\mathscr{C}_\ell=\mathscr{C}_\ell^{\mathfrak{g}}$ has as objects the finite-dimensional $U_q(\widehat{\mathfrak{g}})$-modules $V$ which satisfy the following condition: for every composition factor $S$ of $V$, the highest $l$-weight of $S$ is a monomial in $\mathcal{P}^+_\ell$, cf. \cite{HL10}. The simple modules in $\mathscr{C}_{\ell}$ are of the form $L(M)$ (see \cite{CP94}, \cite{HL10}), where $M \in \mathcal{P}_{\ell}^+$.

Denote by $K_0(\mathscr{C}_{\ell})$ the Grothendieck ring of $\mathscr{C}_{\ell}$. By a slight abuse of notation, we 
sometimes write $L(M)$ or $[M]$ for $[L(M)]$ ($M \in \mathcal{P}^+$) in $K_0(\mathscr{C}_{\ell})$. Hernandez and Leclerc proved that there is a cluster algebra structure on $K_0(\mathscr{C}_{\ell})$ \cite{HL10}. They conjectured that the set of cluster monomials in $K_0(\mathscr{C}_{\ell})$ is the set of equivalence classes of real prime simple modules in $K_0(\mathscr{C}_{\ell})$. The implication ``all cluster variables (resp. cluster monomials) are real prime modules (resp. real modules)'' of their conjecture is proved by Qin in \cite{Qin17} for ADE types, and by Kashiwara, Kim, Oh, and Park for general types in \cite{KKOP20, KKOP22, KKOP21b}. The other direction, also known as the {\it monoidal reachability conjecture}, is widely open (see \cite[Conjecture 5.7]{HL21} and \cite[Remark 5.9]{Qin20}).

\begin{conjecture}[{\cite{HL10, HL21}}] \label{c:monoidal-reach-conj} 
For every $\ell \in \ZZ_{\ge 1}$, the $q$-character of every real prime module (resp. real module) in $\mathscr{C}_{\ell}$ is a cluster variable (resp. cluster monomial) in $K_0(\mathscr{C}_{\ell})$.
\end{conjecture}

\begin{remark}\label{rem:height-functions}
In this paper, we work with two different orientations of the Dynkin diagram of type $\mathbb{A}_{k-1}$: 

(1) If the Dynkin diagram is linearly oriented with arrows $i\to i+1$ (for $i< k-1$), we choose the height function $\xi$: 
\[
\xi:I\to \ZZ, \quad 
\xi(i) = i-2. 
\]

(2) If the orientation of the Dynkin is 
bipartite such that $1$ is a source and $2$ a sink 
($1\to 2 \leftarrow 3 \to\cdots$), we choose the  height function 
$\xi':I\to \ZZ$, 
\[
\xi'(i) = 
\left\{ 
\begin{array}{ll}
0, & \mbox{if $i$ is even,} \\
-1, & \mbox{if $i$ is odd.}
\end{array}
\right.
\]
\end{remark}

To finish this section, we give {the} quiver of an initial seed of $K_0(\mathscr{C}_{\ell})$ for $\mathfrak{g}=\mathfrak{sl}_{k}$ \cite{HL10}. 

\begin{example}\label{ex:quiver-Q-l}
For $k\ge 2$, let $I=\{1,\dots, k-1\}$ {be} the vertices of the Dynkin diagram of type $\mathbb{A}_{k-1}$ and let $\ell$ be {a non-negative integer}. In type $\mathbb{A}$, the maximal element $d$ (of the symmetrizing matrix from Section~\ref{subsec:qu_aff_algebras-qcharacters}) is equal to $1$. 
We orient this Dynkin diagram linearly, $i\to i+1$ for $i<k-1$ and use the height function $\xi$ from Remark~\ref{rem:height-functions}.  
Then the submonoid $\mathcal{P}_{\ell}^+$ is generated by $Y_{i,\xi(i)-2r}$ for $i\in I$ and $r\in[0,\ell]$. With this, we define a quiver $Q_\ell=Q_{\ell,\xi}$ for $\mathcal{P}_{\ell}^+$ 
as follows. It has $(k-1)(\ell+1)$ vertices 
\[
\begin{array}{cl}
&  \{(i,-2),(i,-4),\dots, (i,-2(\ell+1))\mid 1\le i\le k-1, \mbox{$i$ odd} \}  \\
\cup & \{(i,-1),(i,-3),\dots, (i,-2(\ell+1)+1)\mid 1\le i\le k-1, \mbox{$i$ even }\}. 
\end{array}
\]
And the arrows are of 
the form 
\[
\begin{array}{ll}
(i,a)\to (i,a-2), &\mbox{for $1\le i\le k-1$},\\
(i,a)\to (i+1, a+(-1)^{i+1}), &\mbox{for $1\le i< k-1$},\\
(i,a)\to (i-1,a+2+(-1)^{i-1}), &\mbox{for $1< i\le k-1$},
\end{array}
\]
with $a\in \{-2,-4,\dots, -2\ell\}$ for odd $i$ and $a\in \{-1,-3,\dots, -2\ell+1\}$ for even $i$.

The quiver $Q_{\ell}$ gives an initial seed for 
$K_0(\mathscr{C}_{\ell})$ for 
$\mathfrak{g}=\mathfrak{sl}_k$ and it has appeared in~\cite{HL10}. 
An example is on the right hand side 
of~Figure~\ref{fig:initial HL-quiver and initial quiver of Grassmannian cluster algebra} which shows the vertices $(i,m)$ as well as 
%we draw the vertices $(i,m)$ and to the left of them 
the monomials (cluster variables) they correspond to, indicating frozen variables with boxes. 
\end{example} 

\subsection{Grassmannian cluster algebras} \label{subsec:Grassmannian cluster algebras}

%For $k \le n$, the Grassmannian $\Gr(k,n)$ is the set of $k$-dimensional subspaces in an $n$-dimensional vector space. 
Let $k\le n$. 
In this paper, we denote by $\Gr(k,n)$ (the affine cone over) the Grassmannian of $k$-dimensional subspaces in $\CC^n$, and by $\CC[\Gr(k,n)]$ its coordinate ring. %We will always assume $k\le \frac{n}{2}$ as $\Gr(k,n)$ and $\Gr(n-k,n)$ are dual to each other. 
The coordinate ring $\CC[\Gr(k,n)]$ is generated by the Pl\"{u}cker coordinates 
\begin{align*}
P_{i_1, \ldots, i_{k}}, \quad 1 \leq i_1 < \cdots < i_{k} \leq n, 
\end{align*}
subject to the so-called Pl\"ucker relations, see e.g.~\cite[Chapter 9]{Marsh14} for more details. 

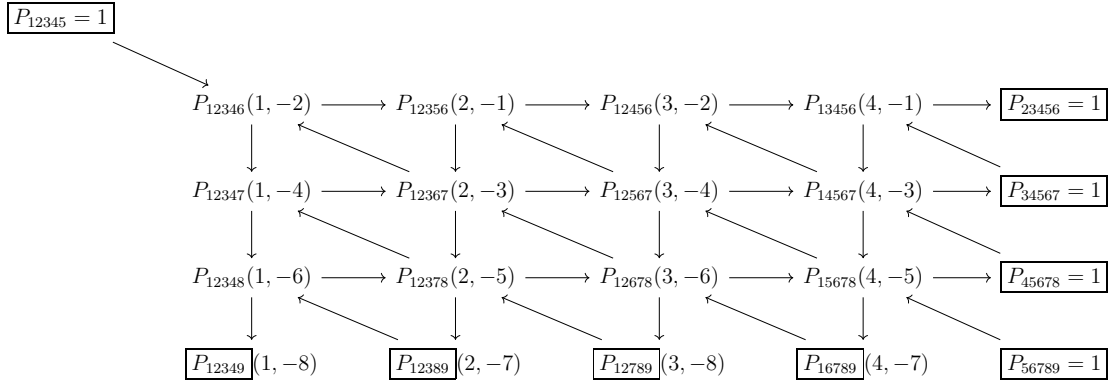
\begin{figure}
    \centering 
\adjustbox{scale=0.7}{\begin{tikzcd}[arrows=<-]
	{\fbox{$P_{12345}=1$}}  \\
	& {P_{12346} (1,-2)} & {P_{12356} (2,-1)} & {P_{12456} (3,-2)} & {P_{13456} (4,-1)} & {\fbox{$P_{23456}=1$}} \\
	& {P_{12347} (1,-4)} & {P_{12367} (2,-3)} & {P_{12567} (3,-4)} & {P_{14567} (4,-3)} & {\fbox{$P_{34567}=1$}} \\
	& {P_{12348} (1,-6)} & {P_{12378} (2,-5)} & {P_{12678} (3,-6)} & {P_{15678} (4,-5)} & {\fbox{$P_{45678}=1$}} \\
	& {\fbox{$P_{12349}$} (1,-8)} & {\fbox{$P_{12389}$} (2,-7)} & {\fbox{$P_{12789}$} (3,-8)} & {\fbox{$P_{16789}$} (4,-7)} & {\fbox{$P_{56789}=1$}}
	\arrow[from=2-2, to=1-1]
	\arrow[from=2-3, to=2-2]
	\arrow[from=2-4, to=2-3]
	\arrow[from=3-5, to=2-5]
	\arrow[from=4-5, to=3-5]
	\arrow[from=5-5, to=4-5]
	\arrow[from=3-2, to=2-2]
	\arrow[from=4-2, to=3-2]
	\arrow[from=5-2, to=4-2]
	\arrow[from=5-3, to=4-3]
	\arrow[from=4-3, to=3-3]
	\arrow[from=3-3, to=2-3]
	\arrow[from=5-4, to=4-4]
	\arrow[from=4-4, to=3-4]
	\arrow[from=3-4, to=2-4]
	\arrow[from=3-4, to=3-3]
	\arrow[from=3-3, to=3-2]
	\arrow[from=3-6, to=3-5]
	\arrow[from=4-4, to=4-3]
	\arrow[from=4-3, to=4-2]
	\arrow[from=4-6, to=4-5]
	\arrow[from=2-5, to=3-6]
	\arrow[from=3-5, to=4-6]
	\arrow[from=4-5, to=5-6]
	\arrow[from=2-2, to=3-3]
	\arrow[from=3-2, to=4-3]
	\arrow[from=4-2, to=5-3]
	\arrow[from=2-3, to=3-4]
	\arrow[from=3-3, to=4-4]
	\arrow[from=4-3, to=5-4]
	\arrow[from=2-5, to=2-4]
	\arrow[from=4-5, to=4-4]
	\arrow[from=3-5, to=3-4]
	\arrow[from=2-4, to=3-5]
	\arrow[from=3-4, to=4-5]
	\arrow[from=4-4, to=5-5]
	\arrow[from=2-6, to=2-5]
\end{tikzcd}
}
% \adjustbox{scale=0.6}{
% \begin{tikzcd}[arrows=<-]
% 	{\fbox{$P_{12345}=1$}}  \\
% 	{P_{12346} (1,-2)} & {P_{12356} (2,-1)} & {P_{12456} (3,-2)} & {P_{13456} (4,-1)} & {\fbox{$P_{23456}=1$}} \\
% 	{P_{12347} (1,-4)} & {P_{12367} (2,-3)} & {P_{12567} (3,-4)} & {P_{14567} (4,-3)} & {\fbox{$P_{34567}=1$}} \\
% 	{P_{12348} (1,-6)} & {P_{12378} (2,-5)} & {P_{12678} (3,-6)} & {P_{15678} (4,-5)} & {\fbox{$P_{45678}=1$}} \\
% 	{\fbox{$P_{12349}$} (1,-8)} & {\fbox{$P_{12389}$} (2,-7)} & {\fbox{$P_{12789}$} (3,-8)} & {\fbox{$P_{16789}$} (4,-7)} & {\fbox{$P_{56789}=1$}}
% 	\arrow[from=2-2, to=2-1]
% 	\arrow[from=2-3, to=2-2]
% 	\arrow[from=2-4, to=2-3]
% 	\arrow[from=3-1, to=2-1]
% 	\arrow[from=4-1, to=3-1]
% 	\arrow[from=5-1, to=4-1]
% 	\arrow[from=3-2, to=2-2]
% 	\arrow[from=4-2, to=3-2]
% 	\arrow[from=5-2, to=4-2]
% 	\arrow[from=5-3, to=4-3]
% 	\arrow[from=4-3, to=3-3]
% 	\arrow[from=3-3, to=2-3]
% 	\arrow[from=5-4, to=4-4]
% 	\arrow[from=4-4, to=3-4]
% 	\arrow[from=3-4, to=2-4]
% 	\arrow[from=3-4, to=3-3]
% 	\arrow[from=3-3, to=3-2]
% 	\arrow[from=3-2, to=3-1]
% 	\arrow[from=4-4, to=4-3]
% 	\arrow[from=4-3, to=4-2]
% 	\arrow[from=4-2, to=4-1]
% 	\arrow[from=2-1, to=3-2]
% 	\arrow[from=3-1, to=4-2]
% 	\arrow[from=4-1, to=5-2]
% 	\arrow[from=2-2, to=3-3]
% 	\arrow[from=3-2, to=4-3]
% 	\arrow[from=4-2, to=5-3]
% 	\arrow[from=2-3, to=3-4]
% 	\arrow[from=3-3, to=4-4]
% 	\arrow[from=4-3, to=5-4]
% 	\arrow[from=2-5, to=2-4]
% 	\arrow[from=4-5, to=4-4]
% 	\arrow[from=3-5, to=3-4]
% 	\arrow[from=2-4, to=3-5]
% 	\arrow[from=3-4, to=4-5]
% 	\arrow[from=4-4, to=5-5]
% 	\arrow[from=2-1, to=1-1]
% \end{tikzcd}
% } 
\caption{The quiver of the initial seed for the Grassmannian cluster algebra $\CC[\Gr(5,9)]$. The coordinates of the vertices are {those of the vertices} of $Q_{3}$.}\label{fig:5-9-quiver}
%It corresponds to the second quiver of Figure~\ref{fig:initial HL-quiver and initial quiver of Grassmannian cluster algebra}. 
\end{figure}

It was shown by Scott \cite{Sco} that $\CC[\Gr(k,n)]$ has a cluster algebra structure. The cluster algebra $\CC[\Gr(k,n)]$ has an initial seed with the initial quiver $Q$ with vertices 
\[
\{(0,0)\} \cup \{(a,b) : a\in [n-k], \ b\in  [k] 
\}
\]
and arrows 
\[
\begin{array}{rl}
(0,0)\to (1,1), \\
(a-1,b)\to (a,b) & \mbox{for } a\in [2,n-k],\ b\in [k], \\
(a,b-1)\to (a,b) & \mbox{for } a\in [n-k],\ b\in [2,k], \\
(a+1,b+1)\to (a,b) & \mbox{for } a\in [n-k-1],\ b\in [k-1]. 
\end{array}
\]
The quiver can be identified with the quiver $Q_{n-k-1}$ 
from Section~\ref{ss:monoidal-cat} by removing the vertices $(0,0)$, $(1,k),\dots, (n-k,k)$ and identifying the coordinate $(a,b)$ with $(b, -2a)$ if $b$ is odd and with $(b,-2a+1)$ if $b$ is even. 
The quiver $Q$ is formed by oriented triangles. 
%$(0, 0)$, and $(a,b)$, $a \in [n-k]$, $b \in [k]$, and arrows $(0,0) \to (1,1)$, $(a-1,b) \to (a,b)$, $a \in [2, n-k]$, $b \in [k]$, $(a, b-1) \to (a,b)$, $a \in [n-k]$, $b \in [2, k]$, $(a+1,b+1) \to (a,b)$, $a \in [n-k-1]$, $b \in [k-1]$, 
% See 
% %Figure \ref{fig:initial-cluster-for-Gr(3,7)} for $k=3$, $n=7$ or 
% Sections~\ref{subsec:non-rigid-module-example-Gr39} and~\ref{subsec:non-rigid-module-example-Gr48}
% for (the stable parts of) these quivers in the cases $(3,9)$ and $(4,8)$. 
All the initial cluster variables are Pl\"{u}cker coordinates. The frozen variable at $(0,0)$ is $P_{1,\ldots,k}$. The cluster variables (including frozen variables) in the column with $b=1$ are $P_{1,2,\ldots, k-1, k+1}$, $\ldots$, $P_{1,2,\ldots, k-1, n}$. The cluster variables (including frozen variables) in the column with $b=2$ are $P_{1,2,\ldots, k-2, k, k+1}$, $\ldots$, $P_{1,2,\ldots, k-2, n-1, n}$. The cluster variables (including frozen variables) in column $b=k$ are $P_{2,\ldots, k+1}$, $\ldots$, $P_{n-k+1, \ldots, n}$, see for example Figure \ref{fig:5-9-quiver}.

% \begin{figure}
% \begin{center}
% \adjustbox{scale=0.72}{
% \begin{tikzcd}[arrows=<-]
% {\fbox{$P_{123} \ (0,0)$}} \\
%  & {P_{124} \ (1,1)} & {P_{134} \ (1,2)} & {\fbox{$P_{234} \ (1,3)$}} \\
%  & {P_{125} \ (2,1)} & {P_{145} \ (2,2)} & {\fbox{$P_{345} \ (2,3)$}}\\
% & {P_{126} \ (3,1)} & {P_{156} \ (3,2)} & {\fbox{$P_{456} \ (3,3)$}} \\
% & {\fbox{$P_{127} \ (4,1)$}} & {\fbox{$P_{167} \ (4,2)$}} & {\fbox{$P_{567} \ (4,3)$}}
% \arrow[from=2-2, to=1-1]
% \arrow[from=3-2, to=2-2]
% \arrow[from=2-3, to=2-2]
% \arrow[from=2-4, to=2-3]
% \arrow[from=3-4, to=2-4]
% \arrow[from=4-4, to=3-4]
% \arrow[from=3-4, to=3-3]
% \arrow[from=3-3, to=3-2]
% \arrow[from=3-3, to=2-3]
% \arrow[from=4-3, to=3-3]
% \arrow[from=4-4, to=4-3]
% \arrow[from=4-3, to=4-2]
% \arrow[from=4-2, to=3-2]
% \arrow[from=2-2, to=3-3]
% \arrow[from=2-3, to=3-4]
% \arrow[from=3-3, to=4-3]
% \arrow[from=3-3, to=4-4]
% \arrow[from=5-2, to=4-2]
% \arrow[from=5-3, to=4-3]
% \arrow[from=3-2, to=4-3]
% \arrow[from=4-2, to=5-3]
% \arrow[from=5-4, to=4-4]
% \arrow[from=4-3, to=5-4]
% \arrow[from=5-4, to=5-3]
% \arrow[from=5-3, to=5-2]
% \end{tikzcd}  
% }
% \end{center}
% \caption{The initial cluster for $\CC[\Gr(3,7)]$, where $(a,b)$'s at the vertices are used to denote the positions of cluster variables and frozen variables. \textcolor{red}{change the coordinates of vertices}} 
% \label{fig:initial-cluster-for-Gr(3,7)}
% %{fig:initial cluster for a quotient of Gr(3,7)}
% \end{figure}

We will use the notation 
$\CC[\Gr(k,n,\sim)]$ to denote the quotient of $\CC[\Gr(k,n)]$ by the ideal generated by  
\[
P_{i,i+1, \ldots, k+i-1}-1, \quad i \in [n-k+1].
\] 

%%%%
%
\subsection{Isomorphism between $K_0(\mathscr{C}_{\ell}^{\mathfrak{sl}_k})$ and $\CC[\Gr(k,n,\sim)]$} \label{subsec:correspondence between modules and tableaux}

In the case of $\mathfrak{g}=\mathfrak{sl}_k$, we choose the 
height function 
$\xi(i)=i-2$, $i \in I = [k-1]$ for the linear orientation, see Remark~\ref{rem:height-functions}. 
The Grothendieck ring  $K_0(\mathscr{C}_{\ell}^{\mathfrak{sl}_k})$ is isomorphic to the cluster algebra $\CC[\Gr(k,n,\sim)]$, for $n=k+\ell+1$,~\cite[Section 13.9]{HL10}. 
The quiver $Q_{\ell}$ of the initial seed in Section~\ref{ss:monoidal-cat} (Example~\ref{ex:quiver-Q-l}) is the same as the quiver  of the initial seed of $\CC[\Gr(k,n,\sim)]$ described in Section~\ref{subsec:Grassmannian cluster algebras}. 
This appears on the right hand side of Figure~\ref{fig:initial HL-quiver and initial quiver of Grassmannian cluster algebra} for $n=k+\ell+1=9$.

We recall a few definitions related to the theory of Young tableaux. A {\em Young tableau} is a diagram of finitely many square boxes arranged in rows which are top-left adjusted and where each box is filled with positive integers. In this paper, all such diagrams are rectangular (i.e. each row has the same number of boxes). A {\em semistandard Young tableau} is a Young tableau where the entries are weakly increasing in each row and strictly increasing in each column. For $k,n \in \ZZ_{\ge 1}$, we denote by ${\rm SSYT}(k, [n])$ the set of rectangular semistandard Young tableaux with $k$ rows and with entries in $[n] = \{1,\ldots, n\}$ (with arbitrarily many columns). For $\mathbf{S},\bT \in {\rm SSYT}(k, [n])$, let $\mathbf{S} \cup \bT$ be the row-increasing tableau whose $i$th row is the union of the $i$th rows of $\mathbf{S}$ and $\bT$ (as multisets), for any $i$, \cite{CDFL20}. 
By Lemma 3.6 in \cite{CDFL20}, $\mathbf{S} \cup \bT$ is in ${\rm SSYT}(k, [n])$. We call $\mathbf{S}$ a factor of $\bT$, and write $\mathbf{S} \subset \bT$, if the $i$th row of $\mathbf{S}$ is contained in that of $\bT$ (as multisets), for every $i \in [k]$. In this case, we define $\frac{\bT}{\mathbf{S}}=\mathbf{S}^{-1}\bT=\bT \mathbf{S}^{-1}$ to be the row-increasing tableau whose $i$th row is obtained by removing that of $\mathbf{S}$ from that of $\bT$ (as multisets), for every $i \in [k]$. A tableau $\bT \in {\rm SSYT}(k, [n])$ is {\em trivial} if each entry of $\bT$ is one less than the entry below it. For any $\bT \in {\rm SSYT}(k, [n])$, we  denote by $\bT_{\text{red}} \subset \bT$ the semistandard tableau obtained by removing a maximal trivial factor from $\bT$, i.e. $\bT_{\text{red}}=\bT \bT_{max}^{-1}$ for $\bT_{max}$ the (unique!) maximal trivial factor of $\bT$. If $\bT$ is itself trivial, then $\bT_{\text{red}}$ is the empty tableau. Let ``$\sim$'' be the equivalence relation on $\mathbf{S}, \bT \in {\rm SSYT}(k, [n])$ defined as follows: $\mathbf{S} \sim \bT$ if and only if $\mathbf{S}_{\text{red}} = \bT_{\text{red}}$. We denote by ${\rm SSYT}(k, [n],\sim)$ the set of $\sim$-equivalence classes. For example, consider $\mathbf{S} = \ytableausetup{centertableaux}
\scalebox{0.6}{\begin{ytableau}
1 & 2 \\
3 & 4 \\
4 & 6 
\end{ytableau}}$, $\bT = \ytableausetup{centertableaux}
\scalebox{0.6}{\begin{ytableau}
1 & 1 \\
2 & 4 \\
3 & 6 
\end{ytableau}}$ in ${\rm SSYT}(3, [6])$. Then $\mathbf{S}_{\rm red} = \bT_{\rm red} = \ytableausetup{centertableaux}
\scalebox{0.6}{\begin{ytableau}
1 \\
4 \\
6
\end{ytableau}}$ and therefore $\mathbf{S} \sim \bT$.

{Let $L_n(\lambda)$ denote the irreducible polynomial representation of $GL_n(\mathbb{C})$ with highest weight $\lambda \in \mathbb{Z}_{\ge 0}^n$. By the foundational results of Lusztig \cite{Lus90} and Kashiwara \cite{Kas91}, each $L_n(\lambda)$ admits a dual canonical basis as a $U_q(\mathfrak{gl}_n)$-module, and specializing at $q \to 1$ yields the corresponding basis of the classical $GL_n(\mathbb{C})$-module. The degree-$d$ component $\mathbb{C}[\mathrm{Gr}(k,n)]_{(d)}$, spanned by degree-$d$ monomials in Plücker coordinates, is an irreducible $GL_n(\mathbb{C})$-module with highest weight $d \omega_k$, and $\mathbb{C}[\mathrm{Gr}(k,n)]_{(d)} \cong L_n(d \omega_k)$. Therefore, $\mathbb{C}[\mathrm{Gr}(k,n)]_{(d)}$ inherits a dual canonical basis from $L_n(d \omega_k)$, and combining these bases for all $d \ge 1$ endows $\mathbb{C}[\mathrm{Gr}(k,n)]$ with a graded dual canonical basis.}

{The quotient $\CC[\Gr(k,n,\sim)]$ inherits a basis from the dual canonical basis of $\mathbb{C}[\mathrm{Gr}(k,n)]$, which we also refer to as the dual canonical basis of $\CC[\Gr(k,n,\sim)]$}. The elements of the dual canonical basis of $\CC[\Gr(k,n,\sim)]$ are in one to one correspondence with equivalence classes of semistandard Young tableaux in $\SSYT(k,[n],\sim)$, \cite[Section 3]{CDFL20}.
For simplicity, we also call the elements in $\SSYT(k,[n],\sim)$ tableaux instead of equivalence classes. So simple modules in $\mathscr{C}_{\ell}^{\mathfrak{sl}_k}$ are in one to one correspondence with tableaux in $\SSYT(k,[n],\sim)$, for $n=k+\ell+1$. 
Recall that simple modules correspond to dominant monomials (Section~\ref{subsec:qu_aff_algebras-qcharacters}). 
To indicate these correspondence, we write ${\bf T}_M$ to denote the tableau corresponding to the dominant monomial $M$ and $M_{\bf T}$ to denote the dominant monomial corresponding to the tableau ${\bf T}$. 

{
By the dual canonical basis of $K_{0}(\mathscr{C}_{\ell}^{\mathfrak{sl}_k})$, we mean the basis given by the equivalence classes of simple modules. The Grothendieck ring $K_0(\mathscr{C}_{\ell}^{\mathfrak{sl}_k})$ can be realized as a sub-ring of the Laurent polynomial ring $\mathbb{Z}[Y_{i,s}^{\pm 1}]_{i \in I,\; s \le \xi(i)}$. Under this realization, the basis elements of $K_0(\mathscr{C}_{\ell}^{\mathfrak{sl}_k})$ correspond to the truncated $q$-characters of simple modules.
}

\begin{remark}\label{rem:corresp}
In the sequel, we will often freely use the bijections of the following diagram (let $n=k+\ell+1$ as before): 
\[
\begin{tikzcd}
\{\mbox{simple modules in $\mathscr{C}_{\ell}^{\mathfrak{sl}_k}$}\} && \mbox{dual canonical basis of $K_0(\mathscr{C}_{\ell}^{\mathfrak{sl}_k})$} \\
%\\
\{\mbox{dominant monomials in $\mathcal{P}^+_{\ell}$}\} && \mbox{$\SSYT(k,[n],\sim)$} 
\arrow[tail reversed, from=1-1, to=1-3]
\arrow[tail reversed, from=1-1, to=2-1]
\arrow[tail reversed, from=1-3, to=2-3]
\arrow[tail reversed, from=2-1, to=2-3]
\end{tikzcd}
\]
The first vertical correspondence is due to Chari-Pressley~\cite{CP95a}, see Section~\ref{subsec:qu_aff_algebras-qcharacters}. {According to the first paragraph of this section 
(Section~\ref{subsec:correspondence between modules and tableaux}), 
the dual canonical basis of 
$K_0(\mathscr{C}_{\ell}^{\mathfrak{sl}_k})$ 
is in bijection with the dual canonical basis of 
$\CC[\Gr(k,n,\sim)]$, where $n=k+\ell+1$. 
By the second paragraph preceding this remark, the latter is in 
bijection with $\SSYT(k,[n],\sim)$. 
Combining these two bijections gives the second vertical correspondence.}
{The first horizontal correspondence follows from 
\cite[\S 5]{Rou12}, \cite[\S 4]{VV11}, \cite[\S 13]{HL10}, \cite[Theorem 5.1]{HL16}. This correspondence is given by sending a simple module to its truncated $q$-character.} 
We explain in Section~\ref{sec:simple-tableaux} below how to get the 
second horizontal correspondence, i.e. how to go from dominant monomials to tableaux. {The diagram is commutative in the appropriate sense.} 
\end{remark}

\begin{remark}\label{rem:g-vector_to_tableau}
{Using the correspondence between dual canonical basis of $\CC[\Gr(k,n)]$ and semistandard Young tableaux, the definition of (extended) $\mathbf{g}$-vector can be extended to arbitrary elements in the dual canonical basis of $\CC[\Gr(k,n)]$.}
% {\color{blue}The definition of ${\mathbf g}$-vectors can be extended to arbitrary elements in the dual canonical basis of $\CC[\Gr(k,n)]$, see Definition 4.5 in \cite{Qin20}. We recall this here.  An element of the dual 
% %The $\mathbf{g}$-vector of an element in the dual 
% canonical basis of $\CC[\Gr(k,n)]$ corresponds to a semistandard tableau $\mathbf{T}$, see~\cite[\S 3]{CDFL20}.} 
The (extended) ${\bf g}$-vector of the tableau $\mathbf{T}$ (and the corresponding element in the dual canonical basis) can be computed using the method in Section 7 of \cite{CDFL20}: 
Any tableau $\mathbf{T}$ can be written uniquely in the form $\mathbf{S}_1^{g_1} \cup \cdots \cup \mathbf{S}_m^{g_m}$ (\cite[Corollary 7.3]{CDFL20}) where the $\bf{S}_i$ are the tableaux corresponding to the elements of the initial cluster (choosing an order for them) and with $g_i\in \mathbb{Z}$. The vector $\hat{\bf{g}}:=(g_1,\dots, g_m){\in \mathbb{Z}^m}$ is called the {\em {extended} $\bf{g}$-vector of $\bf{T}$} {and ${\rm pr}_n(\hat{\mathbf{g}})\in \ZZ^n$ is called the {\em $\mathbf{g}$-vector} of $\bf{T}$. It is known that a semistandard Young tableau $\bT$ is uniquely determined by it (extended) $\mathbf{g}$-vector.
In the following,  we write $\mathbf{g}(\bT)$ for the $\mathbf{g}$-vector of a semistandard Young tableau $\bT$, and denote the tableau associated with an integer vector $\mathbf{g}$ by $\bT_{\mathbf{g}}$.}
\end{remark}

Using the above, the notions `real' and `prime' extend to tableaux and to {(extended)} $\mathbf{g}$-vectors: 

\begin{definition}\label{def:real-tableau_vector}
Consider a tableau $\bT$.  Recall that $M_{\bf T}$ is the dominant monomial corresponding to $\bT$. 
If  $L(M_{\bf T})$ is real (resp. prime), we say that $\bT$ is a {\em real} (resp. {\em prime}) tableau. 
If {$L(M_{\bf T})$} is non-real, we say that $\bT$ is non-real.

Let $\hat{\bf{g}}$ be the {extended}  $\bf{g}$-vector of $\bf{T}$ of the tableau $\bf T$ in $\SSYT(k,[n])$. If $\bf T$ is a real (resp. a prime) tableau, we say that it is a {\em real} (resp. {\em prime}) {extended} $\mathbf{g}$-vector. 
If $\bT$ is non-real, we say that 
$\hat{\bf{g}}$ is non-real. {Consequently, the $\mathbf{g}$-vector ${\rm pr}_n(\hat{\mathbf{g}})$ is also referred to as real or non-real accordingly.}
\end{definition}

%%%%%
%
\subsection{Correspondence between simple modules and tableaux}\label{sec:simple-tableaux}

Here, we explain the correspondence between {dominant monomials in $\mathcal{P}_\ell^+$} and semistandard Young tab\-leaux, recalling Section 3 of~\cite{CDFL20}. We first go from monomials to tableaux.

A {\em fundamental tableau} is a one-column semistandard tableau whose entries are of the form $\{i,i+1, \ldots, \widehat{j}, \ldots, i+k\}$, where $1\le i\le n-k$, $i < j < i+k$, and where $\widehat{j}$ means the number $j$ is omitted. We write $\bT_{i,s}$ to denote the tableau with entries 
$[\frac{i-s}{2}, k+ \frac{i-s}{2}]\setminus \{ k-\frac{i+s}{2}\}$ 
(note that by the choice of the height function $\xi$, $i\mapsto i-2$, for $i\in I=[k-1]$ and the definition of $\mathscr{C}_{\ell}$, the fraction $\frac{i-2}{2}$ is indeed an integer).

The 
correspondence between dominant monomials  and tableaux in $\SSYT(k,[n],\sim)$ is induced from associating  
the dominant monomial $Y_{i,s}$ with the fundamental tableau $\bT_{i,s}$.

If $M$ is an arbitrary dominant monomial, we write it as $M = \prod_{i,s} Y_{i,s}^{u_{i,s}}$, where the $u_{i,s}$ are positive integers. 
%Let 
%${\bf T}_{i,s}$ be the fundamental tableau corresponding to $Y_{i,s}$. 
Then the tableau ${\bf T}_M$ associated to $M$ is obtained by removing all trivial tableaux from $\cup_{i,s} {\bf T}_{i,s}^{\cup u_{i,s}}$. 

\vskip.2cm

For example, the simple module $L(M) = L(Y_{1,-5}Y_{1,-3}Y_{2,-2}Y_{2,0})$ corresponds to the semistandard Young tableau $\scalemath{0.6}{\begin{ytableau}
1 & 2 \\ 3 & 4 \\ 5 & 6
\end{ytableau}}$: the monomial  
$Y_{1,-5}$ gives $\ytableausetup{centertableaux}
\scalebox{0.6}{\begin{ytableau}
3 \\
4 \\
6
\end{ytableau}}$, the monomial 
$Y_{1,-3}$ gives $\ytableausetup{centertableaux}
\scalebox{0.6}{\begin{ytableau}
2 \\
3 \\
5
\end{ytableau}}$, the monomial 
$Y_{2,-2}$ gives $\ytableausetup{centertableaux}
\scalebox{0.6}{\begin{ytableau}
2 \\
4 \\
5
\end{ytableau}}$ and $Y_{2,0}$ gives $\ytableausetup{centertableaux}
\scalebox{0.6}{\begin{ytableau}
1 \\
3 \\
4
\end{ytableau}}$. Then we take the union of these tableaux and obtain $\ytableausetup{centertableaux}
\scalebox{0.6}{\begin{ytableau}
1&2&2&3 \\
3&3&4&4 \\
4&5&5&6
\end{ytableau}} \sim \ytableausetup{centertableaux}
\scalebox{0.6}{\begin{ytableau}
1&2 \\
3&4 \\
5&6
\end{ytableau}} = \bT_M$. 

Now we explain how to go from tableaux to monomials.
Given a tableau $\bT$, there is a unique semistandard tableau $\bT'$ such that $\bT \sim \bT'$ and each column of $\bT'$ is a fundamental tableau. The dominant monomial $M_{\bT}$ associated to $\bT$ is the product of the $Y_{i,s}$'s which correspond to the columns of $\bT'$. For example, let $\bT = \ytableausetup{centertableaux}
\scalebox{0.6}{\begin{ytableau}
1&2 \\
3&4 \\
5&6
\end{ytableau}}$. Then $\bT \sim \bT' = \ytableausetup{centertableaux}
\scalebox{0.6}{\begin{ytableau}
1&2&2&3 \\
3&3&4&4 \\
4&5&5&6
\end{ytableau}}$. 
The columns of $\bT'$ are the fundamental tableaux 
$\bT_{2,0}$, $\bT_{1,-3}$, $\bT_{2,-2}$ and $\bT_{1,-5}$ in that order. From this, we obtain 
the dominant monomial corresponding to $\bT$, it is  $M_{\bT}=Y_{2,0}Y_{1,-3}Y_{2,-2}Y_{1,-5}$.

\vskip.3cm

\section{Correspondence between additive and monoidal categorifications} \label{sec:correspondence between additive and monoid categorifications}

%%
%
%\subsection{The setting}

Let $\widetilde{B}\in M_{m\times n}(\mathbb{Z})$ whose principal part is skew-symmetric and $\mathcal{A}:=\mathcal{A}(\widetilde{B})$ be the associated cluster algebra with coefficients. 
We fix an initial seed $(\mathbf{x}_{t_0},\widetilde{B})$ of $\mathcal{A}$.
  
Assume that $\mathcal{A}$ admits an additive categorification by $(\mathcal{F}, \hat{T})$ as in Section \ref{ss:additive-cat}. In particular, $\mathcal{F}$ is a $\mathbb{k}$-linear Frobenius (extriangulated) category which is stably $2$-Calabi-Yau.
Write $\mathcal{C}$ for its stable category and let $T$ be the basic cluster-tilting object in $\mathcal{C}$ corresponding to $\hat{T}$. We have $B_{T}=\widetilde{B}$.

%, with cluster structure which is stably $2$-Calabi-Yau and  $T$ is a basic cluster-tilting object of $\mathcal{C}$ such that $B_T=\tilde{B}$.

We assume moreover that $\mathcal{A}$ also admits a monoidal categorification $(\mathscr{C}_\ell,\otimes)$, where $\mathscr{C}_\ell$ is the full subcategory of the category $\mathscr{C}$ of finite dimensional $U_q(\hat{\mathfrak{g}})$-modules {introduced in Section \ref{ss:monoidal-cat}} and $\mathfrak{g}$ is a complex simple Lie algebra. 
In particular, we have 
\begin{itemize}
\item $\mathcal{A}\cong K_0(\mathscr{C}_\ell)$, which is injective-reachable in the sense of \cite{Qin17};
\item The set $\{\chi_q(L(M))~|~\text{$M$ is a dominant monomial in $\mathcal{P}_\ell^+$}\}$ is a basis of $\mathcal{A}$ which contains cluster monomials;
\item {The basis $\{\chi_q(L(M))~|~\text{$M$ is a dominant monomial in $\mathcal{P}_\ell^+$}\}$ is parametrized by a subset of $\mathbb{Z}^m$. In particular, for each dominant monomial $M\in \mathcal{P}_\ell^+$, there is a unique vector $\hat{\mathbf{g}}_M\in \ZZ^m$ and a polynomial $F_M(z_1,\ldots, z_n)\in \ZZ[z_1,\ldots, z_n]$ such that 
\[
\chi_q(L(M))=\mathbf{x}^{\hat{\mathbf{g}}_M}F_M(\hat{y_1},\ldots, \hat{y}_n).
\]
The vector $\hat{\mathbf{g}}_M$ is called the {\em extended $\mathbf{g}$-vector} of $\chi_q(L(M))$ (or of $M$) with respect to the initial cluster of $\mathcal{A}$. Its projection $\mathbf{g}_M:={\rm pr}_n(\hat{\mathbf{g}}_M)\in \mathbb{Z}^n$ is called the {\em $\mathbf{g}$-vector} of $\chi_q(L(M))$ (or of $M$). Furthermore, there exists a method to compute the vector $\hat{\mathbf{g}}_M$ by expressing $M$ using the highest $l$-weight monomials of the modules corresponding to the initial cluster variables, see \cite[Section 5.2.2]{HL16} and \cite[Section 7]{CDFL20}.
If $M,M'\in \mathcal{P}_\ell^+$ such that $\mathbf{g}_M=\mathbf{g}_{M'}$, then there exist monomials $\mathbf{x}^a,\mathbf{x}^b$ in frozen variables such that $\chi_q(L(M))\mathbf{x}^a=\chi_q(L(M'))\mathbf{x}^b$, where $a,b\in \mathbb{N}^m$ such that ${\rm pr}_n(a)=0={\rm pr}_n(b)$ (cf. \cite[Proposition 9.1.5]{Qin17}).}
\end{itemize}
\begin{remark}{
We may localize the cluster algebra $\mathcal{A}$ at all frozen variables to obtain the cluster algebra $\mathcal{A}_{\rm loc}$.
  The localization of $\{\chi_q(L(M))~|~\text{$M$ is a dominant monomial in $\mathcal{P}_\ell^+$}\}$ at the frozen variables yields a basis of $\mathcal{A}_{\rm loc}$, which is parametrized by $\ZZ^m$, see \cite[Section 9]{Qin17}.}
\end{remark}

In this section, we study the connection between the additive categorification $(\mathcal{F},\hat{T})$ and the monoidal categorification $\mathscr{C}_\ell$ of $\mathcal{A}$. 

\begin{remark}
The statements in this section  can be formulated for an arbitrary injective-reachable cluster algebra of skew-symmetric type with an (abstract) monoidal categorification. Indeed, the construction only relies on the existence of an additive categorification and such an additive categorification always exists by Remark \ref{rem:existence-additive-cat}. However, the monoidal reachability conjecture is only formulated for monoidal categorifications arising from representation theory of quantum affine algebras, hence we restrict ourselves to the  special case.
\end{remark}

\subsection{The function $\mathfrak{e}(-,-)$}\label{ss:function-e}
We recall the function 
$\mathfrak{e}:\mathbb{Z}^n\times \mathbb{Z}^n\to \mathbb{Z}$ from~\cite{DF15} and~\cite{P13}.

%Following \cite{DF15} and \cite{P13}, we introduce a function  $\mathfrak{e}:\mathbb{Z}^n\times \mathbb{Z}^n\to \mathbb{Z}$.
Let $T$ be a cluster-tilting object of $\mathcal{C}$. Fix an integer vector $\mathbf{g}\in \mathbb{Z}^n$. 
Let $\mathbf{g}_+$ and $\mathbf{g}_-$ be as defined in Section~\ref{sec:notation}. 
Consider the finite dimensional 
%and we have the finite dimensional 
vector space $\Hom_\mathcal{C}(T^{\mathbf{g}_-},T^{\mathbf{g}_+})$. There is an open dense subset  $\mathcal{O}_\mathbf{g}$ of $\Hom_\mathcal{C}(T^{\mathbf{g}_-},T^{\mathbf{g}_+})$ such that for each $f\in \mathcal{O}_{\mathbf{g}}$, $\Hom_\mathcal{C}(T,\Sigma \operatorname{cone}(f))$ has the same $F$-polynomial (cf. \cite[Proposition 3.16]{P13})
\[
\sum_e\chi(\text{Gr}_e\Hom_\mathcal{C}(T,\Sigma \operatorname{cone}(f))\mathbf{x}^e,
\]
where $\operatorname{cone}(f)$ is the mapping cone of $f$ in $\mathcal{C}$.  We remark that the open dense subset $\mathcal{O}_\mathbf{g}$ is not unique.

Let $K^b(\add T)$ be the bounded homotopy category of $\add T$. For $T_0,T_{-1}\in \add T$ and $f\in \Hom_\mathcal{C}(T_{-1},T_0)$, we may identify $f$ with a two term complex $T_f^\bullet:T_{-1}\xrightarrow{f}T_0\in K^b(\add T)$, where $T_0$ is the zero component of $T_f^\bullet$.
\begin{definition}\cite[Definition 3.1]{DF15}
    Let $T_1,T_2,T_3,T_4\in \add T$ and $f\in \Hom_\mathcal{C}(T_2,T_1), g\in \Hom_\mathcal{C}(T_4,T_3)$. The $\operatorname{E}$-invariant $\operatorname{E}(f,g)$ from $f$ to $g$ is defined as
    \[
    \operatorname{E}(f,g):=\dim_\mathbb{k}\Hom_{K^b(\add T)}(f,\Sigma g).
    \]
    Define $\mathbb{E}(-,-):\Hom_\mathcal{C}(T_2,T_1)\times \Hom_\mathcal{C}(T_4,T_3)\to \mathbb{Z}$ by 
    \[
    \mathbb{E}(f,g):=\operatorname{E}(f,g)+\operatorname{E}(g,f).
    \]
   
\end{definition}
\begin{lemma}\label{l:E-invariant-ext}
    Let $T_1,T_2,T_3,T_4\in \add T$ and $f\in \Hom_\mathcal{C}(T_2,T_1), g\in \Hom_\mathcal{C}(T_4,T_3)$. Then $\mathbb{E}(f,g)=\dim_\mathbb{k}\Ext^1_{\mathcal{C}}(\operatorname{cone}(f),\operatorname{cone}(g))$.
\end{lemma}
\begin{proof}
Denote by $X:=\operatorname{cone}(f)$ and $Y:=\operatorname{cone}(g)$.
    According to \cite[Proposition 3.10]{P13}, $\operatorname{E}(f,g)$ is isomorphic to  the subspace $(\Sigma T)(X, \Sigma Y)$  of $\Hom_\mathcal{C}(X, \Sigma Y)$ consisting of morpshims which are factorizing through objects in $\add \Sigma T$. 
    On the other hand, by \cite[Lemma 3.3]{Pa08}, we have $(\Sigma T)(Y, \Sigma X)\cong D\Hom_{\mathcal{C}/(\Sigma T)}(X,\Sigma Y)$, where $\Hom_{\mathcal{C}/(\Sigma T)}(X,\Sigma Y):=\Hom_\mathcal{C}(X,\Sigma Y)/(\Sigma T)(X,\Sigma Y)$.
    If follows that
    \[
    \mathbb{E}(f,g)=\dim_\mathbb{k} (\Sigma T)(X,\Sigma Y)+\dim_\mathbb{k} \Hom_{\mathcal{C}/(\Sigma T)}(X,\Sigma Y)=\dim_\mathbb{k} \Hom_\mathcal{C}(X,\Sigma Y).
    \]
\end{proof}
\begin{remark}\label{rem:E-inv-history}
$E$-invariants were first introduced by Derksen, Weyman and Zelevinsky~\cite{DWZ} for decorated representations of Jacobian algebras. Later, Derksen and Fei \cite{DF15} introduced $E$-invariants for any finite-dimensional algebras (cf. also \cite{AIR14}).
% Adachi, Iyama and Reiten~\cite{AIR14} introduced $E$-invariants for modules for $2$-Calabi-Yau tilted algebras.
Lemma~\ref{l:E-invariant-ext} shows that the symmetrized $E$-invariant $\mathbb{E}(f,g)$ coincides with the $E$-invariant of $\operatorname{cone}(f)$ and $\operatorname{cone}(g)$.
\end{remark}

\begin{lemma}\label{l:constructible-fuction}
Let $T_1,T_2,T_3,T_4\in \add T$.
    The function \[
    \mathbb{E}(-,-):\Hom_\mathcal{C}(T_2,T_1)\times \Hom_\mathcal{C}(T_4,T_3)\to \mathbb{Z}\] is upper semi-continuous and constructible. 
\end{lemma}
\begin{proof}
    The function $\operatorname{E}(-,-)$ is upper semi-continuous (cf. \cite[Section 3]{DF15}), which implies that $\mathbb{E}(-,-)$ is upper semi-continuous.  By definition, for any $f\in \Hom_{\mathcal{C}}(T_2,T_1)$ and $g\in \Hom_{\mathcal{C}}(T_4,T_3)$, $\operatorname{E}(f,g)\leq \dim_k\Hom_\mathcal{C}(T_2,T_3)$ and $\operatorname{E}(g,f)\leq \dim_k\Hom_\mathcal{C}(T_4,T_1)$. As a consequence, we have $0\leq \mathbb{E}(f,g)\leq \dim_k\Hom_\mathcal{C}(T_2,T_3)+\dim_k\Hom_\mathcal{C}(T_4,T_1)$, which implies that the image $\operatorname{im} \mathbb{E}$ of $\mathbb{E}$ is a finite subset of $\mathbb{Z}$. Since $\mathbb{E}$ is upper continuous, we conclude that $\mathbb{E}^{-1}(<t):=\{(f,g)~|{\mathbb{E}}(f,g)<t\}$ is an open subset  and $\mathbb{E}^{-1}(\geq t):=\{(f,g)~|~{\mathbb{E}}(f,g)\geq t\}$ is a closed subset of $\Hom_\mathcal{C}(T_2,T_1)\times \Hom_\mathcal{C}(T_4,T_3)$ for any $t\in \mathbb{Z}$. Consequently, $\mathbb{E}^{-1}(t)=\mathbb{E}^{-1}(<t+1)\cap \mathbb{E}^{-1}(\geq t)$ is a locally closed subset, which implies that $\mathbb{E}$ is constructible.
\end{proof}
\begin{definition}
    {Let $\mathbf{g}, \mathbf{h}\in \mathbb{Z}^n$ and let $\mathcal{O}_\mathbf{g}$ and $\mathcal{O}_\mathbf{h}$ be two open dense subsets} associated with $\mathbf{g}$ and $\mathbf{h}$ respectively. Consider the restriction of $\mathbb{E}$ to 
    $\mathcal{O}_\mathbf{g}\times\mathcal{O}_\mathbf{h}$, %restriction of 
    %$\mathbb{E}(-,-):\mathcal{O}_\mathbf{g}\times \mathcal{O}_\mathbf{h}\to \mathbb{Z}$. 
    Then $\mathbb{E}(-,-)$ admits a generic value, that is, there is a dense open subset of the domain of the function on which the function is constant. We denote the generic value of $\mathbb{E}$ on $\mathcal{O}_\mathbf{g}\times \mathcal{O}_\mathbf{h}$ by $\mathfrak{e}(\mathbf{g}, \mathbf{h})$. 
    %the generic value of $\mathbb{E}(-,-)$ on $\mathcal{O}_\mathbf{g}\times \mathcal{O}_\mathbf{h}$.

\end{definition}
\begin{remark}\label{r:e-function-minimal}
    It is obvious that the value $\mathfrak{e}(\mathbf{g},\mathbf{h})$ does not depend on the choice of the dense open subset of $\mathcal{O}_\mathbf{g}\times \mathcal{O}_\mathbf{h}$ nor the choice of the dense subsets $\mathcal{O}_\mathbf{g}$ and $\mathcal{O}_\mathbf{h}$. In fact, 
    \[
    \mathfrak{e}(\mathbf{g},\mathbf{h})=\min \{\dim_\mathbb{k} \Ext^1_\mathcal{C}(M,N)\},
    \]
    where $M$ and $N$ range over the mapping cones of morphisms in $\Hom_\mathcal{C}(T^{\mathbf{g}_-},T^{\mathbf{g}_+})$ and $\Hom_\mathcal{C}(T^{\mathbf{h}_-},T^{\mathbf{h}_+})$ respectively.
\end{remark}

\begin{lemma}\label{l:ext-dim}
     {Let $\mathbf{g}, \mathbf{h}\in \mathbb{Z}^n$ and let $\mathcal{O}_\mathbf{g}$ and $\mathcal{O}_\mathbf{h}$ be two open dense subsets} associated with $\mathbf{g}$ and $\mathbf{h}$, respectively. Assume that there exist $g\in \mathcal{O}_{\mathbf{g}}$ and $h\in \mathcal{O}_\mathbf{h}$ such that $\operatorname{cone}(g)$ and $\operatorname{cone}(h)$ are rigid, then 
    \[\mathfrak{e}(\mathbf{g},\mathbf{h})=\dim_\mathbb{k}\Ext_\mathcal{C}^1(\operatorname{cone}(g),\operatorname{cone}(h)). 
    \]
    In particular, $\mathfrak{e}(\mathbf{g},\mathbf{g})=0$.
\end{lemma}

\begin{proof}
    There is a free $\operatorname{Aut}(T^{\mathbf{g}_-})\times \operatorname{Aut}(T^{\mathbf{g}_+})$ action on $\Hom_\mathcal{C}(T^{\mathbf{g}_-},T^{\mathbf{g}_+})$. Since $\operatorname{cone}(g)$ is rigid,  the orbit $\mathcal{O}(g)$ of $g$ under the group action is an open dense subset of the space $\Hom_\mathcal{C}(T^{\mathbf{g}_-},T^{\mathbf{g}_+})$. Hence $\mathcal{O}(g)\cap \mathcal{O}_\mathbf{g}\neq \emptyset$. Similarly, the orbit $\mathcal{O}(h)$ is also open and dense such that $\mathcal{O}(h)\cap \mathcal{O}_\mathbf{h}\neq \emptyset$. In particular, $\mathcal{Z}:=(\mathcal{O}(g)\times \mathcal{O}(h))\cap (\mathcal{O}_{\mathbf{g}}\times \mathcal{O}_{\mathbf{h}})$ is an open dense subset of $\mathcal{O}_{\mathbf{g}}\times \mathcal{O}_{\mathbf{h}}$. It is clear that $\mathbb{E}(-,-)$ is constant on $\mathcal{Z}$, which implies that $\mathfrak{e}(\mathbf{g},\mathbf{h})=\dim_\mathbb{k}\Ext_\mathcal{C}^1(\operatorname{cone}(g),\operatorname{cone}(h))$.
\end{proof}

\subsection{The function $\mathfrak{e}(-)$} \label{subsec:the function e minus}

Let $\mathbf{g}\in \mathbb{Z}^n$ and $\mathcal{O}_\mathbf{g}$ an open dense subset as in Section \ref{ss:function-e}.
Recall that for any $f\in \Hom_\mathcal{C}(T^{\mathbf{g}_-},T^{\mathbf{g}_+})$, we have a two term complex $T^\bullet_f:T^{\mathbf{g}_-}\xrightarrow{f} T^{\mathbf{g}_+}\in K^b(\add T)$. 
The function \begin{eqnarray*}
\mathfrak{e}:&\Hom_\mathcal{C}(T^{\mathbf{g}_-},T^{\mathbf{g}_+})&\to \mathbb{Z}\\
&f&\mapsto \dim_\mathbb{k}\Hom_{K^b(\add T)}(T_f^\bullet,\Sigma T_f^\bullet)
\end{eqnarray*}
is upper semi-continuous.
Similarly as for $\mathbb{E}$ (Lemma \ref{l:constructible-fuction}), one can show that the function $\mathfrak{e}$ is constructible.
Therefore, the restriction of $\mathfrak{e}$ to $\mathcal{O}_\mathbf{g}$ admits a generic value, that is, there is a dense open subset $U$ of $\mathcal{O}_\mathbf{g}$ on which the function is constant. It is clear that the value is independent of the choice of the open subset $U$ nor the choice of $\mathcal{O}_\mathbf{g}$. We denote the generic value on $\mathcal{O}_\mathbf{g}$ by $\mathfrak{e}(\mathbf{g})$.

\begin{lemma}\label{l:rigid-e-fuction}
Let $\mathbf{g}\in \mathbb{Z}^n$ and let $\mathcal{O}_\mathbf{g}$ be an open dense subset associated with $\mathbf{g}$. There exists a $g\in \mathcal{O}_\mathbf{g}$ such that $\operatorname{cone}(g)$ is rigid if and only if $\mathfrak{e}(\mathbf{g})=0$.
\end{lemma}
\begin{proof}
Let $U$ be \textcolor{red}{a} open dense subset of $\mathcal{O}_{\mathbf{g}}$ such that $\mathfrak{e}(-)$ is constant on $U$. Denote by $\mathcal{O}(g)$ the orbit  of $g$ under the group action of $\operatorname{Aut}(T^{\mathbf{g}_-})\times \operatorname{Aut}(T^{\mathbf{g}_+})$. 

    If $\operatorname{cone}(g)$ is rigid, then $\mathcal{O}(g)$ is open dense, thus $\mathcal{O}(g)\cap U\neq \emptyset$. It follows that $\mathfrak{e}(\mathbf{g})=\dim_\mathbb{k}\Hom_{K^b(\add T)}(T_g^\bullet, \Sigma T_g^\bullet)$. By Lemma \ref{l:E-invariant-ext}, \[\dim_\mathbb{k}\Hom_{K^b(\add T)}(T_g^\bullet, \Sigma T_g^\bullet)=\frac{1}{2}\dim_{\mathbb{k}} \Ext^1_\mathcal{C}(\operatorname{cone}(g),\operatorname{cone}(g))=0,\] which implies that $\mathfrak{e}(\mathbf{g})=0$.

    Assume that $\mathfrak{e}(\mathbf{g})=0$. By definition, there is a morphism $g\in \Hom_{\mathcal{C}}(T^{\mathbf{g}_-},T^{\mathbf{g}_+})$ such that $\dim_\mathbb{k}\Hom_{K^b(\add T)}(T^\bullet_g, \Sigma T^\bullet_g)=0$. Again by Lemma \ref{l:E-invariant-ext}, we conclude that $\dim_\mathbb{k}\Ext^1_\mathcal{C}(\operatorname{cone}(g), \operatorname{cone}(g))=0$.
\end{proof}

\begin{remark}
Let $\mathbf{g}\in \mathbb{Z}^n$. It is clear that $\mathfrak{e}(\mathbf{g})=0$ implies that $\mathfrak{e}(\mathbf{g},\mathbf{g})=0$. However, the converse is not true in general, cf. Remark \ref{r:e(g)-neq-e(g,g)}.
\end{remark}

\subsection{From monoidal categorification to additive categorification}
 
% Let $\mathcal{A}$ be a cluster algebra with initial cluster $\mathbf{x}_{t_0}=\{x_1,\dots,x_m\}$). Assume that $\mathcal{A}$ has a monoidal categorification 
%  $\mathscr{C}_\ell$ and an additive categorification 
% $(\mathcal{F},\hat{T})$ and let $\mathcal{C}$ be the stable category of $\mathcal{F}$. 

Let $M$ be a dominant monomial in $\mathcal{P}_\ell^+$ and $L(M)$ the associated simple $U_q(\hat{\mathfrak{g}})$-module in $\mathscr{C}_\ell$. Recall that $\hat{\mathbf{g}}_{M}\in \mathbb{Z}^m$ is the extended $\mathbf{g}$-vector of $\chi_q(L(M))$ in the  cluster algebra $\mathcal{A}$ with initial cluster $\mathbf{x}_{t_0}=\{x_1,\dots,x_m\}$. {According to the discussion at the beginning of Section \ref{sec:correspondence between additive and monoid categorifications}, the following map is an injection
\begin{eqnarray*}
    &\{\text{simple modules of $\mathscr{C}_\ell$}\}&\to \mathbb{Z}^m.\\
    &L(M)&\mapsto \hat{\mathbf{g}}_M
\end{eqnarray*}
Recall that $\mathbf{g}_M={\rm pr}_n(\hat{\mathbf{g}}_M)\in \ZZ^n$.}
Combining it with the construction in Section \ref{ss:function-e}, we obtain a map
\begin{eqnarray*}
    &\{\text{simple modules of $\mathscr{C}_\ell$}\}&\to \{\mathcal{O}_{\mathbf{g}}\subset \Hom_\mathcal{C}(T^{\mathbf{g}_-},T^{\mathbf{g}_+})\mid \mathbf{g}\in\mathbb{Z}^n\}.\\
    &L(M)&\mapsto \mathcal{O}_{\mathbf{g}_M}
\end{eqnarray*}
Note that $\mathcal{O}_\mathbf{g}$ can be understood as a `generic' subcategory of $\mathcal{C}$ by taking mapping cone of each morphism in $\mathcal{O}_\mathbf{g}$.

%\begin{remark}
%For any $\mathbf{g}\in \mathbb{Z}^n$, fix a generic homomorphism $\eta\in \mathcal{O}_\mathbf{g}$. Denote by $M(\mathbf{g}):=M(\mathbf{g},\eta)$  the mapping cone of $\eta$.
%Under a mild condition (i.e., full rank of exchange matrix and injective-reachable),
    %$\{\text{CC}(M(\mathbf{g}))~|~\mathbf{g}\in \mathbb{Z}^n\}$ is a basis of $\mathcal{U}$, which is called the {\it generic basis}, where $\mathcal{U}$ is the upper cluster algebra associated with $\mathcal{A}$. On the other hand, Both the bases contains all the cluster monomials. It has been proved that every cluster monomial is the character $\chi_q(L(\mathbf{m}))$ for a real module $L(\mathbf{m})$.
%\end{remark}

\begin{definition}
A dominant monomial $M\in \mathcal{P}_\ell^+$ is {\it reachable} if $\chi_q(L(M))$ is a cluster monomial of $\mathcal{A}$. {In this case we also call the module $L(M)$ is reachable.}
\end{definition}
The following is a direct consequence of Lemmas \ref{l:ext-dim} and \ref{l:rigid-e-fuction}.

\begin{proposition}
For every reachable dominant monomial  $M\in \mathcal{P}_\ell^+$, we have 
\[ \mathfrak{e}(\mathbf{g}_M,\mathbf{g}_M)=0 \mbox{ and }\mathfrak{e}(\mathbf{g}_M)=0.\] 
\end{proposition}

In general, we conjecture that the following holds: 
\begin{conjecture}\label{c:monoidal-additive-correspondence}
Let $M$ be a dominant monomial in $\mathcal{P}_\ell^+$. Then $L(M)$ is real if and only if $\mathfrak{e}(\mathbf{g}_M)=0$.
  
% \item[(3)] Let $\mathbf{m}_1,\mathbf{m}_2$ be dominant monomials such that $\mathbf{g}_{\mathbf{m}_1}\not\in \mathbb{N} \mathbf{g}_{\mathbf{m}_2}$ and $\mathbf{g}_{\mathbf{m}_2}\not\in \mathbb{N} \mathbf{g}_{\mathbf{m}_1}$. Then $\chi_q(L(\mathbf{m}_1))\chi_q(L(\mathbf{m}_2))=\chi_q(L(\mathbf{m}_1\mathbf{m}_2))$ if and only if 
\end{conjecture}

The following statement reveals the relationship between Conjectures \ref{c:addtive-reach-conj}, \ref{c:monoidal-reach-conj} and \ref{c:monoidal-additive-correspondence}.

\begin{proposition}\label{p:equivalence-add-mult}
Any two of Conjecture \ref{c:addtive-reach-conj}, Conjecture \ref{c:monoidal-reach-conj} and Conjecture \ref{c:monoidal-additive-correspondence} imply the other one.
%Assume that Conjecture \ref{c:monoidal-additive-correspondence} holds. Then Conjecture \ref{c:addtive-reach-conj} is equivalent to Conjecture \ref{c:monoidal-reach-conj}.
\end{proposition}
\begin{proof}
We first show that Conjecture \ref{c:addtive-reach-conj} is equivalent to Conjecture \ref{c:monoidal-reach-conj} under Conjecture \ref{c:monoidal-additive-correspondence}.
According to Conjecture \ref{c:monoidal-additive-correspondence} and Lemma~\ref{l:rigid-e-fuction}, { a simple module $L(M)$ of $\mathscr{C}_\ell$} is real if and only if there is a morphism $f\in \mathcal{O}_{\mathbf{g}_M}$ such that $\operatorname{cone}(f)$ is rigid. 

First assume that Conjecture \ref{c:addtive-reach-conj} is true. {Let $L(M)$ be a fixed real module.} Since the cluster algebra $\mathcal{A}$ is injective-reachable, it follows that the exchange graph $\mathcal{E}(\mathcal{C})$ is connected by Conjecture \ref{c:addtive-reach-conj}. As a consequence, $\operatorname{cone}(f)$ corresponds to a cluster monomial. Hence $\mathbf{g}_{M}$ is a $\mathbf{g}$-vector of a cluster monomial. Let $\hat{\mathbf{g}}_{\operatorname{cone}(f)}$ be the extended $\mathbf{g}$-vector of the cluster monomial $\mathbf{X}_{\operatorname{cone}(f)}^{\hat{T}}$. It follows that ${\rm pr}_n(\hat{\mathbf{g}}_{\operatorname{cone}(f)})={\rm pr}_n(\hat{\mathbf{g}}_M)$. Consequently, there are projective objects $P,Q$ of $\mathcal{F}$ such that 
\[
\chi_q(L(M))\mathbf{X}_P^{\hat{T}}=\mathbf{X}_{\operatorname{cone}(f)}^{\hat{T}} \mathbf{X}_Q^{\hat{T}}.
\]
We conclude that  $\chi_q(L(M))$ is a cluster monomial by \cite[Lemma 5.1]{CKQ24}.

Conversely, for any rigid object $N\in \mathcal{C}$, let $\mathbf{g}_N:={\rm pr}_n(\operatorname{ind}_{\hat{T}}N)$. Then $\mathfrak{e}(\mathbf{g}_N)=0$ { by Lemma \ref{l:rigid-e-fuction}.} Let $M_{\mathbf{g}_N}$ be a dominant monomial such that ${\rm pr}_n(\hat{\mathbf{g}}_M)=\mathbf{g}_N$. Consequently, $L(M_{\mathbf{g}_N})$ is real by Conjecture \ref{c:monoidal-additive-correspondence}. According to Conjecture \ref{c:monoidal-reach-conj}, $L(M_{\mathbf{g}_N})$ is a cluster monomial, which implies that $\mathbf{g}_N$ is the $\mathbf{g}$-vector of a cluster monomial. By Theorem \ref{t:additve-cat}, there exists a rigid object $L \oplus P$ in $\mathcal{F}$ that is reachable from $T$, where $P$ is projective and $L$ has no non-zero projective direct summands, such that\[\operatorname{pr}_n(\operatorname{ind}_{\widehat{T}}(L\oplus P)) = \mathbf{g}_N.
\]
It follows from \cite[Theorem 2.3]{DK08} that $L \cong N$, which implies that $N$ is reachable from $T$.

It remains to show that Conjecture \ref{c:addtive-reach-conj} and Conjecture \ref{c:monoidal-reach-conj} imply Conjecture \ref{c:monoidal-additive-correspondence}. Let $M$ be a dominant monomial in $\mathcal{P}_l^+$. Assume that $L(M)$ is real, then by Conjecture \ref{c:monoidal-reach-conj}, $\chi_q(L(M))$ is a cluster monomial. It follows from Theorem \ref{t:additve-cat} that there is a rigid object $N\in \mathcal{F}$ such that $\mathbf{X}^{\hat T}_N=\chi_q(L(M))$ and $\operatorname{ind}_{\hat{T}}N=\hat{\mathbf{g}}_M$. We conclude that $\mathfrak{e}(\mathbf{g}_M)=0$ by Lemma \ref{l:rigid-e-fuction}.

Now assume that $\mathfrak{e}(\mathbf{g}_M)=0$. Again by Lemma \ref{l:rigid-e-fuction}, there exists $g\in \mathcal{O}_{\mathbf{g}_M}$ such that $\Ext^1_\mathcal{C}(\operatorname{cone}(g),\operatorname{cone}(g))=0$.
     By Conjecture \ref{c:addtive-reach-conj} and Theorem \ref{t:additve-cat}, we know that $\mathbf{X}^{\hat T}_{\operatorname{cone}(g)}$ is a cluster monomial. Let $\hat{\mathbf{g}}_{\operatorname{cone}(f)}$ be the extended $\mathbf{g}$-vector of the cluster monomial $\mathbf{X}^{\hat T}_{\operatorname{cone}(g)}$. It follows that ${\rm pr}_n(\hat{\mathbf{g}}_{\operatorname{cone}(f)})=\mathbf{g}_M$. Consequently, there are projective objects $P,Q$ of $\mathcal{F}$ such that 
     \[
     \chi_q(L(M))\mathbf{X}_P^{\hat{T}}=\mathbf{X}_{\operatorname{cone}{f}}^{\hat{T}}\mathbf{X}_Q^{\hat{T}}.
     \]
     We conclude that $\chi_q(L(M))$ is a cluster monomial, which implies that $L(M)$ is real.

\end{proof}

\begin{remark}
    The proof of Proposition \ref{p:equivalence-add-mult} also shows the following implications: 
    \begin{itemize}
        \item Conjecture \ref{c:addtive-reach-conj} implies the``$\Leftarrow$" part of Conjecture \ref{c:monoidal-additive-correspondence};
        \item Conjecture \ref{c:monoidal-reach-conj} implies the``$\Rightarrow$" part of Conjecture \ref{c:monoidal-additive-correspondence}.
    \end{itemize}
\end{remark}

% \begin{theorem}
% Suppose that the exchange graph $\mathcal{E}(\mathcal{C})$ %of $\mathcal{C}$ 
% is connected.
% If $M$ is a dominant monomial such that $\mathfrak{e}(\mathbf{g}_M)=0$, then $L(M)$ is real.
% \end{theorem}
% \begin{proof}
% By the assumption, there exists $g\in \mathcal{O}_{\mathbf{g}_M}$ such that $\Ext^1_\mathcal{C}(\operatorname{cone}(g),\operatorname{cone}(g))=0$.
%      By the connectedness of the exchange graph of $\mathcal{C}$, we know that $\mathbf{X}^T_{\operatorname{cone}(g)}$ is a cluster monomial. Note that $\mathbf{X}^T_{\operatorname{cone}(g)}$ and $\chi_q(L(M))$ have the same $\mathbf{g}$-vector $\mathbf{g}_M$. It follows that $\chi_q(L(M))=\mathbf{X}^T_{\operatorname{cone}(g)}$, which implies that $L(M)$ is real.
% \end{proof}

We can now formulate a conjecture analogous to the statement of Lemma \ref{l:prod-cluster-monomial} and of Corollary \ref{c:prod-cluster-monomial} on the compatibility of cluster monomials.
\begin{conjecture}\label{c:monoidal-additive-corr-compatible}
Let $M_1,M_2\in \mathcal{P}_\ell^+$ be dominant monomials such that $L(M_1)$ and $L(M_2)$ are real. Then we have 
\[
\chi_q(L(M_1))\chi_q(L(M_2))=\chi_q(L(M_1M_2)) \mbox{ if and only if } \mathfrak{e}(\mathbf{g}_{M_1},\mathbf{g}_{M_2})=0.
\]

\end{conjecture}

\begin{proposition} \label{thm:conjecture-true-finite-type} 
Assume that $\mathcal{A}$ is of finite type. Then Conjectures \ref{c:monoidal-additive-correspondence} and \ref{c:monoidal-additive-corr-compatible} hold. %for $\mathcal{A}$. 
\end{proposition}
\begin{proof}
Since $\mathcal{A}$ is of finite type,  the cluster monomials form a basis of $\mathcal{A}$. {By the discussion at the beginning of Section \ref{sec:correspondence between additive and monoid categorifications}, $\chi_q(L(M))$ is a cluster monomial for any dominant monomial $M\in \mathcal{P}_{\ell}^+$. As a consequence, every dominant monomial $M\in \mathcal{P}_\ell^+$ is real and Conjecture \ref{c:monoidal-reach-conj} holds for $\mathcal{A}$. On the other hand, Conjecture \ref{c:addtive-reach-conj} also holds for $\mathcal{A}$ by Theorem \ref{t:additve-cat} and a well-known fact that the exchange graph $\mathcal{E}(\mathcal{C})$ is connected if it has a finite connected component.}
Now Conjecture \ref{c:monoidal-additive-correspondence} follows by Proposition \ref{p:equivalence-add-mult}.

% for any $\mathbf{g}\in \mathbb{Z}^n$, there exists %a morphism 
% $g\in \Hom_\mathcal{C}(T^{\mathbf{g}_-},T^{\mathbf{g}_+})$ such that $\operatorname{cone}(g)$ is rigid.  Consequently, $\mathfrak{e}(\mathbf{g})=0$. One the other hand, the cluster monomials form a basis of $\mathcal{A}$. In particular, $\chi_q(L(M))$ is a cluster monomial for any dominant monomial $M$. This implies that Conjecture \ref{c:monoidal-additive-correspondence} holds.

We now turn to Conjecture \ref{c:monoidal-additive-corr-compatible}. 
Let $M_1, M_2$ be dominant monomials in $\mathcal{P}_\ell^+$. 
Let $g_i\in \Hom_\mathcal{C}(T^{{\mathbf{g}_{M_i}}_-},T^{{\mathbf{g}_{M_i}}_+})$  such that $\operatorname{cone}(g_i)$ is rigid for $i=1,2$. Consequently, $\chi_q(L(M_i))=\mathbf{X}^{\hat T}_{\operatorname{cone}(g_i)\oplus P_i}$ for $i=1,2$, where $P_1$ and $P_2$ are projective objects of $\mathcal{F}$.

For the first implication, assume that 
$\chi_q(L(M_1))\chi_q(L(M_2))=\chi_q(L(M_1M_2))$. 
Since $L(M_1M_2)$ is real simple, 
$\chi_q(L(M_1M_2))$ is a cluster monomial and so the left hand side 
$\chi_q(L(M_1))\chi_q(L(M_2))=\mathbf{X}^{\hat T}_{\operatorname{cone}(g_1)\oplus P_1}\mathbf{X}^{\hat T}_{\operatorname{cone}(g_2)\oplus P_2}$ is a cluster monomial.
%This implies that $\mathbf{X}^T_{\operatorname{cone}(g_1)}\mathbf{X}^T_{\operatorname{cone}(g_2)}$ is a cluster monomial.
Therefore there exist a rigid object $N$ of $\mathcal{F}$ such that $\mathbf{X}^{\hat T}_{\operatorname{cone}(g_1)\oplus P_1}\mathbf{X}^{\hat T}_{\operatorname{cone}(g_2)\oplus P_2}=\mathbf{X}_N^{\hat T}$. By Corollary \ref{c:prod-cluster-monomial}, we conclude that $\operatorname{cone}(g_i)\in \add N$ for $i=1,2$. It follows that 
$\Ext_\mathcal{C}^1(\operatorname{cone}(g_1),\operatorname{cone}(g_2))=0$. Hence $\mathfrak{e}(\mathbf{g}_{M_1},\mathbf{g}_{M_2})=0$. 

For the other implication, assume that $\mathfrak{e}(\mathbf{g}_{M_1},\mathbf{g}_{M_2})=0$. This implies that \\
$\Ext_\mathcal{C}^1(\operatorname{cone}(g_1),\operatorname{cone}(g_2))=0$. By the $2$-Calabi-Yau property of $\mathcal{C}$, we conclude that $\operatorname{cone}(g_1)\oplus \operatorname{cone}(g_2)$ is rigid, hence $\mathbf{X}^{\hat T}_{\operatorname{cone}(g_1)\oplus P_1}\mathbf{X}^{\hat T}_{\operatorname{cone}(g_2)\oplus P_2}$ is a cluster monomial, which implies that $\chi_q(L(M_1))\chi_q(L(M_2))=\chi_q(L(M_1M_2))$. This completes the proof.
\end{proof}

%The following is a corollary of Conjecture \ref{c:monoidal-additive-correspondence}. Let ${\bf m}$ be a dominant monomial such that $L(\mathbf{m})$ is real. Then $L({\bf m})$ is prime if and only if $\mathfrak{e}(\mathbf{g}_{\mathbf{m}'},\mathbf{g}_{\mathbf{m}''}) \ne 0$ for any non-trivial dominant monomials ${\bf m}', {\bf m}''$.

\begin{proposition}\label{p:exchange-pair}
 Let $M, N$ be dominant monomials in $\mathcal{P}_\ell^+$ such that $\chi_q(L(M))$ and $\chi_q(L(N))$  are cluster variables of $\mathcal{A}$. Then $\chi_q(L(M))$ and $\chi_q(L(N))$ form an exchange pair if and only if $\mathfrak{e}(\mathbf{g}_{M},\mathbf{g}_{N})=1$.   
\end{proposition}
\begin{proof}
{We may assume that neither $\chi_q(L(M))$ nor $\chi_q(L(N))$  is a frozen variable.}
  Since $\chi_q(L(M))$ and $\chi_q(L(N))$ are mutable variables, there are $f\in \mathcal{O}_{\mathbf{g}_{M}}$ and $h\in \mathcal{O}_{\mathbf{g}_{N}}$ such that $\operatorname{cone}(f)$ and $\operatorname{cone}(h)$ are rigid. By Lemma \ref{l:ext-dim}, \[\mathfrak{e}(\mathbf{g}_{M},\mathbf{g}_{N})=\dim_\mathbb{k}\Ext^1_\mathcal{C}(\operatorname{cone}(f),\operatorname{cone}(h)).\]
  Now the result follows from \cite[Theorem 5.6 \& 5.8]{FGy19}.
\end{proof}

In general, we expect the following to be true: 
\begin{conjecture}\label{c:exchange-pair}
    Let $M, N$ be dominant monomials such that $L(M)$ and $L(N)$  are real prime. Then $\chi_q(L(M))\chi_q(L(N))=\chi_q(L(U))+\chi_q(L(V))$ for certain dominant monomials $U\neq V$ if and only if $\mathfrak{e}(\mathbf{g}_{M},\mathbf{g}_{N})=1$.
\end{conjecture}

\section{Evidence in the Grassmannian cluster categories of tame type} \label{sec:Evidences in Grassmannian cluster categories of tame type}

In this section, we focus on the Grassmannian cluster categories $\operatorname{CM}(B_{k,n})$. In particular, we establish the additive reachability conjecture for Grassmannian cluster categories of tame type, thus providing evidence for Conjecture~\ref{c:monoidal-additive-correspondence} and \ref{c:monoidal-additive-corr-compatible}.

\subsection{Grassmannian cluster categories} \label{subsec:grassmannian cluster categories}

\begin{figure}
\centering
    \begin{minipage}{.45\textwidth}
\begin{tikzpicture}[
    myedge/.style={thick,draw=black, postaction={decorate},
      decoration={markings,mark=at position .6 with {\arrow[black]{triangle 45}}}},
    myshorten/.style={shorten <= 2pt, shorten >= 2pt}]
    
    \node (A1) at (360/6*1:2cm) [circle,fill,inner sep=2pt, label=above:6] {}; 
    \node (A2) at (360/6*2:2cm) [circle,fill,inner sep=2pt, label=above:5] {}; 
    \node (A3) at (360/6*3:2cm) [circle,fill,inner sep=2pt, label=left:4] {};
    \node (A4) at (360/6*4:2cm) [circle,fill,inner sep=2pt, label=below:3] {}; 
    \node (A5) at (360/6*5:2cm) [circle,fill,inner sep=2pt, label=below:2] {}; 
    \node (A6) at (360/6*6:2cm) [circle,fill,inner sep=2pt, label=right:1] {}; 

  \draw (A1) to [bend right=25] node[midway,above] {6} (A2);
  \draw (A2) to [bend right=25] node[midway,left] {5} (A3);
  \draw (A3) to [bend right=25] node[midway,left] {4} (A4);
  \draw (A4) to [bend right=25] node[midway,below] {3} (A5);
  \draw (A5) to [bend right=25] node[midway,right] {2} (A6);
  \draw (A6) to [bend right=25] node[midway,right] {1} (A1);

\end{tikzpicture}
\end{minipage}
\centering
    \begin{minipage}{.45\textwidth}
\begin{tikzpicture}[
    myedge/.style={thick,draw=black, postaction={decorate},
      decoration={markings,mark=at position .6 with {\arrow[black]{triangle 45}}}},
    myshorten/.style={shorten <= 2pt, shorten >= 2pt}]
    
    \node (A1) at (360/6*1:2cm) [circle,fill,inner sep=2pt, label=above:6] {}; 
    \node (A2) at (360/6*2:2cm) [circle,fill,inner sep=2pt, label=above:5] {}; 
    \node (A3) at (360/6*3:2cm) [circle,fill,inner sep=2pt, label=left:4] {};
    \node (A4) at (360/6*4:2cm) [circle,fill,inner sep=2pt, label=below:3] {}; 
    \node (A5) at (360/6*5:2cm) [circle,fill,inner sep=2pt, label=below:2] {}; 
    \node (A6) at (360/6*6:2cm) [circle,fill,inner sep=2pt, label=right:1] {}; 
    
  \draw[<-] (A1) to [bend right=25] node[midway,above] {$x_6$} (A2);
  \draw[<-] (A2) to [bend right=25] node[midway,left] {$x_5$} (A3);
  \draw[<-] (A3) to [bend right=25] node[midway,left] {$x_4$} (A4);
  \draw[<-] (A4) to [bend right=25] node[midway,below] {$x_3$} (A5);
  \draw[<-] (A5) to [bend right=25] node[midway,right] {$x_2$} (A6);
  \draw[<-] (A6) to [bend right=25] node[midway,right] {$x_1$} (A1);
  
  \draw[<-] (A2) to [bend right=25] node[midway,below] {$y_6$} (A1);
  \draw[<-] (A3) to [bend right=25] node[midway,right] {$y_5$} (A2);
  \draw[<-] (A4) to [bend right=25] node[midway,right] {$y_4$} (A3);
  \draw[<-] (A5) to [bend right=25] node[midway,above] {$y_3$} (A4);
  \draw[<-] (A6) to [bend right=25] node[midway,left] {$y_2$} (A5);
  \draw[<-] (A1) to [bend right=25] node[midway,left] {$y_1$} (A6);

\end{tikzpicture}
\end{minipage}

\caption{The graph $C$ and the quiver $Q_C$, $n=6$.}
\label{fig:graph C and quiver Q_C}
\end{figure}
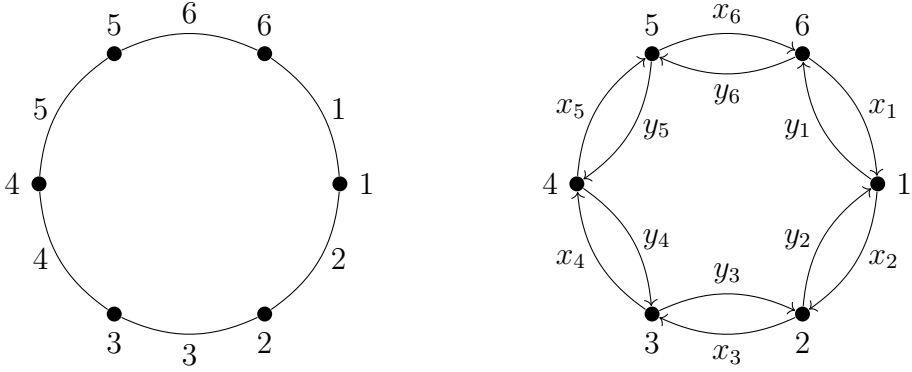

Denote by $C=(C_0, C_1)$ the circular graph with 
vertex set $C_0=\ZZ_{n}$ clockwise around the circle, and with the edge set $C_1=\ZZ_n$, with edge 
$i$ joining vertices $i-1$ and $i$, cf. Figure \ref{fig:graph C and quiver Q_C}. 
Denote by $Q_C$ the quiver with the same vertex set $C_0$ and with arrows 
$x_i: i-1 \to i$, $y_i: i \to i-1$ for every $i \in C_0$, cf. Figure \ref{fig:graph C and quiver Q_C}.

Denote by $B_{k,n}$ the quotient of the complete path algebra $\widehat{\CC Q_C}$ by the ideal generated 
by the $2n$ relations $x y = y x$, $x^{k} = y^{n-k}$, where $x, y$ are arrows of the form $x_i, y_j$ 
for appropriate $i,j$ (two relations for each vertex of $Q_C$). 

The center $Z$ of $B_{k,n}$ is the ring of formal power series $\CC[[t]]$, where $t = \sum_{i=1}^n x_i y_i$. A (maximal) Cohen-Macaulay $B$-module is given by a representation $\{M_i: i \in C_0\}$ of $Q_C$, where each $M_i$ is a free $Z$-module of the same rank, see \cite[Section 3]{JKS16}.
For any $B_{k,n}$-module $M$ and $K$ the field of fractions of $Z$, the rank of $M$, denoted by ${\rm rk}(M)$, is defined to be ${\rm rk}(M) = {\rm len}(M \otimes_Z K)$, see \cite[Definition 3.5]{JKS16}. Jensen, King, and Su \cite{JKS16} proved that the category ${\rm CM}(B_{k,n})$ 
is an additive categorification of the cluster algebra structure on $\CC[\Gr(k,n)]$. The category ${\rm CM}(B_{k,n})$ is exact and Frobenius with projective-injective objects given by the $B_{k,n}$ projective modules, and it has an Auslander-Reiten quiver (\cite[Remark 3.3]{JKS16}). We 
denote by $\tau(M)$ the Auslander-Reiten translation of $M$ and by $\tau^{-1}(M)$ the inverse Auslander-Reiten translation of $M$. 

A module $M$ in ${\rm CM}(B_{k,n})$ is rigid if $ \Ext^1_{{\rm CM}(B_{k,n})} (M,M)=0 $. A special class of objects of ${\rm C}(B_{k,n})$ are the rank 1 modules which are known to be rigid, \cite[Proposition 5.6]{JKS16}. For any $k$-subset $I$ of $C_1$, we define a rank 1 module $L_I$ in ${\rm CM}(B_{k,n})$ 
as follows: 
\begin{align*}
L_I=(U_i, i \in C_0; x_i, y_i, i \in C_1),
\end{align*}
where $U_i = \CC[[t]]$, $i \in C_0$, and 
\begin{align*}
x_i: U_{i-1} \to U_i \text{ is given by multiplication by $1$ if $i \in I$, and by $t$ if $i \not\in I$,} \\
y_i: U_{i} \to U_{i-1} \text{ is given by multiplication by $t$ if $i \in I$, and by $1$ if $i \not\in I$,}
\end{align*}
see \cite[Definition 5.1]{JKS16} 

By \cite[Proposition 5.2]{JKS16}, every rank 1 module is isomorphic to $L_I$ for some 
$k$-subset $I$ of $[n]$. So there are bijections between the rank 1 modules in ${\rm CM}(B_{k,n})$, the $k$-subsets of $[n]$ and the cluster variables of $\CC[\Gr(k,n)]$ which are Pl\"ucker coordinates. 

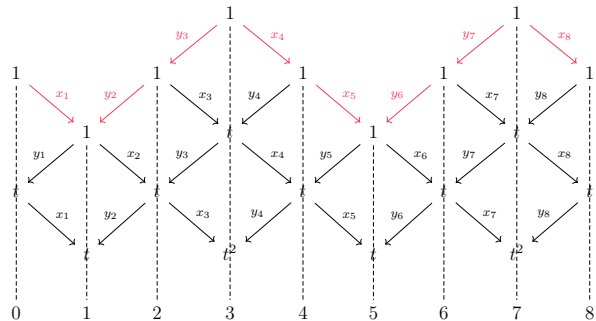
\begin{figure}
\centering
\adjustbox{scale=0.6, center}{
\begin{tikzcd}
	&&& 1 &&&& 1 \\
	1 && 1 && 1 && 1 && 1 \\
	& 1 && t && 1 && t \\
	t && t && t && t && t \\
	& t && {t^2} && t && {t^2} \\
	0 & 1 & 2 & 3 & 4 & 5 & 6 & 7 & 8
	\arrow["{x_1}", color={rgb,255:red,244;green,62;blue,86}, from=2-1, to=3-2]
	\arrow["{y_3}"', color={rgb,255:red,244;green,62;blue,86}, from=1-4, to=2-3]
	\arrow["{y_2}"', color={rgb,255:red,244;green,62;blue,86}, from=2-3, to=3-2]
	\arrow["{x_4}", color={rgb,255:red,244;green,62;blue,86}, from=1-4, to=2-5]
	\arrow["{x_5}", color={rgb,255:red,244;green,62;blue,86}, from=2-5, to=3-6]
	\arrow["{y_6}"', color={rgb,255:red,244;green,62;blue,86}, from=2-7, to=3-6]
	\arrow["{y_7}"', color={rgb,255:red,244;green,62;blue,86}, from=1-8, to=2-7]
	\arrow["{x_8}", color={rgb,255:red,244;green,62;blue,86}, from=1-8, to=2-9]
	\arrow["{y_1}"', from=3-2, to=4-1]
	\arrow[dashed, no head, from=2-1, to=6-1]
	\arrow["{x_2}", from=3-2, to=4-3]
	\arrow["{x_3}", from=2-3, to=3-4]
	\arrow["{x_4}", from=3-4, to=4-5]
	\arrow["{x_3}", from=4-3, to=5-4]
	\arrow["{x_6}", from=3-6, to=4-7]
	\arrow["{x_5}", from=4-5, to=5-6]
	\arrow["{x_7}", from=4-7, to=5-8]
	\arrow["{x_7}", from=2-7, to=3-8]
	\arrow["{x_8}", from=3-8, to=4-9]
	\arrow["{y_8}"', from=4-9, to=5-8]
	\arrow[dashed, no head, from=2-9, to=6-9]
	\arrow["{y_8}"', from=2-9, to=3-8]
	\arrow["{y_7}"', from=3-8, to=4-7]
	\arrow["{y_6}"', from=4-7, to=5-6]
	\arrow["{y_5}"', from=3-6, to=4-5]
	\arrow["{y_4}"', from=4-5, to=5-4]
	\arrow["{y_4}"', from=2-5, to=3-4]
	\arrow["{y_3}"', from=3-4, to=4-3]
	\arrow["{y_2}"', from=4-3, to=5-2]
	\arrow["{x_1}", from=4-1, to=5-2]
	\arrow[dashed, no head, from=3-2, to=6-2]
	\arrow[dashed, no head, from=2-3, to=6-3]
	\arrow[dashed, no head, from=1-4, to=6-4]
	\arrow[dashed, no head, from=2-5, to=6-5]
	\arrow[dashed, no head, from=3-6, to=6-6]
	\arrow[dashed, no head, from=2-7, to=6-7]
	\arrow[dashed, no head, from=1-8, to=6-8]
\end{tikzcd}
}
\caption{The lattice diagram of $L_{\{1,4,5,8\}}$ in ${\rm CM}(B_{4,8})$ 
with its rim indicated by the red arrows.}
\label{fig:lattice diagram of L1458 in B48}
\end{figure}

It is convenient to represent the module $L_I$ by a lattice diagram, see 
Figure~\ref{fig:lattice diagram of L1458 in B48}. The spaces $U_0, \ldots, U_n$ are represented by 
columns from left to right and $U_0$ and $U_n$ are identified. The vertices in each column correspond 
to the natural monomial $\CC$-basis of $\CC[t]$. The column corresponding to $U_{i+1}$ is displaced half a step vertically downwards (resp. upwards) in 
relation to $U_i$ if $i + 1 \in I$ (resp. $i+1 \not\in I$). The upper boundary of the lattice diagram of $L_I$ is called the rim of $L_I$. The $k$-subset $I \subset [n]$ of the rank $1$ module $L_I$ can be read off as the set of labels on the edges going down to the right which are on the rim of $L_I$, i.e. the labels of the $x_i$'s appearing in the rim.  

The rank 1 modules can be viewed as building blocks for the category as every module in ${\rm CM}(B_{k,n})$ has a filtration with factors which are rank $1$ modules (cf. \cite[Proposition 6.6]{JKS16}). {We refer to such a filtration as a {\em generic filtration}, since it depends on a generic choice of a full flag in the limiting vector space of the module.} Let $M$ be a rank $m$ module in ${\rm CM}(B_{k,n})$ with factors $L_{I_1},\dots, L_{I_m}$ in its generic filtration, where $L_{I_m}$ is a submodule of $M$. We write $M=L_{I_1}|L_{I_2}| \cdots |L_{I_m}$ or 
$M=\scalemath{0.86}{\begin{array}{c} L_{I_1}\\ \hline \vdots \\ \hline L_{I_m} \end{array}}$\,. The ordered 
collection of $k$-subsets $I_1, \ldots, I_m$ in the generic filtration of $M$ is called the profile of $M$, denoted $P_M$. We write 
$P_M= \scalemath{0.86}{\begin{array}{c} I_1\\ \hline \vdots \\ \hline I_m \end{array}}$ 
or $P_M=I_1| \cdots |I_m$ if $M$ has a filtration having factors $L_{I_1}, \ldots, L_{I_m}$ (in this order). 
We sometimes write $M=P_M$ to indicate that $M$ is a module with profile $P_M$. Note that in general, such a filtration is not unique, but in case $M$ is rigid, the filtration is unique in the sense that it gives a canonical ordered set of rank 1 composition factors (\cite[Lemma 6.2]{JKS16}).

%In the case of $\Gr(2,n)$ or $(k,n) \in \{(3,6), (3,7), (3,8), (4,6), (4,7)\}$, the cluster algebra $\CC[\Gr(k,n)]$ is of finite type. The cluster algebras $\CC[\Gr(3,9)]$ and $\CC[\Gr(4,8)]$ are of tame type. In this section, we ...

%\begin{definition}
    %Let $\mathcal{E}_{k,n}$ be the exchange graph of basic cluster-tilting objects of $\operatorname{CM}(B_{k,n})$ and $\underline{\mathcal{E}}_{k,n}$ the exchange graph of basic cluster-tilting objects of $\underline{\operatorname{CM}}(B_{k,n})$.
%\end{definition}
%\begin{lemma}
    %There is a one-to-one correspondence between the set of basic cluster-tilting objects of $\operatorname{CM}(B_{k,n})$ and the set of basic cluster-tilting objects of $\underline{\operatorname{CM}}(B_{k,n})$. Moreover, the bijection is compatible with mutations. As a consequence, $\mathcal{E}_{k,n}\cong \underline{\mathcal{E}}_{k,n}$.
%\end{lemma}

%%
%
\subsection{Connectivity of the exchange graphs $\mathcal{E}(\operatorname{CM}(B_{3,9}))$ and $\mathcal{E}(\operatorname{CM}(B_{4,8}))$.} 

If we assume that $k\le \frac{n}{2}$, then the category $\operatorname{CM}(B_{k,n})$ is of finite type if and only if either $k=2$ or $(k,n)\in \{(3,6), (3,7), (3,8)\}$  and that it is of tame type if and only if $(k,n)\in \{(3,9),(4,8)\}$, see~\cite[Section 2]{JKS16}.

%It is well-known that  $\operatorname{CM}(B_{k,n})$ is of finite type if and only if  either $k=2$ or $(k,n)\in \{(3,6), (3,7), (3,8), (4,6), (4,7)\}$, and it is of tame type if and only if $(k,n)\in \{(3,9),(4,8)\}$. 
In the finite types, the exchange graph $\mathcal{E}(\operatorname{CM}(B_{k,n}))$ is always connected. 
In this Section, we show that we also have connectedness in the tame types. 

We begin with an observation, which is a direct consequence of \cite[Theorem 4.1]{AIR14}.
\begin{lemma}\label{l:tau-tilting}
 Let $T$ be a basic cluster-tilting object in a $2$-Calabi-Yau triangulated category $\mathcal{C}$  and $\End_\mathcal{C}(T)$ the endomorphism algebra of $T$. Then the exchange graph of support $\tau$-tilting modules of $\End_\mathcal{C}(T)$ is isomorphic to the exchange graph $\mathcal{E}(\mathcal{C})$.
\end{lemma}

We recall a result of~\cite{BKM16} which gives a combinatorial construction of (endomorphism algebras of) cluster-tilting objects in $\operatorname{CM}(B_{k,n})$. 

Let $D_{k,n}$ be a $(k,n)$-Postnikov diagram and $(\widehat{Q}(D_{k,n}), \widehat{W}(D_{k,n}))$ the associated quiver with potential (cf. \cite[Section 2]{BKM16}). 
The quiver $\widehat{Q}(D_{k,n})$ is a so-called dimer quiver: it is a quiver with faces. {Its faces are given by alternating positive and negative oriented minimal cycles (cycles which do not contain any shorter cycles), such that each internal arrow belongs to exactly two faces.}
Furthermore, the quiver is 
embedded in a disk with $n$ vertices on the boundary.
An example for $(k,n)=(5,9)$ is in Figure~\ref{fig:5-9-quiver}. 
The quiver {$\widehat{Q}(D_{k,n})$} comes with a natural potential 
$\widehat{W}(D_{k,n})$ which is the sum of all the positive unit cycles minus the sum of all the negative unit cycles.

The Jacobian algebra $J(\widehat{Q}(D_{k,n}),\widehat{W}(D_{k,n}))$ of such a quiver is the completed path algebra with relations given by the cyclic derivatives of the potential.
It has been shown in~\cite[Theorem 10.3]{BKM16} that the Jacobian algebra $J(\widehat{Q}(D_{k,n}),\widehat{W}(D_{k,n}))$ is isomorphic to the endomorphism algebra $\End_{B_{k,n}}(\widehat{T}_{k,n})$ for a basic cluster-tilting object $\widehat{T}_{k,n}$ of $\operatorname{CM}(B_{k,n})$ which is given explicitly from the Postnikov diagram $D_{k,n}$ (and is a direct sum of rank 1 modules). 

Let $({Q}(D_{k,n}),{W}(D_{k,n}))$ be the  quiver with potential we obtain from the original quiver $(\widehat{Q}(D_{k,n}),\widehat{W}(D_{k,n}))$ by deleting {the} vertices on the boundary and the arrows incident with them. 
Examples of the dimer quivers {$Q(D_{k,n})$} are given in Section~\ref{subsec:non-rigid-module-example-Gr39} (for $(3,9)$) and~\ref{subsec:non-rigid-module-example-Gr48} (for $(4,8)$).  
Denote the (basic) cluster-tilting object of $\underline{\operatorname{CM}}(B_{k,n})$ obtained from $\widehat{T}_{k,n}$ by removing the projective-injective summands by $T_{k,n}$ (its summands correspond to the non-boundary vertices of {$\widehat{Q}(D_{k,n})$}). 
We then have $J({Q}(D_{k,n}),{W}(D_{k,n}))\cong \End({T}_{k,n})$. One has: 
%The following was proved in \cite{CZ22}.
\begin{theorem}\cite[Theorem 0.2]{CZ22}\label{t:uniqueness-deg-potential}
    The potential ${W}(D_{k,n})$ (resp. $\widehat{W}(D_{k,n})$) is the unique nondegenerate potential on ${Q}(D_{k,n})$ (resp, $\widehat{Q}(D_{k,n})$) up to right equivalence.
\end{theorem}

Now we are ready to prove the connectedness of the exchange graphs in the tame types. 

\begin{theorem} \label{thm:CMGr39-48-connected}
%\label{thm:CMGr39 amd CMGr48 are connected}
The exchange graphs $\mathcal{E}(\operatorname{CM}(B_{3,9}))$ and $\mathcal{E}(\operatorname{CM}(B_{4,8}))$ are connected.
\end{theorem}

\begin{proof}
By Lemma~\ref{l:isom-exchange-graph}, it suffices to prove the statement for the stable categories, i.e. to show that $\mathcal{E}(\underline{\operatorname{CM}}(B_{3,9}))$ and $\mathcal{E}(\underline{\operatorname{CM}}(B_{4,8}))$ are connected. 

Let $D_{k,n}$ be a $(k,n)$-Postnikov diagram with quiver $Q(D_{k,n})$. We first observe that $Q(D_{4,8})$ is mutation-equivalent to the quiver $Q_7^{(1,1)}$ from Figure~\ref{f:quiver-tubular-type} and that  
${Q}(D_{3,9})$ is mutation-equivalent to the quiver $Q^{(1,1)}_8$. This is checked using Keller's quiver mutation Java applet \cite{Keller_mutation}. 
By~\cite[Lemma 4.1]{FG21}, the quiver $Q_7^{(1,1)}$ (resp. $Q_8^{(1,1)}$) admits a unique nondegenerate potential which we denote by $W_7$ (resp. by $W_8$). 
 
Therefore, we get that 
$(Q_7^{(1,1)},W_7)$ (resp. $(Q_8^{(1,1)},W_8)$) is mutation-equivalent to $({Q}(D_{4,8}),{W}(D_{4,8}))$ (resp. $({Q}(D_{3,9}),{W}(D_{3,9}))$) by Theorem \ref{t:uniqueness-deg-potential}.

Denote by $\mathcal{C}(Q_7^{(1,1)},W_7)$ and $\mathcal{C}(Q_8^{(1,1)},W_8)$ the generalized cluster categories associated with 
these quivers with potential %$(Q_7^{(1,1)},W_7)$ and $(Q_8^{(1,1)},W_8)$ respectively 
in the sense of \cite{Amiot09}. 
These two categories are triangle equivalent to cluster categories of weighted projective lines: 
let $\mathbb{X}(p_1,p_2,p_3)$ be the weighted projective line with weights $(p_1,p_2,p_3)$ and $\mathcal{C}_{\mathbb{X}(p_1,p_2,p_3)}$ the associated cluster category. 
Then we have $\mathcal{C}(Q_7^{(1,1)},W_7)\cong \mathcal{C}_{\mathbb{X}(2,4,4)}$ and $\mathcal{C}(Q_8^{(1,1)},W_8)\cong \mathcal{C}_{\mathbb{X}(2,3,6)}$ (cf. \cite[Section 4]{FG21}). 
It follows that $\End({T}_{4,8})$ (resp. $\End({T}_{3,9})$) is isomorphic to the endomorphism algebra of a basic cluster-tilting object of $\mathcal{C}_{\mathbb{X}(2,4,4)}$ (resp. $\mathcal{C}_{\mathbb{X}(2,3,6)}$). Since the exchange graphs of $\mathcal{C}_{\mathbb{X}(2,4,4)}$  and $\mathcal{C}_{\mathbb{X}(2,3,6)}$ are connected by~\cite[Section 8]{BKL10}, we get the desired result using Lemma~\ref{l:tau-tilting}. %the exchange graphs of $\mathcal{C}_{\mathbb{X}(2,4,4)}$  and $\mathcal{C}_{\mathbb{X}(2,3,6)}$ are connected. We conclude the desired result by Lemma \ref{l:tau-tilting}.

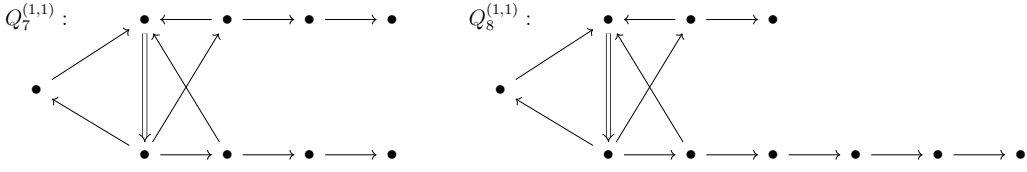
\begin{figure}
    \centering
\adjustbox{scale=0.7}{ \begin{tikzcd}
	{Q_7^{(1,1)}:} & \bullet & \bullet & \bullet & \bullet & {Q_8^{(1,1)}:} & \bullet & \bullet & \bullet \\
	\bullet &&&&& \bullet \\
	& \bullet & \bullet & \bullet & \bullet && \bullet & \bullet & \bullet & \bullet & \bullet & \bullet
	\arrow[from=2-1, to=1-2]
	\arrow[from=1-3, to=1-2]
	\arrow[from=1-3, to=1-4]
	\arrow[from=1-4, to=1-5]
	\arrow[from=3-2, to=2-1]
	\arrow[Rightarrow, from=1-2, to=3-2]
	\arrow[from=3-2, to=3-3]
	\arrow[from=3-3, to=3-4]
	\arrow[from=3-4, to=3-5]
	\arrow[from=3-3, to=1-2]
	\arrow[from=3-2, to=1-3]
	\arrow[from=3-8, to=3-9]
	\arrow[from=3-9, to=3-10]
	\arrow[from=3-10, to=3-11]
	\arrow[from=3-11, to=3-12]
	\arrow[from=2-6, to=1-7]
	\arrow[from=3-7, to=2-6]
	\arrow[Rightarrow, from=1-7, to=3-7]
	\arrow[from=1-8, to=1-7]
	\arrow[from=1-8, to=1-9]
	\arrow[from=3-7, to=3-8]
	\arrow[from=3-8, to=1-7]
	\arrow[from=3-7, to=1-8]
\end{tikzcd}}
    \caption{Quivers of the tubular types $E_7^{(1,1)}$ and $E_8^{(1,1)}$.}
    \label{f:quiver-tubular-type}
\end{figure}

\end{proof}

We conclude Section~\ref{sec:Evidences in Grassmannian cluster categories of tame type} by providing examples for non-real $\mathbf{g}$-vectors for the two tame cases, thus giving evidence for Conjecture~\ref{c:monoidal-additive-corr-compatible}. 
\subsection{An example of non-real $\mathbf{g}$-vector for $\Gr(3,9)$} \label{subsec:non-rigid-module-example-Gr39}

Consider 
\[
\widehat{T}=L_{124} \oplus L_{125}\oplus L_{126}\oplus L_{127}
\oplus L_{128}\oplus L_{134} \oplus L_{145} \oplus L_{156} \oplus L_{167} \oplus L_{178} \oplus(\oplus_{i=1}^9 L_{i,i+1,i+2})
\]
(where we reduce subscripts modulo $9$). The nine summands $L_{i,i+1,i+2}$ are projective-injective. The module $\widehat{T}$ is a basic cluster-tilting object of 
$\text{CM}(B_{3,9})$. 
%Let $\widehat{T}=\boxed{123}\oplus 124\oplus 125\oplus 126\oplus 127\oplus 128\oplus \boxed{129}\oplus 134\oplus 145\oplus 156\oplus 167\oplus 178\oplus 189\oplus \boxed{234}\oplus \boxed{345}\oplus \boxed{456}\oplus \boxed{567}\oplus \boxed{678}\oplus \boxed{789}$ be a basic cluster-tilting object of $\text{CM}(B_{3,9})$, where $\boxed{*}$ means that $*$ is projective-injective. 
We denote the associated (basic) cluster-tilting object of $\underline{\text{CM}}(B_{3,9})$ by $T$. The endomorphism algebra $A:=\End_{\underline{\text{CM}}(B_{3,9})}(T)$ is isomorphic to the opposite of the Jacobian algebra of the following quiver 
(to abbreviate notation, we write $I$ instead of $L_I$ for the vertices of the quiver). 
    \[\adjustbox{scale=0.7}
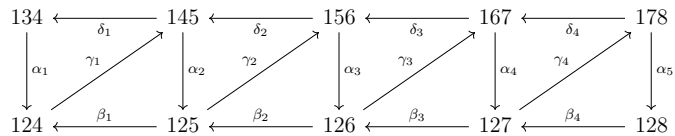
{\begin{tikzcd}
	134 && 145 && 156 && 167 && 178 \\
	\\
	124 && 125 && 126 && 127 && 128
	\arrow["{\alpha_1}", from=1-1, to=3-1]
	\arrow["{\alpha_2}", from=1-3, to=3-3]
	\arrow["{\alpha_4}", from=1-7, to=3-7]
	\arrow["{\alpha_3}", from=1-5, to=3-5]
	\arrow["{\alpha_5}", from=1-9, to=3-9]
	\arrow["{\beta_1}"',   from=3-3, to=3-1]
	\arrow["{\beta_2}"', from=3-5, to=3-3]
	\arrow["{\beta_3}"', from=3-7, to=3-5]
	\arrow["{\beta_4}"', from=3-9, to=3-7]
	\arrow["{\delta_1}", from=1-3, to=1-1]
	\arrow["{\delta_2}", from=1-5, to=1-3]
	\arrow["{\delta_3}", from=1-7, to=1-5]
	\arrow["{\delta_4}", from=1-9, to=1-7]
	\arrow["{\gamma_1}", from=3-1, to=1-3]
	\arrow["{\gamma_2}", from=3-3, to=1-5]
	\arrow["{\gamma_3}", from=3-5, to=1-7]
	\arrow["{\gamma_4}", from=3-7, to=1-9]
\end{tikzcd}}\]
with potential $W=\sum_{i=1}^4\beta_i\gamma_i\alpha_{i+1}-\sum_{i=1}^4\delta_i\alpha_i\gamma_i$. 

Note that right $A$-modules identify with representations over the quiver with relations. In the following, we work with representations of the quiver with relations.

We consider the $\mathbf{g}$-vector 
$\mathbf{g}=
\scalemath{0.6}{\begin{pmatrix}
    -1&0&1&1&0\\ 0&-1&-1&0&1
\end{pmatrix}}\in \mathbb{Z}^{10}$, writing it according to the shape of the quiver of $A$. 

One can check that this $\mathbf{g}$-vector corresponds to the tableau (see Remark~\ref{rem:g-vector_to_tableau}) 
\[
\ytableausetup{centertableaux}
\scalemath{0.6}{\begin{ytableau}
1 & 2 & 3 \\
4 & 5 & 6 \\
7 & 8 & 9
\end{ytableau}}.
\]
One can check that the dominant monomial associated to the tableau/the $\mathbf{g}$-vector is not real, see~\cite[Section 8]{CDFL20}.

We expect that no positive multiple of $\mathbf{g}$ corresponds to a real simple module. To provide evidence for Conjecture~\ref{c:monoidal-additive-correspondence}, we  show that 
$\mathfrak{e}(t\mathbf{g})>0$ for any positive integer $t$. 

We have 
$T^{\mathbf{g}_-}=P(125)\oplus P(126)\oplus P(134)$ and 
$T^{\mathbf{g}_+}=P(128)\oplus P(156)\oplus P(167)$. 
The indecomposable projective $A$-modules appearing in these modules are the following: 
\[\adjustbox{scale=0.76}{ \begin{tikzcd}
	{P{(125):}} & 0 & \mathbb{k} & \mathbb{k} & 0 & 0 & {P{(128):}} & 0 & 0 & 0 & 0 & 0 \\
	& \mathbb{k} & \mathbb{k} & 0 & 0 & 0 && \mathbb{k} & \mathbb{k} & \mathbb{k} & \mathbb{k} & \mathbb{k}\\
	{P{(126):}} & 0 & \mathbb{k} & \mathbb{k} & \mathbb{k} & 0 & {P{(156):}} & \mathbb{k} & \mathbb{k} & \mathbb{k}& 0 & 0 \\
	& \mathbb{k}& \mathbb{k} & \mathbb{k} & 0 & 0 && \mathbb{k} & \mathbb{k} & \mathbb{k} & 0 & 0 \\
	{P{(134):}} & \mathbb{k} & 0 & 0 & 0 & 0 & {P{(167):}} & \mathbb{k} & \mathbb{k} & \mathbb{k}& \mathbb{k} & 0 \\
	& \mathbb{k}& 0 & 0 & 0 & 0 && \mathbb{k} & \mathbb{k} & \mathbb{k} & \mathbb{k} & 0
	\arrow[from=1-2, to=2-2]
	\arrow["0"{description}, from=1-3, to=2-3]
	\arrow[from=1-4, to=2-4]
	\arrow[from=1-5, to=2-5]
	\arrow[from=1-6, to=2-6]
	\arrow[from=1-3, to=1-2]
	\arrow["1"{description}, from=1-4, to=1-3]
	\arrow[from=1-5, to=1-4]
	\arrow[from=1-6, to=1-5]
	\arrow["1"{description}, from=2-3, to=2-2]
	\arrow[from=2-4, to=2-3]
	\arrow[from=2-5, to=2-4]
	\arrow["1"{description}, from=2-2, to=1-3]
	\arrow["1"{description}, from=2-3, to=1-4]
	\arrow[from=2-4, to=1-5]
	\arrow[from=2-5, to=1-6]
	\arrow[from=2-6, to=2-5]
	\arrow[from=3-2, to=4-2]
	\arrow["1"{description}, from=4-3, to=4-2]
	\arrow["1"{description}, from=4-2, to=3-3]
	\arrow[from=3-3, to=3-2]
	\arrow["1"{description}, from=3-4, to=3-3]
	\arrow["1"{description}, from=3-5, to=3-4]
	\arrow["1"{description}, from=4-4, to=4-3]
	\arrow[from=4-5, to=4-4]
	\arrow[from=3-6, to=3-5]
	\arrow[from=3-5, to=4-5]
	\arrow[from=3-6, to=4-6]
	\arrow["0"{description}, from=3-4, to=4-4]
	\arrow["0"{description}, from=3-3, to=4-3]
	\arrow["1"{description}, from=4-3, to=3-4]
	\arrow["1"{description}, from=4-4, to=3-5]
	\arrow[from=4-5, to=3-6]
	\arrow[from=4-6, to=4-5]
	\arrow["1"{description}, from=5-2, to=6-2]
	\arrow[from=5-3, to=6-3]
	\arrow[from=5-4, to=6-4]
	\arrow[from=5-5, to=6-5]
	\arrow[from=5-6, to=6-6]
	\arrow[from=6-6, to=6-5]
	\arrow[from=6-5, to=6-4]
	\arrow[from=6-4, to=6-3]
	\arrow[from=6-3, to=6-2]
	\arrow[from=5-3, to=5-2]
	\arrow[from=5-4, to=5-3]
	\arrow[from=5-5, to=5-4]
	\arrow[from=5-6, to=5-5]
	\arrow[from=6-2, to=5-3]
	\arrow[from=6-3, to=5-4]
	\arrow[from=6-4, to=5-5]
	\arrow[from=6-5, to=5-6]
	\arrow[from=1-8, to=2-8]
	\arrow[from=1-9, to=2-9]
	\arrow[from=1-10, to=2-10]
	\arrow[from=1-11, to=2-11]
	\arrow[from=1-12, to=2-12]
	\arrow[from=1-9, to=1-8]
	\arrow[from=1-10, to=1-9]
	\arrow[from=1-11, to=1-10]
	\arrow[from=1-12, to=1-11]
	\arrow["1"{description}, from=2-9, to=2-8]
	\arrow["1"{description}, from=2-10, to=2-9]
	\arrow["1"{description}, from=2-11, to=2-10]
	\arrow["1"{description}, from=2-12, to=2-11]
	\arrow[from=2-8, to=1-9]
	\arrow[from=2-9, to=1-10]
	\arrow[from=2-10, to=1-11]
	\arrow[from=2-11, to=1-12]
	\arrow["1"{description}, from=3-8, to=4-8]
	\arrow["1"{description}, from=3-9, to=4-9]
	\arrow["1"{description}, from=3-10, to=4-10]
	\arrow[from=3-11, to=4-11]
	\arrow[from=3-12, to=4-12]
	\arrow[from=4-12, to=4-11]
	\arrow[from=4-11, to=4-10]
	\arrow["1"{description}, from=4-10, to=4-9]
	\arrow["1"{description}, from=4-9, to=4-8]
	\arrow["1"{description}, from=3-9, to=3-8]
	\arrow["1"{description}, from=3-10, to=3-9]
	\arrow[from=3-11, to=3-10]
	\arrow[from=3-12, to=3-11]
	\arrow["0"{description}, from=4-8, to=3-9]
	\arrow["0"{description}, from=4-9, to=3-10]
	\arrow[from=4-10, to=3-11]
	\arrow[from=4-11, to=3-12]
	\arrow["1"{description}, from=5-8, to=6-8]
	\arrow["1"{description}, from=5-9, to=6-9]
	\arrow["1"{description}, from=5-10, to=6-10]
	\arrow["1"{description}, from=5-11, to=6-11]
	\arrow[from=5-12, to=6-12]
	\arrow["1"{description}, from=5-9, to=5-8]
	\arrow["1"{description}, from=5-10, to=5-9]
	\arrow["1"{description}, from=5-11, to=5-10]
	\arrow["0"{description}, from=5-12, to=5-11]
	\arrow["1"{description}, from=6-9, to=6-8]
	\arrow["1"{description}, from=6-10, to=6-9]
	\arrow["1"{description}, from=6-11, to=6-10]
	\arrow[from=6-12, to=6-11]
	\arrow["0"{description}, from=6-8, to=5-9]
	\arrow["0"{description}, from=6-9, to=5-10]
	\arrow["0"{description}, from=6-10, to=5-11]
	\arrow[from=6-11, to=5-12]
\end{tikzcd}}\]
The $\Hom$-spaces between these modules are listed in Table~\ref{ta:hom-space-39}. 
%$\Hom(P(i),P(j))$ are listed  by the following table:
\begin{table}[ht]

\begin{tabular}{|c|c|c|c|c|c|c|}
\hline
     & $P{(125)}$&$P{(126)}$&$P{(134)}$&$P{(128)}$&$P{(156)}$&$P{(167)}$  \\
     \hline
     $P{(125)}$&$\mathbb{k}$&$\mathbb{k}$&0&$\mathbb{k}$&$\mathbb{k}$&$\mathbb{k}$ \\
     \hline
     $P{(126)}$&0&$\mathbb{k}$&0&$\mathbb{k}$&$\mathbb{k}$&$\mathbb{k}$\\ \hline
     $P{(134)}$&0&0&$\mathbb{k}$&0&$\mathbb{k}$&$\mathbb{k}$\\ \hline
     $P{(128)}$&0&0&0&$\mathbb{k}$&0&0\\ \hline
     $P{(156)}$&$\mathbb{k}$&$\mathbb{k}$&0&0&$\mathbb{k}$&$\mathbb{k}$\\ \hline
     $P{(167)}$&0&$\mathbb{k}$&0&0&0&$\mathbb{k}$\\\hline
\end{tabular}
\caption{$\Hom$-spaces for 
$\text{CM}(B_{3,9})$}\label{ta:hom-space-39}
\end{table}

\begin{proposition}\label{p:non-real-39}
    We have  $\mathfrak{e}(t\mathbf{g})>0$ for any positive integer $t$. 
\end{proposition}
\begin{proof}
We first consider $t=1$ and show $\mathfrak{e}(\mathbf{g})>0$.

Let $B\in \Hom(T^{\mathbf{g}_-},T^{\mathbf{g}_+})$ be arbitrary. We show that there exists $X\in \Hom(T^{\mathbf{g}_-},T^{\mathbf{g}_+})$ such that 
\[
X\neq B\circ U+V\circ B.
\]
for any $U\in \End(T^{\mathbf{g}_-})$ and for any $V\in \End(T^{\mathbf{g}_+})$. 
%
%Let $B\in \Hom(P{(125)}\oplus P{(126)}\oplus P{(134)}, P{(128)}\oplus P{(156)}\oplus P{(167)})$ be an arbitrary morphism. We are going to show that there is a morphism $X\in \Hom(P{(125)}\oplus P{(126)}\oplus P{(134)}, P{(128)}\oplus P{(156)}\oplus P{(167)})$ such that for any morphism 
%
%$U\in \End(P{(125)}\oplus P{(126)}\oplus P{(134)})$ and 
%$V\in \Hom(P{(128)}\oplus P{(156)}\oplus P{(167)},P{(128)}\oplus P{(156)}\oplus P{(167)})$, 
%$V\in \End(P{(128)}\oplus P{(156)}\oplus P{(167)})$, 
%we have 
%\[
%X\neq B\circ U+V\circ B.
%\]
If such a morphism $X$ exists, it is non-zero in 
%In other words, $X$ is non-zero in 
$\Hom_{K^b(\add A)}(P_B^\bullet, \Sigma P_B^\bullet)$, where $P_B^\bullet$ is the $2$-term complex 
$T^{\mathbf{g}_-}\xrightarrow{B}T^{\mathbf{g}_+}$ and 
hence $\mathfrak{e}(\mathbf{g})>0$.

%$P{(125)}\oplus P{(126)}\oplus P{(134)}\xrightarrow{B}P{(128)}\oplus P{(156)}\oplus P{(167)}$ and hence $\mathfrak{e}(\mathbf{g})>0$. 

To see that $X$ exists, we write the morphisms $X,B,U,V$ as $3\times 3$ matrices over $\mathbb{k}$, we can do this using Table~\ref{ta:hom-space-39}: 
%According to Table \ref{ta:hom-space-39}, the morphisms $X,B,U,V$ can all be represented by $3\times 3$ matrices over $\mathbb{k}$: 
\[
\scalemath{0.96}{
X=\begin{bmatrix}
    x_{11}&x_{12}&0\\ x_{21}&x_{22}&x_{23}\\ x_{31}&x_{32}&x_{33}
\end{bmatrix},\,
B=\begin{bmatrix}
    b_{11}&b_{12}&0\\ b_{21}&b_{22}&b_{23}\\ b_{31}&b_{32}&b_{33}
\end{bmatrix}, \,
U=\begin{bmatrix}
 u_{11}&0&0\\u_{21}&u_{22}&0\\ 0&0&u_{33}
\end{bmatrix},\,
V=\begin{bmatrix}
    v_{11}&0&0\\ 0&v_{22}&0\\ 0&v_{32}&v_{33}
\end{bmatrix}.
}
\] 
Regarding $u_{ij}$ and $v_{ij}$ as indeterminates, we can find $X$ such that $BU+VB=X$ has no solution by the fact that the determinant of the coefficient matrix of this equation is always zero.

Now suppose that $\mathfrak{e}(t\mathbf{g})=0$ for some positive integer $t$. Then there is a morphism $f\in\Hom(P^{t\mathbf{g}_-},P^{t\mathbf{g}_+})$ such that $\Hom_{K^b(\add A)}(P_f^\bullet,\Sigma P_f^\bullet)=0$. Equivalently, the cokernel of $f$ is a $\tau$-rigid $A$-module. According to \cite{AIR14}, there exist $\mathbb{Z}$-linearly independent integer vectors $h_1,\dots, h_{10}\in \mathbb{Z}^{10}$  such that
\[
t\mathbf{g}=\sum_{i=1}^{10}a_ih_i,
\]
where $a_i\in \mathbb{N}$ for $1\leq i\leq 10$. Since $h_1,\dots, h_{10}$  is a $\mathbb{Z}$-basis of $\mathbb{Z}^{10}$, we also have integers $b_1,\dots, b_{10}$ such that
\[\mathbf{g}=\sum_{i=1}^{10}b_ih_i.\]
Hence, 
$t\mathbf{g}=\sum_{i=1}^{10}tb_ih_i=\sum_{i=1}^{10}a_ih_i$. It follows that $a_i=tb_i$ and $b_i\geq 0$ for each $i$. 
As a consequence, there exists $h\in \Hom(P^{\mathbf{g}_-},P^{\mathbf{g}_+})$ such that $\Hom_{K^b(\add A)}(T_h^\bullet,\Sigma T_h^\bullet)\neq 0$ and hence $\mathfrak{e}(\mathbf{g})=0$, a contradiction. In particular, $\mathfrak{e}(t\mathbf{g})>0$ for any positive integer $t$.
\end{proof}

\begin{remark}\label{r:e(g)=1}
Keeping the notation from Proposition~\ref{p:non-real-39}, we can even show that $\mathfrak{e}(\mathbf{g})=1$.
Consider the subspace $L_B$ of $\Hom(T^{\mathbf{g}_1},T^{\mathbf{g}_+})$  consisting of $BU+VB$ for any $U\in \End(T^{\mathbf{g}_-})$ and $V\in \End(T^{\mathbf{g}_+})$. We have shown that $\dim L_B<8$ for any $B$. On the other hand, one can check that there is a $B_0\in\Hom(T^{\mathbf{g}_1},T^{\mathbf{g}_+})$  such that $\dim L_{B_0}=7$. It follows that $\dim \Hom_{K^b(\add A)}(P^\bullet_{B_0},\Sigma P^\bullet_{B_0})=\dim \Hom(T^{\mathbf{g}_1},T^{\mathbf{g}_+})-\dim L_{B_0}=1$. Consequently, $\mathfrak{e}(\mathbf{g})=\dim \Hom_{K^b(\add A)}(P^\bullet_{{B_0}},\Sigma P^\bullet_{B_0})=1$.
\end{remark}

\begin{remark}\label{r:e(g)-neq-e(g,g)}
Keeping the notation from Proposition~\ref{p:non-real-39}, one can show that  $\mathfrak{e}(\mathbf{g},\mathbf{g})=0$ as follows. Let 
$B_1=
\begin{bmatrix}
0&1&0\\ 0&1&1\\ 1&0&0
\end{bmatrix}$ and $B_2=
\begin{bmatrix}
0&1&0\\ 1&1&0\\ 0&0&1
\end{bmatrix}$ 
be morphisms in 
$\Hom(T^{\mathbf{g}_-},T^{\mathbf{g}_+})$. 
%$\Hom(P{(125)}\oplus P{(126)}\oplus P{(134)}, P{(128)}\oplus P{(156)}\oplus P{(167)})$. 
It is straightforward to check that for any $X\in \Hom(T^{\mathbf{g}_-},T^{\mathbf{g}_+})$
%$\Hom(P{(125)}\oplus P{(126)}\oplus P{(134)}, P{(128)}\oplus P{(156)}\oplus P{(167)})$, 
there exist $U_1,U_2\in 
\End(T^{\mathbf{g}_-})$ 
%$\Hom(P{(125)}\oplus P{(126)}\oplus P{(134)},P{(125)}\oplus P{(126)}\oplus P{(134)})$ 
and $V_1,V_2\in \End(T^{\mathbf{g}_+})$ 
%$\Hom(P{(128)}\oplus P{(156)}\oplus P{(167)},P{(128)}\oplus P{(156)}\oplus P{(167)})$ 
such that \[X=B_2\circ U_1+V_1\circ B_1=B_1\circ U_2+V_2\circ B_2.\] Consequently, \[\Hom_{K^b(\add A)}(P^\bullet_{B_1},\Sigma P_{B_2}^\bullet)=0=\Hom_{K^b(\add A)}(P^\bullet_{B_2},\Sigma P_{B_1}^\bullet).\] We conclude that $\mathbb{E}(B_1,B_2)=0$ and hence $\mathfrak{e}(\mathbf{g},\mathbf{g})=0$ by Lemma \ref{l:E-invariant-ext} and Remark \ref{r:e-function-minimal}. On the other hand, let $M$ be the dominant monomial of $\mathcal{P}_\ell^+$ whose $\mathbf{g}$-vector is $\mathbf{g}$. One can show that $\chi_q(L(M))^2\neq \chi_q(L(M^2))$. Hence the requirement $L(M_1)$ and $L(M_2)$ are real in Conjecture~\ref{c:monoidal-additive-corr-compatible} is a necessary condition. 
\end{remark}

\subsection{An example of non-real $\mathbf{g}$-vector for $\Gr(4,8)$} \label{subsec:non-rigid-module-example-Gr48}
Now we provide an example for the tame case $\Gr(4,8)$. 
Consider 
\[
\widehat{T}= L_{1235}\oplus L_{1245}\oplus L_{1345}\oplus L_{1236}\oplus L_{1256}\oplus L_{1456}\oplus L_{1237}\oplus L_{1267}\oplus L_{1567}\oplus (\oplus_{i=1^8} L_{i,i+1,i+2,i+3}),
\]
reducing subscripts modulo $8$. The eight summands $L_{i,i+1,i+2,i+3}$ are projective-injective. We denote by $T$ the associated %basic 
cluster-tilting object in $\underline{\text{CM}}(B_{4,8})$. 

%Let $T=L_{1235}\oplus L_{1245}\oplus L_{1345}\oplus L_{1236}\oplus L_{1256}\oplus L_{1456}\oplus L_{1237}\oplus L_{1267}\oplus L_{1567}$ be a basic cluster-tilting object in $\underline{\text{CM}}(B_{4,8})$. 
The endomorphism algebra $A=\End_{\underline{\text{CM}}(B_{4,8})}(T)$ is isomorphic to the opposite of the Jacobian algebra of the following quiver 
\[
\adjustbox{scale=0.8}{\begin{tikzcd}
	1235 & 1245 & 1345 \\
	1236 & 1256 & 1456 \\
	1237 & 1267 & 1567
	\arrow["{\alpha_1}"', from=1-2, to=1-1]
	\arrow["{\alpha_2}"', from=1-3, to=1-2]
	\arrow["{\beta_1}"', from=2-2, to=2-1]
	\arrow["{\beta_2}"', from=2-3, to=2-2]
	\arrow["{\gamma_2}", from=3-3, to=3-2]
	\arrow["{\gamma_1}", from=3-2, to=3-1]
	\arrow["{\delta_2}", from=3-1, to=2-1]
	\arrow["{\delta_1}", from=2-1, to=1-1]
	\arrow["{u_1}", from=2-2, to=1-2]
	\arrow["{v_1}"', from=2-3, to=1-3]
	\arrow["{u_2}", from=3-2, to=2-2]
	\arrow["{v_2}"', from=3-3, to=2-3]
	\arrow["{x_1}", from=1-1, to=2-2]
	\arrow["{x_2}", from=1-2, to=2-3]
	\arrow["{y_1}", from=2-1, to=3-2]
	\arrow["{y_2}", from=2-2, to=3-3]
\end{tikzcd}}\]
with potential 
\[
W=\alpha_1 x_1u_1+\alpha_2x_2v_1+\beta_1y_1u_2+\beta_2y_2v_2-\delta_1x_1\beta_1-u_1x_2\beta_2-\delta_2y_1\gamma_1-u_2y_2\gamma_2.
\]

Let $\mathbf{g}= \scalemath{0.6}{ \begin{pmatrix}
    0&-1&0\\ -1&0&1\\ 0&1&0
\end{pmatrix} } \in \mathbb{Z}^9$. This is the $\mathbf{g}$-vector of the tableau $\ytableausetup{centertableaux}
\scalebox{0.6}{\begin{ytableau}
1 & 2 \\
3 & 4 \\
5 & 6 \\
7 & 8
\end{ytableau}}$ (see Remark~\ref{rem:g-vector_to_tableau}).
It is known that $\mathbf{g}$ is non-real, see \cite[Section 8]{CDFL20}. We also expect that no positive multiple of $\mathbf{g}$ corresponds to a real simple module. 

We proceed to show that $\mathfrak{e}(t\mathbf{g})>0$ for any $t\in \mathbb{Z}_{>0}$. 

We have $T^{\mathbf{g}_-}=L_{1236}\oplus L_{1245}$ and 
$T^{\mathbf{g}_+}=L_{1267}\oplus L_{1456}$
The indecomposable projective 
$A$-modules appearing in these summands are listed here: 
\[\adjustbox{scale=0.76}{
\begin{tikzcd}
{P{(1236)}:} & \mathbb{k} & 0 & 0 & {P{(1245)}:} & \mathbb{k} & \mathbb{k} & 0 \\
 & \mathbb{k}& \mathbb{k} & 0 && 0 & \mathbb{k} & \mathbb{k} \\
 & 0 & \mathbb{k}& 0 && 0 & 0 & 0 \\
{P{(1267)}:} & \mathbb{k} & \mathbb{k} & 0 & {P{(1456)}:} & \mathbb{k} & \mathbb{k} & \mathbb{k} \\
 & \mathbb{k} & \mathbb{k} & 0 && \mathbb{k} & \mathbb{k} &\mathbb{k} \\
 & \mathbb{k} & \mathbb{k} & 0 && 0 & 0 & 0
\arrow["1", from=2-2, to=1-2]
\arrow[from=3-2, to=2-2]
\arrow[from=2-3, to=1-3]
\arrow[from=2-4, to=1-4]
\arrow[from=3-4, to=2-4]
\arrow["1", from=3-3, to=2-3]
\arrow[from=1-3, to=1-2]
\arrow[from=1-4, to=1-3]
\arrow["0"{description}, from=2-3, to=2-2]
\arrow[from=2-4, to=2-3]
\arrow[from=3-3, to=3-2]
\arrow[from=3-4, to=3-3]
\arrow["1"{description}, from=1-2, to=2-3]
\arrow[from=1-3, to=2-4]
\arrow["1"{description}, from=2-2, to=3-3]
\arrow[from=2-3, to=3-4]
\arrow[from=3-7, to=3-6]
\arrow[from=2-7, to=2-6]
\arrow["1"', from=1-7, to=1-6]
\arrow[from=1-8, to=1-7]
\arrow["1"', from=2-8, to=2-7]
\arrow[from=3-8, to=3-7]
\arrow[from=3-6, to=2-6]
\arrow[from=2-6, to=1-6]
\arrow["0"{description}, from=2-7, to=1-7]
\arrow[from=2-8, to=1-8]
\arrow[from=3-7, to=2-7]
\arrow[from=3-8, to=2-8]
\arrow["1"{description}, from=1-6, to=2-7]
\arrow["1"{description}, from=1-7, to=2-8]
\arrow[from=2-6, to=3-7]
\arrow[from=2-7, to=3-8]
\arrow["1", from=5-2, to=4-2]
\arrow["1"', from=5-3, to=4-3] 
\arrow[from=5-4, to=4-4]
\arrow["1", from=6-2, to=5-2]
\arrow["1"', from=6-3, to=5-3]
\arrow[from=6-4, to=5-4]
\arrow["1"', from=6-3, to=6-2]
\arrow[from=6-4, to=6-3]
\arrow[from=5-4, to=5-3]
\arrow["1"', from=5-3, to=5-2]
\arrow["1"', from=4-3, to=4-2]
\arrow[from=4-4, to=4-3]
\arrow["0"{description}, from=4-2, to=5-3]
\arrow[from=4-3, to=5-4]
\arrow["0"{description}, from=5-2, to=6-3]
\arrow[from=5-3, to=6-4]
\arrow["1", from=5-6, to=4-6]
\arrow[from=6-6, to=5-6]
\arrow["1"', from=4-7, to=4-6]
\arrow["1"', from=4-8, to=4-7]
\arrow["1"', from=5-7, to=5-6]
\arrow[from=6-7, to=6-6]
\arrow["1", from=5-7, to=4-7]
\arrow[from=6-7, to=5-7]
\arrow["1", from=5-8, to=4-8]
\arrow[from=6-8, to=5-8]
\arrow["1"', from=5-8, to=5-7]
\arrow[from=6-8, to=6-7]
\arrow["0"{description}, from=4-6, to=5-7]
\arrow[from=5-6, to=6-7]
\arrow["0"{description}, from=4-7, to=5-8]
\arrow[from=5-7, to=6-8]
\end{tikzcd}}
\]

The $\Hom$-space between these projective indecomposables are listed in Table~\ref{ta:hom-space-48}. 

\begin{table}[ht]
\centering
\begin{tabular}{|c|c|c|c|c|}
\hline
 &$P{(1236)}$&$P{(1245)}$&$P{(1267)}$&$P{(1456)}$  \\ \hline
        $P{(1236)}$ &$\mathbb{k}$ &0&$\mathbb{k}$&$\mathbb{k}$\\ \hline
       $ P{(1245)}$&0&$\mathbb{k}$&$\mathbb{k}$&$\mathbb{k}$\\ \hline
$P{(1267)}$&$\mathbb{k}$&0&$\mathbb{k}$&0\\ \hline
$P{(1456)}$ &0&$\mathbb{k}$&0&$\mathbb{k}$\\ \hline   
\end{tabular}\caption{$\Hom$-spaces for $\text{CM}(B_{4,8})$}\label{ta:hom-space-48}
\end{table}

Similarly as in 
Proposition~\ref{p:non-real-39}, one can prove the following, providing further evidence for Conjecture \ref{c:monoidal-additive-correspondence}.
\begin{proposition}
We have $\mathfrak{e}(t\mathbf{g})>0$ for any positive integer $t$. 
\end{proposition}

\section{%Correspondence between 
Reachable real prime modules and reachable rigid indecomposable modules} \label{sec:correspondence between real modules and rigid modules}

{In this section, we give an explicit method to construct reachable rigid indecomposable modules from reachable prime real modules (equivalently, reachable prime real tableaux).}

%We thus establish a correspondence between the two.

%%
%
\subsection{Mutations of tableaux in $\SSYT(k,[n])$}\label{ss:mutation-tableaux}

Mutation in the Grassmannian cluster algebra $\CC[\Gr(k,n)]$ can be described using tableaux as has been explained in Section 4 of \cite{CDFL20}, using the fact that the cluster variables of 
$\CC[\Gr(k,n)]$ bijectively correspond to (reachable) prime real tableaux in $\SSYT(k,[n])$ (\cite[Theorem 3.25]{CDFL20}), see below for the definition of reachable tableaux. 

There is a partial order on the set $\SSYT(k,[n])$ of semistandard Young tableaux, induced from the dominance order on partitions. 

Let $\lambda = (\lambda_1,\dots,\lambda_r)$ and $\mu = (\mu_1,\dots,\mu_r)$ be two partitions, i.e. $\lambda_1 \ge \cdots \ge \lambda_r \ge 0$, $\mu_1 \ge \cdots \ge \mu_r \ge 0$. Then we say that 
{\em $\lambda \geq \mu$ in the dominance order} if $\sum_{j \leq i}\lambda_j \geq \sum_{j \leq i}\mu_j$ for $i=1,\dots,r$. 
{The {\em shape} of a tableau $\mathbf{T}$, denoted by ${\rm sh}(\mathbf{T})$, is the partition corresponding to the Young diagram obtained from $\mathbf{T}$ by forgetting its entries.}

For $i \in [n]$, let $\mathbf{T}[i]$ denote the restriction of $\mathbf{T} \in {\rm SSYT}(k,[n])$ to the entries in $[i]$. The {\em content} of a tableau $\mathbf{T}$ is the vector $(\nu_1,\dots,\nu_n) \in {\mathbb{N}^n}$, where $\nu_i$ is the number of $i$-filled boxes in $\mathbf{T}$. 

Let $\mathbf{T},\mathbf{T}' \in {\rm SSYT}(k, [n])$. We say that $\mathbf{T} \geq \mathbf{T}'$ if $\mathbf{T}$ and $\mathbf{T}'$ have the same content and if ${\rm sh}(\mathbf{T}[i]) \geq {\rm sh}(\mathbf{T}'[i])$ in the dominance order on partitions, for $i=1,\dots,n$.  

\begin{example}
Consider the tableaux $\mathbf{T}=\ytableausetup{centertableaux}
\scalebox{0.6}{\begin{ytableau}
1&3\\2&5\\4&6
\end{ytableau}}$ and $\mathbf{T'}=\ytableausetup{centertableaux}
\scalebox{0.6}{\begin{ytableau}
1&2\\3&4\\5&6
\end{ytableau}}$ in ${\rm SSYT}(3, [6])$. We have that 
\begin{align*}
& {\rm sh}(\mathbf{T}[1])= (1), \ {\rm sh}(\mathbf{T}[2])= (1,1), \ {\rm sh}(\mathbf{T}[3])= (2,1), \\
&{\rm sh}(\mathbf{T}[4])= (2,1,1), \  {\rm sh}(\mathbf{T}[5])= (2,2,1), \ {\rm sh}(\mathbf{T}[6])= (2,2,2), 
\end{align*}
and
\begin{align*}
& {\rm sh}(\mathbf{T'}[1])= (1), \ {\rm sh}(\mathbf{T'}[2])= (2), \ {\rm sh}(\mathbf{T'}[3])= (2,1), \\
& {\rm sh}(\mathbf{T'}[4])= (2,2), \ {\rm sh}(\mathbf{T'}[5])= (2,2,1), \ {\rm sh}(\mathbf{T'}[6])= (2,2,2).
\end{align*}
Since ${\rm sh}(\mathbf{T'}[i]) \ge {\rm sh}(\mathbf{T}[i] )$ for all $i$, we have that ${\bf T'} \ge {\bf T}$. 
\end{example}

We now recall the mutation sequence 
from~\cite[Section 4]{CDFL20}. 
%
%Every cluster variable $\mathbf{T}$ in $\SSYT(k, [n])$ corresponds to a cluster variable (denoted by $\ch(\mathbf{T})$) in $\CC[\Gr(k,n)]$. 
% 
To start with, take a seed consisting of Pl\"ucker coordinates (for example, the one from Section~\ref{subsec:Grassmannian cluster algebras} with triangular faces, cf. Figure~\ref{fig:5-9-quiver}) which correspond to one-column tableaux $\mathbf{S}_1,\dots, \mathbf{S}_m$, where  
$m=k(n-k)+1$. We perform the mutations using the following rule. Let $\mathbf{T}_1,\dots, \mathbf{T}_m$ be the tableaux in a seed (it can be the initial seed or a seed obtained from a sequence of mutations). Mutation at a (non-frozen) vertex $r$ 
changes the tableaux $\mathbf{T}_r$ as follows: 
%At each step, when we mutate at the vertex $r$, the cluster variable $\mathbf{T}_r$ is changed as follows:
\begin{align}
\mathbf{T}'_r = \mathbf{T}^{-1}_r \max\{\cup_{i \to r} \mathbf{T}_i, \cup_{r \to i} \mathbf{T}_i \}, \label{eq:mutation of tableaux}
\end{align}
where for two tableaux $\mathbf{T}$, $\mathbf{T'}$, 
$\max\{\mathbf{T}, \mathbf{T'}\}$ 
%$\max\{\cup_{i \to r} \mathbf{T}_i, \cup_{r \to i} \mathbf{T}_i \}$ 
takes the larger tableau with respect to 
$\ge$ and where the union is taken over arrows incident with $r$ of the quiver of the seed. 
Then we iterate this (doing quiver mutations along the way). 

{A tableau $\mathbf{T}$ in $\SSYT(k, [n])$ is called {\em reachable} if $\ch({\bf T})$ is a cluster monomial. If $\ch({\bf T})$ is a cluster monomial (resp. cluster variable), we also call ${\bf T}$ itself a cluster monomial (resp. cluster variable). Note that by \cite{Qin17} and \cite{KKOP20}, any cluster monomial is real. %a cluster monomial real. 
Therefore, any reachable tableau is also real (see Definition~\ref{def:real-tableau_vector} for the notion of real and prime tableaux).}

%%%%%%%%
%%
\subsection{Cluster variables in $\CC[\Gr(k,n)]$ and reachable rigid indecomposable objects of ${\rm CM}(B_{k,n})$}

To show that the category ${\rm CM}(B_{k,n})$ give an additive categorification of $\CC[\Gr(k,n)]$, Jensen, King and Su give a cluster character map in~\cite[\S 9]{JKS16}, 
giving the bijection between reachable rigid indecomposable modules in ${\rm CM}(B_{k,n})$ and cluster variables in $\CC[\Gr(k,n)]$.
We remark that the cluster character $\psi_{?}^{\hat{T}}$ in \cite{JKS16} is the composition of the Caldero Chapoton map ${\bf X}_{?}^{\hat{T}}$ in Section \ref{ss:additive-cat} with the isomorphism from the cluster algebra $\mathcal{A}$ to $\CC[\Gr(k,n)]$ in \cite{Sco} sending initial cluster variables to the corresponding Pl\"{u}cker coordinates.

Here we show how this bijection induces a bijection between (reachable) rigid indecomposable modules and (reachable) prime real semistandard Young tableaux.

% Recall the bijection between simple $U_q(\widehat{\mathfrak{sl}_k})$-module and tableaux in $\SSYT(k, [n], \sim)$ from Remark \ref{rem:corresp}. 
% A simple $U_q(\widehat{\mathfrak{sl}_k})$-module $L(M)$ is prime (resp. real) if and only if the corresponding tableau $T_M \in \SSYT(k, [n], \sim)$ is prime (resp. real). If a simple $U_q(\widehat{\mathfrak{sl}_k})$-module $L(M)$ is prime (resp. real), then the representative of $T_M$ in $\SSYT(k, [n])$ is prime. If a tableau $T$ in $\SSYT(k, [n])$ is prime (resp. real), then $L(M_T)$ is prime (resp. real).

We will need the following
\begin{theorem} [{\cite[Theorem 9.5]{JKS16}}] \label{thm:correspondence reachable rigid indecomposable modules and cluster variables}
{There is a one to one correspondence between the isomorphism classes of reachable rigid indecomposable modules in ${\rm CM}(B_{k,n})$ and cluster variables in $\CC[{\rm Gr}(k,n)]$.}
\end{theorem}

Under the notation of Definition~\ref{def:real-tableau_vector}, the following result relates cluster variables, reachable rigid indecomposable modules in ${\rm CM}(B_{k,n})$, reachable prime real modules, and reachable prime real tableaux.
\begin{theorem} \label{thm: cluster variables in Grkn and CMBkn} 
We have the following. 
(1) The map $N\mapsto \mathbf{T}_{\mathbf{g}(N)}$ gives a bijection 
between isomorphism classes of reachable rigid indecomposable modules in ${\rm CM}(B_{k,n})$ and reachable prime real semistandard Young tableaux in $\SSYT(k, [n])$. 

(2) The map $N \mapsto L(M_{ \mathbf{T}_{\mathbf{g}(N)} })$ gives a bijection between isomorphism classes of reachable rigid indecomposable modules in ${\rm CM}(B_{k,n})$ and isomorphism classes of reachable prime real modules in $\mathscr{C}_{\ell}$, $n=k+\ell+1$. 
\end{theorem} 

\begin{proof}
Part (1) follows from combining Theorem \ref{thm:correspondence reachable rigid indecomposable modules and cluster variables} with~\cite[Theorem 3.25]{CDFL20}. The former gives the bijection between isomorphism classes of reachable rigid indecomposables in ${\rm CM}(B_{k,n})$ and cluster variables of $\CC[\Gr(k,n)]$, the latter the bijection between cluster variables of $\CC[\Gr(k,n)]$ and reachable prime real semistandard Young tableaux in $\SSYT(k, [n])$. 

For part (2), one uses Theorem \ref{thm:correspondence reachable rigid indecomposable modules and cluster variables}, combined with the fact that  
~\cite[Theorem 1.2.1]{Qin17} and~\cite[Theorem 6.10]{KKOP22} give a one to one correspondence between isomorphism classes of reachable prime real modules in $\mathscr{C}_{\ell}$ ($n = k + \ell + 1$) and cluster variables in $\CC[\Gr(k,n)]$.
\end{proof}

{By Theorem \ref{thm: cluster variables in Grkn and CMBkn}, reachable rigid indecomposable module and reachable real prime tableau which correspond to each other have the same $\mathbf{g}$-vector. We have the following mutations in the cluster category. We start with the cluster-tilting object $S= S_1 \oplus \cdots \oplus S_m$ of ${\rm CM}(B_{k,n})$, where $S_1, \ldots, S_m$ are rank $1$ profiles corresponding to the Pl\"{u}cker coordinates in the initial cluster of $\CC[\Gr(k,n)]$ (with the initial cluster as in Section \ref{subsec:Grassmannian cluster algebras}). Let $M= M_1 \oplus \cdots \oplus M_m$ be the cluster-tilting object corresponding to a seed (it can be the initial seed or the seed obtained from a sequence of mutations). Each mutation in ${\rm CM}(B_{k,n})$ corresponds to mutation as described in (\ref{eq:mutation of tableaux}) in $\CC[\Gr(k,n)]$, where in ${\rm CM}(B_{k,n})$, the rigid indecomposable module corresponding to $\mathbf{T}'_r$ is the rigid indecomposable module which has the $\mathbf{g}$-vector $\mathbf{g}(\mathbf{T}'_r)$. 
Under the mutation, we replace $M_r$ by $M_r'$, where $M_r'$ has the same $\mathbf{g}$-vector as $\mathbf{T}_r'$.}

{Recall that the exchange graphs of ${\rm CM}(B_{3,9})$ and ${\rm CM}(B_{4,8})$ are connected (Theorem~\ref{thm:CMGr39-48-connected}). }
So here, we can remove the assumption of reachability from Theorem~\ref{thm: cluster variables in Grkn and CMBkn}: 

\begin{corollary}\label{cor:bijection-between-Gr3948-tableaux}
    Let $(k,n) \in \{(3,9), (4,8)\}$. Then there is a one to one correspondence between rigid indecomposable modules in ${\rm CM}(B_{k,n})$ and reachable real prime tableaux in $\SSYT(k,[n])$.
\end{corollary}

It is expected that assumption of reachability on $\SSYT(k,[n])$ is also not needed in Corollary~ \ref{cor:bijection-between-Gr3948-tableaux}, see Conjecture~\ref{c:monoidal-reach-conj}.

%%%%
%
\subsection{Rigid indecomposable objects and profiles}\label{sec:profiles}

In the rest of this section, we illustrate a method to obtain the indecomposable rigid (reachable) object of the Grassmannian cluster category for a given cluster variable and the corresponding profile (Section~\ref{subsec:grassmannian cluster categories}).

Given a cluster variable in $\CC[\Gr(k,n)]$, or equivalently, a tableau $\mathbf{T}$ in $\SSYT(k,[n])$ which corresponds to a cluster variable in $\CC[\Gr(k,n)]$, we can recover the reachable indecomposable object of ${\rm CM}(B_{k,n})$ it corresponds to as a mapping cone as we now explain. 
We first fix an initial cluster: It is convenient to work with the one from Section~\ref{subsec:Grassmannian cluster algebras}. It consists of Pl\"ucker coordinates. 
Let $T_1,\dots, T_m$ be the rank one modules corresponding to the chosen initial cluster. Then $\hat{T}=\bigoplus_{j\in [m]}T_j$ is the basic cluster-tilting object of ${\rm CM}(B_{k,n})$ corresponding to the initial cluster. We write $T$ for the associated basic cluster-tilting object of $\underline{\rm CM}(B_{k,n})$. 
 
We first compute the extended $\mathbf{g}$-vector $\hat{\mathbf{g}}:=(g_1, \ldots, g_m)$ of $\mathbf{T}$ using the method described in Remark~\ref{rem:g-vector_to_tableau}. Let $\mathbf{g}=\mathbf{g}(\bT)$, which is the truncation of $\hat{\mathbf{g}}$ obtained by deleting the components corresponding to the coefficients.
Let $J_+$ be the set of indices of the positive entries in $(g_1, \ldots, g_m)$ and let $J_-$ be the set of indices of the negative entries in $(g_1, \ldots, g_m)$. 
%Denote by $T$ the associated basic cluster-tilting object of $\underline{\rm CM}(B_{k,n})$ and $\mathbf{g}$ be the truncation of $\hat{\mathbf{g}}$ by deleting the components corresponding to the coefficient variables.
According to Section \ref{subsec:the function e minus}, there is a generic morphism from $T^{\mathbf{g}_-}$ to $T^{\mathbf{g}_+}$ whose mapping cone $M$ is rigid, and this is the rigid module in ${\rm CM}(B_{k,n})$ corresponding to $\mathbf{T}$. 
The associated short exact sequence in ${\rm CM}(B_{k,n})$ ending at $M$ is 
\[
0\to \bigoplus_{j\in J_-}T_j^{-g_j}\to \bigoplus_{j\in J_+}T_j^{g_j} \to M\to 0
\]

\begin{figure} 
\scalebox{0.7}{
\begin{tikzpicture}[scale=1] 
    % labels
    \foreach \i in {0,...,6}
      \path[black] (\i,-2) node{\i};
    % loop over the lattice points
    % \foreach \i in {0,...,8}
    %   \foreach \j in {0,...,8}{
    %     \draw (\i,\j) circle(3pt);    
  %    }; 

    \coordinate (a0m1u) at (0, -0.75);
    \coordinate (a6m1u) at (6, -0.75);

    \coordinate (a02) at (0, 2);
    \coordinate (a00) at (0, 0); 
    \coordinate (a11) at (1, 1);
    \coordinate (a1m1) at (1, -1);
    \coordinate (a22) at (2, 2);
    \coordinate (a20) at (2, 0); 
    \coordinate (a31) at (3, 1); 
    \coordinate (a42) at (4, 2);
    \coordinate (a40) at (4, 0);
    \coordinate (a51) at (5, 1);  
    \coordinate (a5m1) at (5, -1);   
    \coordinate (a62) at (6, 2); 
    \coordinate (a60) at (6, 0);
    \coordinate (a20) at (2, 0);  
    \coordinate (a40) at (4, 0);  

\node at (0.25, 0) {$1$}; 
\node at (1.25, -1) {$1$};
\node at (3.25, 1) {$1$};
\node at (5.25, -1) {$1$}; 
\node at (6.25, 0) {$1$}; 

    \draw (a00)--(a1m1);
    \draw (a1m1)--(a31);
    \draw (a31)--(a5m1);
    \draw (a5m1)--(a60);

    \draw[dashed] (a00)--(a0m1u);
    \draw[dashed] (a60)--(a6m1u);
    
\end{tikzpicture}   }
\scalebox{0.7}{
\begin{tikzpicture}[scale=1] 
    % labels
    \foreach \i in {0,...,6}
      \path[black] (\i,-2) node{\i};
    % loop over the lattice points
    % \foreach \i in {0,...,8}
    %   \foreach \j in {0,...,8}{
    %     \draw (\i,\j) circle(3pt);    
  %    }; 

    \coordinate (a00u) at (0, 0.25);
    \coordinate (a60u) at (6, 0.25);

    \coordinate (a02) at (0, 2);
    \coordinate (a00) at (0, 0);
    \coordinate (a0m1) at (0, -0.75);
     \coordinate (a00s) at (0, 0.1);
    \coordinate (a11) at (1, 1);
    \coordinate (a1m1) at (1, -1);
     \coordinate (a1m1s) at (0.9, -0.9);
    \coordinate (a22) at (2, 2);
    \coordinate (a22s) at (2.1, 2.1);
    \coordinate (a20) at (2, 0); 
    \coordinate (a31) at (3, 1); 
    \coordinate (a42) at (4, 2);
    \coordinate (a42s) at (3.9, 2.1);
    \coordinate (a40) at (4, 0);
    \coordinate (a51) at (5, 1);  
    \coordinate (a5m1) at (5, -1);
    \coordinate (a5m1s) at (5.1, -0.9);
    \coordinate (a62) at (6, 2); 
    \coordinate (a60) at (6, 0);
    \coordinate (a6m1) at (6, -0.75);
    \coordinate (a20) at (2, 0);  
    \coordinate (a40) at (4, 0);

   \node at (0.25, 2) {$1$};
   \node at (2.25, 2) {$1$};
   \node at (4.25, 2) {$1$};
   \node at (6.25, 2) {$1$};
   \node at (1.25, 1) {$2$};
   \node at (3.25, 1) {$3$};
   \node at (5.25, 1) {$2$};
   \node at (0.25, 0) {$3$};
   \node at (2.25, 0) {$3$};
   \node at (4.25, 0) {$3$};
   \node at (6.25, 0) {$3$};
   \node at (1.25, -1) {$3$};
   \node at (5.25, -1) {$3$};
    
    \draw (a02)--(a20);
    \draw (a20)--(a31);
    \draw (a31)--(a40);
    \draw (a40)--(a62);

    \draw[dashed] (a00)--(a22s);
    \draw[dashed] (a22s)--(a5m1s);
    \draw[dashed] (a5m1s)--(a60);

    \draw[thick, dotted] (a00s)--(a1m1s);
    \draw[thick, dotted] (a1m1s)--(a42s);
    \draw[thick, dotted] (a42s)--(a60); 
    
    \draw[dashed] (a02)--(a0m1);
    \draw[dashed] (a62)--(a6m1);
\end{tikzpicture}   }
\scalebox{0.7}{
\begin{tikzpicture}[scale=1] 
    % labels
    \foreach \i in {0,...,6}
      \path[black] (\i,-2) node{\i};
    % loop over the lattice points
    % \foreach \i in {0,...,8}
    %   \foreach \j in {0,...,8}{
    %     \draw (\i,\j) circle(3pt);    
  %    }; 

    \coordinate (a00) at (0, 0.25);
    \coordinate (a60) at (6, 0.25);

    \coordinate (a02) at (0, 2);
    \coordinate (a00) at (0, 0);
    \coordinate (a00s) at (0.1, 0.1);
    \coordinate (a11) at (1, 1);
    \coordinate (a1m1) at (1, -1);
    \coordinate (a1m1s) at (1.1, -0.9);
    \coordinate (a22) at (2, 2);
    \coordinate (a20) at (2, 0); 
    \coordinate (a31) at (3, 1); 
    \coordinate (a42) at (4, 2);
    \coordinate (a40) at (4, 0);
    \coordinate (a51) at (5, 1);  
    \coordinate (a5m1) at (5, -1);   
    \coordinate (a62) at (6, 2); 
    \coordinate (a60) at (6, 0);
    \coordinate (a20) at (2, 0);  
    \coordinate (a40) at (4, 0);

   \node at (0.25, 2) {$1$};
   \node at (2.25, 2) {$1$};
   \node at (4.25, 2) {$1$};
   \node at (6.25, 2) {$1$}; 
   \node at (1.25, 1) {$2$}; 
   \node at (3.25, 1) {$2$}; 
   \node at (5.25, 1) {$2$}; 
   \node at (0.25, 0) {$2$}; 
   \node at (2.25, 0) {$2$}; 
   \node at (4.25, 0) {$2$}; 
   \node at (6.25, 0) {$2$};

    \draw (a02)--(a11);
    \draw (a11)--(a22);
    \draw (a22)--(a31);
    \draw (a31)--(a42);
    \draw (a42)--(a51);
    \draw (a51)--(a62);

    \draw[dashed] (a00)--(a11);
    \draw[dashed] (a11)--(a20);
    \draw[dashed] (a20)--(a31);
    \draw[dashed] (a31)--(a40);
    \draw[dashed] (a40)--(a51);
    \draw[dashed] (a51)--(a60);

    \draw[dashed] (a0m1)--(a02);
    \draw[dashed] (a62)--(a6m1);
\end{tikzpicture}   }
\caption{A short exact sequence ending at the module $135|246$. 
%resolution of Cohen-Macaulay module $135|246$. The three pictures corresponds to the three terms in $P_{145} \to P_{124} \oplus P_{345} \oplus P_{156} \to 135|246$ respectively.
}
\label{fig:resolution-135-246}
\end{figure}
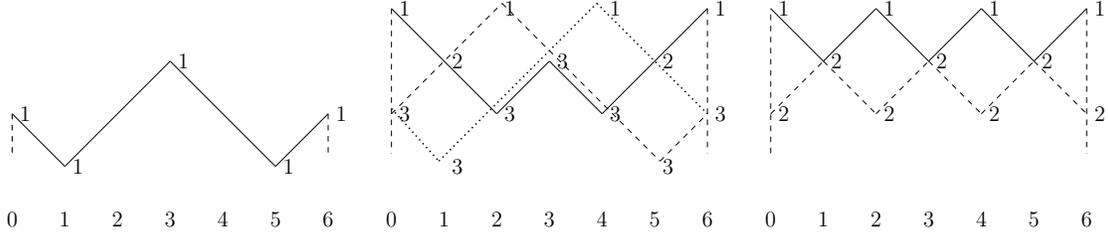

This short exact sequence can be used to in find the profile of the indecomposable $M$, if we know the profiles of the first and middle terms in the sequence. This works particularly well, if they have few summands (for example, if the first term is indecomposable). 
Figures~\ref{fig:resolution-246-135} and~\ref{fig:resolution-135-246} show two examples for this.

The profile of the module $M$ then arises from `subtracting' the profile of the first term from the one of the middle term: we subtract the dimension vectors of the lattice diagrams, see Sections 2.2 and 6.1 of~\cite{BBGL21}.

We remark that one has to draw the profiles of the summands in the appropriate places and that at present, there exists no general rule for doing this. In Section~\ref{sec:application to construction of indecomposable modules}, we provide more examples for this method.

\begin{example}\label{ex:rank2-indec-profile}
In this example, we illustrate the above and explain how to get to the (profiles for the) indecomposable rank two modules $135\mid 246$ and $246\mid 135$ of ${\rm CM}(B_{3,6})$.

(1) The module $246|135$ corresponds to the tableau $\scalemath{0.6}{\begin{ytableau} 
1 & 2 \\ 3 & 4 \\ 5 & 6 \end{ytableau}}$. The extended $\mathbf{g}$-vector of $\scalemath{0.6}{\begin{ytableau} 
1 & 2 \\ 3 & 4 \\ 5 & 6 \end{ytableau}}$ can be read from 
\begin{align} \label{eq:decomposition of tableau in terms of initial tableaux}
\scalemath{0.7}{
\begin{ytableau}
1 & 2  \\
3 & 4 \\
5 & 6 
\end{ytableau} = \left( \begin{ytableau}
1   \\
2  \\
6  
\end{ytableau} \cup \begin{ytableau}
1   \\
4  \\
5  
\end{ytableau} \cup \begin{ytableau}
2   \\
3  \\
4  
\end{ytableau} \right)
\begin{ytableau}
1   \\
2  \\
4  
\end{ytableau}^{-1}, }
\end{align}
where ``$\cup$'' and quotients of tableaux are defined in Section \ref{subsec:correspondence between modules and tableaux}. The expression (\ref{eq:decomposition of tableau in terms of initial tableaux}) corresponds to the following short exact sequence ending at the module $246|135$:
\begin{align*}
P_{124} \to P_{126} \oplus P_{145} \oplus P_{234} \to 246|135,
\end{align*}
see Figure~\ref{fig:resolution-246-135}.
%\ref{fig: resolution 246 135}. 

(2) 
The module $135|246$ corresponds to the tableau $\scalemath{0.6}{\begin{ytableau}
1 & 3 \\ 2 & 5 \\ 4 & 6 \end{ytableau}}$. The extended $\mathbf{g}$-vector of $\scalemath{0.6}{\begin{ytableau} 
1 & 3 \\ 2 & 5 \\ 4 & 6 \end{ytableau}}$ can be read from 
\begin{align} \label{eq:decomposition of tableau in terms of initial tableaux-124356}
\scalemath{0.7}{
\begin{ytableau}
1 & 3  \\
2 & 5 \\
4 & 6 
\end{ytableau} = \left( \begin{ytableau}
1   \\
2  \\
4  
\end{ytableau} \cup \begin{ytableau}
3   \\
4  \\
5  
\end{ytableau} \cup \begin{ytableau}
1   \\
5  \\
6  
\end{ytableau}
\right)
\begin{ytableau}
1   \\
4  \\
5  
\end{ytableau}^{-1}. }
\end{align}
Note that the tableaux appearing on the right hand side in (\ref{eq:decomposition of tableau in terms of initial tableaux}) are obtained from the ones on the right hand side of  (\ref{eq:decomposition of tableau in terms of initial tableaux-124356}) by adding $3 \pmod 6$ to each entry and sort each one-column tableau such that the result is standard (recall that $\cup$ is commutative). Accordingly, the short exact sequence ending at $135\mid 246$ can be obtained from the one in part (1) by adding $3$ to every $3$-subset. The lattice pictures for the short exact sequence can be obtained from Figure~\ref{fig:resolution-246-135} by rotating each lattice picture three steps around. See Figure~\ref{fig:resolution-135-246}.
\end{example}

\begin{figure} 
\scalebox{0.7}{
\begin{tikzpicture}[scale=1] 
    % labels
    \foreach \i in {0,...,6}
      \path[black] (\i,-1) node{\i};
    % loop over the lattice points
    % \foreach \i in {0,...,8}
    %   \foreach \j in {0,...,8}{
    %     \draw (\i,\j) circle(3pt);    
  %    }; 

    \coordinate (a00) at (0, 0.25);
    \coordinate (a60) at (6, 0.25);

    \coordinate (a02) at (0, 2);
    \coordinate (a13) at (1, 3);
    \coordinate (a11) at (1, 1);
    \coordinate (a22) at (2, 2);   
    \coordinate (a33) at (3, 3);   
    \coordinate (a31) at (3, 1); 
    \coordinate (a42) at (4, 2);  
    \coordinate (a51) at (5, 1);  
    \coordinate (a53) at (5, 3);   
    \coordinate (a62) at (6, 2);  
    \coordinate (a20) at (2, 0);  
    \coordinate (a40) at (4, 0);  

\node at (0.25, 2) {$1$}; 
\node at (2.25, 0) {$1$}; 
\node at (3.25, 1) {$1$}; 
\node at (4.25, 0) {$1$}; 
\node at (6.25, 2) {$1$}; 

    \draw (a02)--(a20);
    \draw (a20)--(a31);
    \draw (a31)--(a40);
    \draw (a40)--(a62);

    \draw[dashed] (a00)--(a02);
    \draw[dashed] (a62)--(a60);
    
\end{tikzpicture}   }
\scalebox{0.7}{
\begin{tikzpicture}[scale=1] 
    % labels
    \foreach \i in {0,...,6}
      \path[black] (\i,-1) node{\i};
    % loop over the lattice points
    % \foreach \i in {0,...,8}
    %   \foreach \j in {0,...,8}{
    %     \draw (\i,\j) circle(3pt);    
  %    }; 

    \coordinate (a00) at (0, 0.25);
    \coordinate (a60) at (6, 0.25);

    \coordinate (a02) at (0, 2);
    \coordinate (a02s) at (0.1, 2.1);
    \coordinate (a13) at (1, 3);
    \coordinate (a11) at (1, 1);
    \coordinate (a11s) at (1.1, 1.1);
    \coordinate (a22) at (2, 2);   
    \coordinate (a33) at (3, 3);   
    \coordinate (a31) at (3, 1); 
    \coordinate (a42) at (4, 2);  
    \coordinate (a51) at (5, 1);  
    \coordinate (a51s) at (4.9, 1.1);  
    \coordinate (a53) at (5, 3);   
    \coordinate (a62) at (6, 2); 
    \coordinate (a62s) at (5.9, 2.1); 
    \coordinate (a20) at (2, 0);  
    \coordinate (a40) at (4, 0);  

   \node at (1, 3.25) {$1$};
   \node at (3, 3.25) {$1$};
   \node at (5, 3.25) {$1$};

   \node at (0.25, 2) {$3$};
   \node at (2.25, 2) {$2$};
   \node at (4.25, 2) {$2$};
   \node at (6.25, 2) {$3$};
   \node at (1.25, 1) {$3$};
   \node at (3.25, 1) {$3$};
   \node at (5.25, 1) {$3$};
   \node at (2.25, 0) {$3$};
   \node at (4.25, 0) {$3$};
    
    \draw (a02)--(a13);
    \draw (a13)--(a40);
    \draw (a40)--(a62);

    \draw[dashed] (a02)--(a20);
    \draw[dashed] (a20)--(a53);
    \draw[dashed] (a53)--(a62);

    \draw[thick, dotted] (a02s)--(a11s);
    \draw[thick, dotted] (a11s)--(a33);
    \draw[thick, dotted] (a33)--(a51s);
    \draw[thick, dotted] (a51s)--(a62s);

    \draw[dashed] (a00)--(a02);
    \draw[dashed] (a62)--(a60);
    
\end{tikzpicture}   }
\scalebox{0.7}{
\begin{tikzpicture}[scale=1] 
    % labels
    \foreach \i in {0,...,6}
      \path[black] (\i,-1) node{\i};
    % loop over the lattice points
    % \foreach \i in {0,...,8}
    %   \foreach \j in {0,...,8}{
    %     \draw (\i,\j) circle(3pt);    
  %    }; 

    \coordinate (a0m1) at (0, -0.75);
    \coordinate (a6m1) at (6, -0.75);

    \coordinate (a02) at (0, 2);
    \coordinate (a13) at (1, 3);
    \coordinate (a11) at (1, 1);
    \coordinate (a22) at (2, 2);   
    \coordinate (a33) at (3, 3);   
    \coordinate (a31) at (3, 1); 
    \coordinate (a42) at (4, 2);  
    \coordinate (a51) at (5, 1);  
    \coordinate (a53) at (5, 3);   
    \coordinate (a62) at (6, 2);  

   \node at (1, 3.25) {$1$};
   \node at (3, 3.25) {$1$};
   \node at (5, 3.25) {$1$};

   \node at (0.25, 2) {$2$};
   \node at (2.25, 2) {$2$};
   \node at (4.25, 2) {$2$};
   \node at (6.25, 2) {$2$};

    \draw (a02)--(a13);
    \draw (a13)--(a22);
    \draw (a22)--(a33);
    \draw (a33)--(a42);
    \draw (a42)--(a53);
    \draw (a53)--(a62);

    \draw[dashed] (a02)--(a11);
    \draw[dashed] (a11)--(a22);
    \draw[dashed] (a22)--(a31);
    \draw[dashed] (a31)--(a42);
    \draw[dashed] (a42)--(a51);
    \draw[dashed] (a51)--(a62);

    \draw[dashed] (a00)--(a02);
    \draw[dashed] (a62)--(a60);
\end{tikzpicture}   }
\caption{A short exact sequence ending at the 
%A presentation of the %Cohen-Macaulay 
module $246|135$. 
%The three pictures corresponds to the three terms in $P_{124} \to P_{126} \oplus P_{145} \oplus P_{234} \to 246|135$ respectively.
}
\label{fig:resolution-246-135}%{fig: resolution 246 135}
\end{figure}
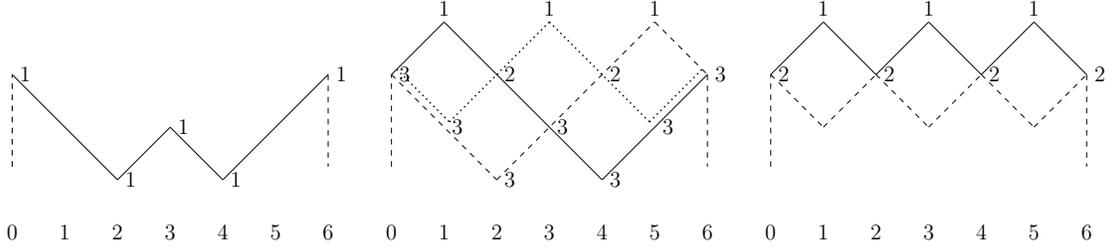

One can produce many new cluster variables from known cluster variables through mutation or braid group actions as we will show later, see Section~\ref{subsec:Quasi-homomorphisms of cluster algebras and braid group actions}.
Using Theorem \ref{thm: cluster variables in Grkn and CMBkn}, one can produce many rigid indecomposable modules from known rigid indecomposable modules. We will apply this result to produce rigid indecomposable modules of the categories ${\rm CM}(B_{3,9})$ and ${\rm CM}(B_{4,8})$ in Sections \ref{subsec:rigid indecomposable modules Gr39} and \ref{subsec:rigid indecomposable modules Gr48} respectively. 

% \begin{example}
% $\Gr(3,6)$
% \end{example}

\section{Realization of Hernandez-Leclerc's generic kernels} \label{sec:realization of generic kernels of Hernandez-Leclerc}
In this section, we apply results of Section \ref{sec:correspondence between additive and monoid categorifications} to give an alternative construction of the generic kernels introduced by Hernandez and Leclerc in \cite{HL16} for type $\mathbb{A}$ via Grassmannian cluster categories. 
As an application, this will give a criterion for checking whether two cluster variables corresponding to Kirillov-Reshetikhin modules belong to a common cluster, see Remark~\ref{rem:criterion-compatibility}.

\subsection{The generic kernels of Hernandez-Leclerc for type $\mathbb{A}$} \label{subsec:Application to generic kernels of Hernandez-Leclerc type A}
Let $C=(c_{ij})$ be an indecomposable Cartan matrix of type $\mathbb{A}_{k-1}$. Denote by $I$ {the set} $\{1,\dots, k-1\}$. So $c_{ii}=2$ for all $i$, $c_{ij}=-1$ whenever $|j-i|=1$ and all other entries are $0$. We first recall the semi-infinite quiver $\Gamma^-$ Hernandez and Leclerc associated with $C$: 
\begin{itemize}
    \item 
    The vertices of $\Gamma^-$ are labeled by $(i,m)$, where $i\in I$ and $m\in \mathbb{Z}_{<0}$ satisfy that $i\equiv 1 \pmod 2, 2\mid m$ or $i\equiv 0\pmod 2, 2\mid m-1$. We write $V^-$ to denote the set of vertices of $\Gamma^-$.
    \item 
    The arrows of $\Gamma^-$ consist of $(i,m+2)\to (i,m)$ for any $(i,m)$ and $(i,m-1)\to (j,m)$ whenever $c_{ij}=-1$.
\end{itemize} 
We draw $\Gamma^-$ in the plane using the vertices as coordinate positions, see Example~\ref{ex:A3-HL-quiver} or the left hand of Figure~\ref{fig:initial HL-quiver and initial quiver of Grassmannian cluster algebra}.
We remark that the quiver $\Gamma^-$ is the opposite of the one defined by Hernandez and Leclerc in \cite{HL16}, since we are working with right modules.

\begin{example}\label{ex:A3-HL-quiver}
The Hernandez-Leclerc quiver $\Gamma^-$ for $C$ of type $\mathbb{A}_3$ is as follows
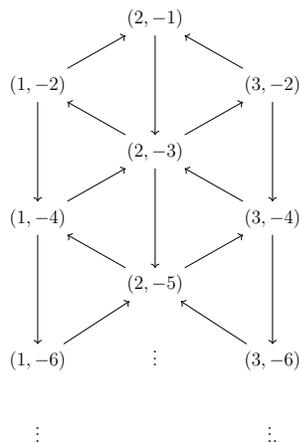
\[\adjustbox{scale=0.6}{\begin{tikzcd}
 & {(2,-1)} \\
 {(1,-2)} && {(3,-2)} \\
 & {(2,-3)} \\
 {(1,-4)} && {(3,-4)} \\
 & {(2,-5)} \\
 {(1,-6)} & \vdots & {(3,-6)} && {} \\
 \vdots && \vdots.
 \arrow[from=2-1, to=4-1]
 \arrow[from=1-2, to=3-2]
 \arrow[from=2-3, to=4-3]
 \arrow[from=4-3, to=6-3]
 \arrow[from=3-2, to=5-2]
 \arrow[from=4-1, to=6-1]
 \arrow[from=2-1, to=1-2]
 \arrow[from=2-3, to=1-2]
 \arrow[from=3-2, to=2-1]
 \arrow[from=4-3, to=3-2]
 \arrow[from=5-2, to=4-1]
 \arrow[from=6-3, to=5-2]
 \arrow[from=3-2, to=2-3]
 \arrow[from=4-1, to=3-2]
 \arrow[from=5-2, to=4-3]
 \arrow[from=6-1, to=5-2]
\end{tikzcd}}
\]
\end{example}
We denote by $W$ the formal sum of all  oriented $3$-cycles of $\Gamma^-$ up to cyclic permutations and by $A=J(\Gamma^-,W)$ the Jacobian algebra associated with $(\Gamma^-,W)$, \cite{DWZ08}. It is defined to be the path algebra of $\Gamma^-$ up to the relations from the potential $W$: whenever $\alpha$ is an internal arrow of the quiver, it is incident with exactly two triangles, $p\alpha$ and $q\alpha$ (where $p,q$ are paths of length $2$). In that case, we have $p=q$ in $A$. 
If $\alpha$ is only incident with one triangle, say $p\alpha$, then the relations of $W$ enforce $p=0$.

Denote by $I(i,m)$ be the indecomposable injective right $A$-module associated with vertex $(i,m)$. One can show that $\Hom_{A}(I(i,m),I(i,m-2v))$ is a finite  dimensional vector space over $\mathbb{k}$ for every $v\ge 1$. Hernandez and Leclerc proved in~\cite{HL16} that there exists a morphism $f\in \Hom_{A}(I(i,m),I(i,m-2v))$ such that $\ker f$ is finite dimensional. Similarly as in 
Section~\ref{sec:correspondence between additive and monoid categorifications} (cf. also \cite{P13}), there is an open dense subset $\mathcal{O}$ of $\Hom_{A}(I(i,m),I(i,m-2v))$ such that the kernels of all elements of $\mathcal{O}$ are finite dimensional and have the same $F$-polynomial. The elements of $\mathcal{O}$ are called the {\it generic homomorphisms} from $I(i,m)$ to $I(i,m-2v)$.
\begin{definition}\cite[Definition 4.5]{HL16}\label{def:generic-kernel}
Let $(i,m)$ be a vertex of $\Gamma^-$ and $v\ge 1$. Then we denote the kernel of a generic $A$-module homomorphism from $I(i,m)$ to $I(i,m-2v)$ by $K_{v,m}^{(i)}$ and we call it a generic kernel. 
%For any $(i,m)$ and $v\geq 1$, denote by $K_{v,m}^{(i)}$  the kernel of a generic $A$-module homomorphism from $I(i,m)$ to $I(i,m-2v)$.
\end{definition}

We want to determine the $K_{v,m}^{(i)}$. To do so, it is enough to consider modules over a truncated finite-dimensional version of $A$. We recall this now. 

% \begin{lemma}
%    The $A$-module $K_{v,m}^{(i)}$ is independent of the choice of the generic homomorphism and it is $\tau^{-1}$-rigid, {\it i.e.}, let $f\in \Hom_A(I(i,m),I(i,m-2v))$  be a generic homomorphism, then $\Hom_A(K_{v,m}^{(i)}, I(i,m))\xrightarrow{\Hom_A(K_{v,m}^{(i)},f)}\Hom_A(K_{v,m}^{(i)}, I(i,m-2v))$ is surjective.
% The following is a consequence of \cite{HL16}.
% \end{lemma}

For $s\in \mathbb{Z}_{<0}$, let $\Gamma_s^-$ be the full subquiver of $\Gamma^-$ with  set of vertices $V_s^-:=\{(i,m)\in V^-|m\geq s\}$, i.e. the finite quiver obtained at truncating $\Gamma^-$ below level $s$. Denote by $W_s$ the sum of the oriented $3$-cycles of $\Gamma_s^-$ and by $A_s=J(\Gamma_s^-,W_s)$ the (finite dimensional) associated Jacobian algebra. 
For $(i,m)\in V_s^-$, we write 
$I_s(i,m)$ to denote the indecomposable injective right $A_s$-module associated with $(i,m)$. 

To determine $K_{v,m}^{(i)}$, it suffices to work with $A_s$ and its injective indecomposable modules for some $s \ll 0$: by~\cite[Section 4.5.3]{HL16}, $K_{v,m}^{(i)}$ is the kernel of a generic right $A_s$-module homomorphism from $I_{s}(i,m)$ to $I_s(i,m-2v)$. Furthermore, it suffices to take $s\le m-2v-2$.

%Note that $A_s$ is finite dimensional. 
Let $\mathcal{C}(s):=\mathcal{C}(\Gamma_s^-,W_s)$ be the generalized cluster category associated with $(\Gamma_s^-,W_s)$ in the sense of \cite{Amiot09}. Let $T=\bigoplus\limits_{(i,m)\in V_s^-}T_{i,m}$ be a basic cluster-tilting object with indecomposable direct summands $T_{i,m}$ of $\mathcal{C}(s)$ such that $A_s\cong \End(T)$ (such $T_{i,m}$ always exist). 
It is known that  $\Hom_{\mathcal{C}(s)}(T,-):\mathcal{C}(s)\to \operatorname{mod} A_s$ induces an equivalence $\mathcal{C}(s)/\add \Sigma T\cong \operatorname{mod} A_s$, where $\mathcal{C}(s)/\add \Sigma T$ is the additive quotient of $\mathcal{C}(s)$ by morphisms which factor through objects in $\add\Sigma T$.

Denote by $\tilde{K}_{v,m}^{(i)}$ a preimage in $\mathcal{C}(s)$ of $K_{v,m}^{(i)}$.
According to \cite[Theorem 4.8 and Proposition 4.16]{HL16},
the image $\mathbf{X}_{\tilde{K}_{v,m}^{(i)}}^{\Sigma T}$ of $\tilde{K}_{v,m}^{(i)}$ under the Caldero-Chapoton map with respect to $\Sigma T$ (Section~\ref{ss:additive-cat}) 
is a cluster variable of the cluster algebra defined by $\Gamma_s^-$. It follows from Theorem~\ref{t:additve-cat} that 
$\tilde{K}_{v,m}^{(i)}$  is an in\-decomposable rigid object. Moreover, 
there is a triangle
\[
\Sigma^{-1}\tilde{K}_{v,m}^{(i)}\to \Sigma T_{i,m}\xrightarrow{g}\Sigma T_{i,m-2v}\to \tilde{K}_{v,m}^{(i)}.
\]

As $\tilde{K}_{v,m}^{(i)}$ is rigid, 
the orbit of $\Sigma g$ is an open dense subset of $\Hom_{\mathcal{C}(s)}(\Sigma^2 T_{i,m}, \Sigma^2 T_{i,m-2v})$ (cf. the proof of Lemma \ref{l:ext-dim}). We conclude that
% According to the proof of Lemma \ref{l:ext-dim}, we deduce that
$\Hom_{\mathcal{C}(s)}(T,\Sigma g)$ is a generic $A_s$-module homomorphism from $I_s(i,m)$ to $I_s(i,m-2v)$. 

\subsection{Generic kernels via Grassmannian cluster categories} \label{subsec:Generic kernels via Grassmannian cluster categories}

In their paper, Hernandez and Leclerc categorified the cluster algebra associated to $\Gamma^-_{-2\ell-2}$ using $U_q(\widehat{\mathfrak{sl}_k})$-modules. 
The modules corresponding to the vertices of the initial quiver $\Gamma^-_{-2\ell-2}$ are certain Kirillov-Reshetikhin modules.  
% Let $\xi:I\to \ZZ$ be the
% height function %$\xi:I \to \ZZ$, 
% \[
% \xi(i) = 
% \left\{ 
% \begin{array}{ll}
%     -1 & \mbox{if $i$ is odd,} \\
%     0 & \mbox{if $i$ is even.} 
% \end{array}
% \right.
% \]

\begin{figure}
    \centering
    \hspace{-3.5cm}
    \begin{minipage}{.35\textwidth}
    \centering
     \adjustbox{scale=0.37}{
\begin{tikzcd}[arrows=<-]
	& {L(2_0)} (2,-1) && {L(4_0)} (4,-1) \\
	{L(1_{-1})} (1,-2) && {L(3_{-1})} (3,-2) \\
	& {L(2_{-2}2_0)} (2,-3) && {L(4_{-2}4_0)} (4,-3) \\
	{L(1_{-3}1_{-1})} (1,-4) && {L(3_{-3}3_{-1})} (3,-4) \\
	& {L(2_{-4}2_{-2}2_0)} (2,-5) && {L(4_{-4}4_{-2}4_0)} (4,-5) \\
	{L(1_{-5}1_{-3}1_{-1})} (1,-6) && {L(3_{-5}3_{-3}3_{-1})} (3,-6) \\
	& {\fbox{$L(2_{-6}\cdots 2_0)$}} (2,-7) && {\fbox{$L(4_{-6}\cdots 4_0)$}} (4,-7) \\
	{\fbox{$L(1_{-7}\cdots 1_{-1})$}} (1,-8) && {\fbox{$L(3_{-7}\cdots 3_{-1})$}} (3,-8)
	\arrow[from=1-2, to=2-1]
	\arrow[from=1-2, to=2-3]
	\arrow[from=1-4, to=2-3]
	\arrow[from=3-4, to=1-4]
	\arrow[from=4-3, to=2-3]
	\arrow[from=2-3, to=3-4]
	\arrow[from=3-4, to=4-3]
	\arrow[from=5-4, to=3-4]
	\arrow[from=2-3, to=3-2]
	\arrow[from=2-1, to=3-2]
	\arrow[from=4-1, to=2-1]
	\arrow[from=3-2, to=4-1]
	\arrow[from=4-1, to=5-2]
	\arrow[from=5-2, to=6-1]
	\arrow[from=6-1, to=7-2]
	\arrow[from=7-2, to=8-1]
	\arrow[from=8-1, to=6-1]
	\arrow[from=6-1, to=4-1]
	\arrow[from=3-2, to=4-3]
	\arrow[from=4-3, to=5-2]
	\arrow[from=5-2, to=6-3]
	\arrow[from=6-3, to=7-2]
	\arrow[from=7-2, to=8-3]
	\arrow[from=8-3, to=6-3]
	\arrow[from=6-3, to=4-3]
	\arrow[from=4-3, to=5-4]
	\arrow[from=5-4, to=6-3]
	\arrow[from=6-3, to=7-4]
	\arrow[from=7-4, to=5-4]
	\arrow[from=7-4, to=8-3]
	\arrow[from=7-2, to=5-2]
	\arrow[from=5-2, to=3-2]
	\arrow[from=3-2, to=1-2]
\end{tikzcd}
}
\end{minipage} 
    \hspace{1.5cm}
 \begin{minipage}{.35\textwidth} 
 \centering
 \adjustbox{scale=0.5}{ 
\begin{tikzcd}[arrows=<-]
	{L(1_{-1})} (1,-2) & {L(2_0)} (2,-1) & {L(3_1)} (3,-2) & {L(4_2)} (4,-1) \\
{L(1_{-3}1_{-1})} (1,-4) & {L(2_{-2}2_0)} (2,-3) & {L(3_{-1}3_1)} (3,-4) & {L(4_04_2)} (4,-3) \\
{L(1_{-5}1_{-3}1_{-1})} (1,-6) & {L(2_{-4}2_{-2}2_0)} (2,-5) & {L(3_{-3}3_{-1}3_1)} (3,-6) & {L(4_{-2}4_04_2)} (4,-5) \\
{\fbox{$L(1_{-7}\cdots 1_{-1})$}} (1,-8) & {\fbox{$L(2_{-6}\cdots 2_0)$}} (2,-7) & {\fbox{$L(3_{-5}\cdots 3_1)$}} (3,-8) & {\fbox{$L(4_{-4}\cdots 4_2)$}} (4,-7) 
	\arrow[from=1-2, to=1-1]
	\arrow[from=1-3, to=1-2]
	\arrow[from=1-4, to=1-3]
	\arrow[from=2-1, to=1-1]
	\arrow[from=3-1, to=2-1]
	\arrow[from=4-1, to=3-1]
	\arrow[from=2-2, to=1-2]
	\arrow[from=3-2, to=2-2]
	\arrow[from=4-2, to=3-2]
	\arrow[from=4-3, to=3-3]
	\arrow[from=3-3, to=2-3]
	\arrow[from=2-3, to=1-3]
	\arrow[from=4-4, to=3-4]
	\arrow[from=3-4, to=2-4]
	\arrow[from=2-4, to=1-4]
	\arrow[from=2-4, to=2-3]
	\arrow[from=2-3, to=2-2]
	\arrow[from=2-2, to=2-1]
	\arrow[from=3-4, to=3-3]
	\arrow[from=3-3, to=3-2]
	\arrow[from=3-2, to=3-1]
	\arrow[from=1-1, to=2-2]
	\arrow[from=2-1, to=3-2]
	\arrow[from=3-1, to=4-2]
	\arrow[from=1-2, to=2-3]
	\arrow[from=2-2, to=3-3]
	\arrow[from=3-2, to=4-3]
	\arrow[from=1-3, to=2-4]
	\arrow[from=2-3, to=3-4]
	\arrow[from=3-3, to=4-4]
\end{tikzcd}
}
\end{minipage} 
\caption{The quiver on the right is the initial quiver of the cluster algebra $K_0(\mathscr{C}_{3}^{\mathfrak{sl}_5})$ (Section \ref{ss:monoidal-cat}). The (mutable part of the) quiver on the left can be obtained from the second quiver by mutating column 4, then column 3 (Remark~\ref{rem:mut-sequence}). {We write $Y_{i,s}$ as $i_s$.}}
\label{fig:initial HL-quiver and initial quiver of Grassmannian cluster algebra}
\end{figure}
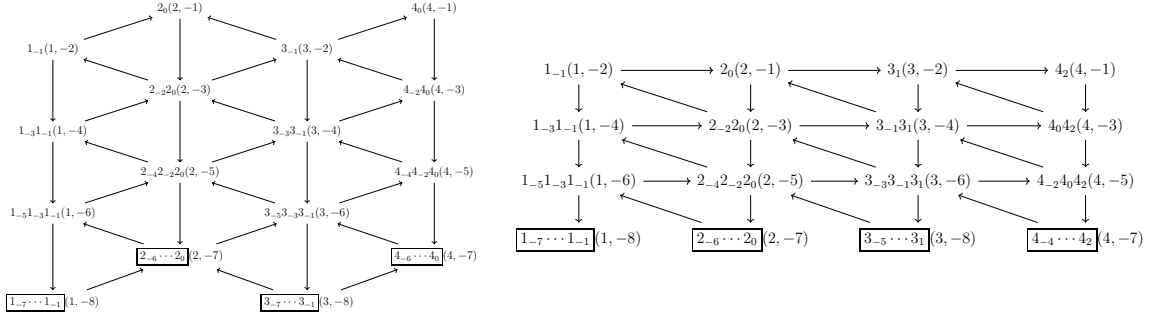

%$\xi(i)=-1$ if $i$ is odd and $\xi(i)=0$ if $i$ is even. 
Let $\xi':I\to \ZZ$ be the height function 
$\xi'(i)=0$ if $i$ is even, $\xi'(i)=-1$ if $i$ is odd. 
Then the $U_q(\widehat{\mathfrak{sl}_k})$-module at position $(i,m)$ is $L(Y_{i,m+1}Y_{i,m+3} \cdots Y_{i,\xi'(i)})$, see the first picture in Figure~\ref{fig:initial HL-quiver and initial quiver of Grassmannian cluster algebra}, for $\Gamma_{-8}^-$, with $U_q(\widehat{\mathfrak{sl}_5})$-modules. 

The quiver $\Gamma_{-2\ell-2}^-$ is mutation equivalent to the quiver of an initial seed for $K_0(\mathscr{C}_{\ell}^{\mathfrak{sl}_k})$ (up to frozen vertices) as we explain now. 

Fix $k,\ell$ as above. Let $Q_{\ell}$ be the quiver with $(k-1)(\ell+1)$ vertices from Example~\ref{ex:quiver-Q-l}. 
The second quiver in Figure~\ref{fig:initial HL-quiver and initial quiver of Grassmannian cluster algebra} shows the quiver $Q_{3}$.

\begin{remark}\label{rem:mut-sequence} 
We compare $\Gamma_{-2\ell-2}^-$ and $Q_{\ell}$, drawing the latter in the plane, {and recording the  coordinates of its vertices} (the columns of $Q_{\ell}$ are indexed by the first entry $i\in I$). 

Both $\Gamma_{-2\ell-2}^-$ and $Q_{\ell}$ are formed by oriented triangles but they are not isomorphic (for example, the arrows on the top of the quiver $Q_{\ell}$ all point to the right). 
One can give an explicit mutation sequence to get the non-frozen part of  $\Gamma_{-2\ell-2}^-$ from non-frozen part of $Q_{\ell}$. 
%The quiver $\Gamma^-_{-2\ell-2}$ can be obtained from the quiver $Q_{\ell}$ by performing the following mutations sequence. 
%
First, one mutates at the vertices of the last column (column $k-1$) of $Q_{\ell}$ (in the order from top to bottom), then the second last column, etc., until the third column of $Q_{\ell}$. Then one starts again with the last column, but now mutating the columns only up to the fifth column. 
Continue this procedure: this ends with mutating the last two columns (columns $k-1$ and $k-2$) if $k$ is odd. It ends with mutating the last column if $k$ is even. 
%the last part of mutations are: if $k$ is odd, then mutate the $(k-1)$th column and then the $(k-2)$th column; if $k$ is even, then mutation the $(k-1)$th column. 
\end{remark}

As before, let $n=k+\ell+1$. 
The cluster algebra $\CC[\Gr(k,n,\sim)]$ has an initial seed corresponding to the initial seed of $K_0(\mathscr{C}_{\ell}^{\mathfrak{sl}_k})$ with quiver $Q_{\ell}$, see Section~\ref{subsec:correspondence between modules and tableaux} and 
Figure~\ref{fig:5-9-quiver}. 
Applying the above mutation sequence to the non-frozen part $Q_{\ell}$, we obtain the non-frozen part of $\Gamma_{-2\ell-2}^-$.

We now describe the indecomposable objects of the Grassmannian cluster category 
${\rm CM}(B_{k,n})$ forming a cluster-tiling object whose quiver is $\Gamma_{-2\ell-2}^-$.

To every vertex of $\Gamma_{-2\ell-2}^-$, we give a $k$-subset (or a Pl\"ucker coordinate). Let $(i,m)$ be a vertex of $\Gamma_{-2\ell-2}^-$. For $(i,m)$, we define a $k$-subset (which depends on the height function) 
$J_{i,m}=J_{i,m}^{\xi'}$ as follows: 
\begin{align}\label{eq:k-subsets-for-i_m}
J_{i,m}= [\frac{i-\xi'(i)}{2}, \frac{i-\xi'(i)}{2}+k-i-1] \cup [\frac{i-m-1}{2}+k-i+1, \frac{i-m-1}{2}+k].
\end{align}
Here, we use the convention that the intervals cyclically wrap around, e.g. for $n=9$, the interval $[7,1]$ is the set $\{7,8,9,1\}$. 
We illustrate the $k$-subsets on our running example:
\begin{example}\label{ex:k-subsets}
Let $k=5$ and $\ell=3$, and so $n=9$. The sets $J_{i,m}$ are the following $5$-subsets of $9$: 
\[
\xymatrix@R4pt{
[1,4]\cup[6] & [1,3]\cup[5,6] & [2,3]\cup [5,7] & [2] \cup [4,7] \\
[1,4]\cup[7] & [1,3]\cup[6,7] & [2,3]\cup [6,8] & [2] \cup [5,8] \\
[1,4]\cup[8] & [1,3]\cup[7,8] & [2,3]\cup [7,9] & [2] \cup [6,9]\\
[1,4]\cup[9] & [1,3]\cup[8,9] 
& [2,3]\cup [8,1] & [2] \cup 
[7,1]\\
}
\]
\end{example}

\begin{lemma} \label{lem:plucker coordinate corresponding to Uqslkhat module at position i,m} 
For every vertex $(i,m)$ of $\Gamma_{-2\ell-2}^-$ let 
$L_{i,m}\in {\rm CM}(B_{k,n})$ %$T_{i,m}$
be the rank one module 
with $k$-subset $J_{i,m}$. Then 
\[
L=\bigoplus_{(i,m)\in\Gamma_{-2\ell}^-} L_{i,m}\]
is a cluster tilting object of the stable category $\underline{{\rm CM}}(B_{k,n})$ 
and the quiver of its endomorphism algebra is isomorphic to the non-frozen part of $\Gamma_{-2\ell-2}^-$. 
% For $m \ge -2\ell-2$ (is it also $0\ge m$?), the module $T_{i,m}$ in ${\rm CM}(B_{k,n})$ corresponding to the $U_q(\widehat{\mathfrak{sl}_k})$-module at position $(i,m)$ of the quiver $\Gamma^-_{-2\ell-2}$ is the rank $1$ module associated to the $k$-subset $J_{i,m}$. 
% \begin{align} \label{eq:entries of rank 1 module corresponding to the module at position im}
% [\frac{i-\xi(i)}{2}, \frac{i-\xi(i)}{2}+k-i-1] \cup [\frac{i-m-1}{2}+k-i+1, \frac{i-m-1}{2}+k]. 
% \end{align}
\end{lemma}

\begin{proof}
Let $(i,m)$ be a non-frozen vertex of $\Gamma_{-2\ell-2}^-$. By Remark \ref{rem:mut-sequence}, the mutable part of the quiver $\Gamma^-_{-2\ell-2}$ is mutation equivalent to the mutable part of the quiver $Q_{\ell}$ in Section \ref{subsec:Grassmannian cluster algebras}.  
As discussed in the beginning of this subsection, the $U_q(\widehat{\mathfrak{sl}_k})$-module at position $(i,m)$ 
is $L(Y_{i,m+1}Y_{i,m+3} \cdots Y_{i,\xi'(i)})$, where $\xi': I \to \ZZ$ is the height function: $\xi'(i)=-1$ if $i$ is odd and $\xi'(i)=0$ if $i$ is even (cf. Remark~\ref{rem:height-functions}). {In particular, $\{\chi_q(L(Y_{i,m+1}Y_{i,m+3}\cdots Y_{i,\xi'(i)}))\mid (i,m)\in \Gamma_{-2\ell}^-\}$  are the mutable variables of a cluster of $\CC[\Gr(k,n,\sim)]$ with exchange quiver $\Gamma_{-2\ell-2}^-$. }

Under the isomorphism $\tilde{\Phi}$ of~\cite[Theorem 3.17]{CDFL20}, the module $L(Y_{i,s})$ is mapped to the one-column Young tableau which corresponds to the $k$-subset $[\frac{i-s}{2}, \frac{i-s}{2}+k] \setminus \{ \frac{i-s}{2}+k-i\}$. Using this, one finds that 
$L(Y_{i,m+1}Y_{i,m+3}\cdots Y_{i,\xi'(i)})$ is mapped to the $k$-subset $J_{i,m}$ of (\ref{eq:k-subsets-for-i_m}), 
{so to the module $L_{i,m}$ of ${\rm CM}(B_{k,n})$ by Theorem \ref{thm: cluster variables in Grkn and CMBkn}. 
We conclude that $L$ is a basic cluster-tilting object in $\underline{\rm CM}(B_{k,n})$ by Theorem \ref{t:additve-cat}, and that the quiver of its endomorphism algebra is isomorphic to the non-frozen part of $\Gamma_{-2\ell-2}^-$.}

\end{proof}

For the statement of the theorem, we use the notation of the previous lemma: $(i,m)$ is a vertex of $\Gamma_{-2\ell-2}^-$, $L_{i,m} \in {\rm CM}(B_{k,n})$ is as described in Lemma~\ref{lem:plucker coordinate corresponding to Uqslkhat module at position i,m}, $n=k+\ell+1$.

We also use the following $k$-subset in the theorem: let $(i,m)$ be a vertex of $\Gamma_{-2\ell-2}^-$ and $v$ be a positive integer such that $v\le \frac{m+2\ell+(-1)^{i+1}}{2}$. Then we set 
$I_{i,m}^{(v)}$ to be 

\begin{align} \label{eq:the rank 1 module corresponding to a KR module}
\scalemath{0.9}{
[\frac{i-m+1}{2}, \frac{i-m+1}{2}+k-i-1] \cup [\frac{i-m+2v-1}{2}+k-i+1, \frac{i-m+2v-1}{2}+k]. }
\end{align}

\begin{theorem}\label{thm:generic-kernel-grass} For every vertex $(i,m)$ of $\Gamma^-_{-2\ell-2}$ and every positive integer $v$ such that $v\le \frac{m+2\ell+(-1)^{i+1}}{2}$,
there is a unique indecomposable rigid module 
$M_{v,m}^{(i)}$ 
in ${\rm CM}(B_{k,n})$ such that we have the following triangle in $\underline{\rm CM}(B_{k,n})$ 
\begin{align*}
 L_{i,m} \to  L_{i,m-2v} \to  M_{v,m}^{(i)}\to \tau L_{i,m},
\end{align*}
%where $n=k+\ell+1$, and $T_{i,m}$ is the module in Lemma \ref{lem:plucker coordinate corresponding to Uqslkhat module at position i,m}. 
Moreover, $M_{v,m}^{(i)}$ is the rank one module with $k$-subset 
$I_{i,m}^{(v)}$. 
\end{theorem}

\begin{proof}

{ According to Theorem \ref{thm: cluster variables in Grkn and CMBkn}, there is a one to one correspondence between reachable indecomposable rigid modules in ${\rm CM}(B_{k,n})$ and reachable real prime simple $U_q(\widehat{\mathfrak{sl}_k})$-modules. By Lemma \ref{lem:plucker coordinate corresponding to Uqslkhat module at position i,m} and its proof, the module $L_{i,m}\in {\rm CM}(B_{k,n})$ corresponds to the simple module $L(Y_{i,m+1}Y_{i,m+3}\cdots Y_{i,\xi'(i)})$, which in turn corresponds to the initial cluster variable at the vertex $(i,m)$. On the other hand, one can compute in terms of $U_q(\widehat{\mathfrak{sl}_k})$-modules that $L(Y_{i,m-2v+1}\cdots Y_{i,m-1})$ corresponds to the cluster variable with the $\mathbf{g}$-vector  $e_{i,m-2v}-e_{i,m}\in \mathbb{Z}^{\ell(k-1)}$, where $e_{i,j}$ is the unit vector corresponding to the vertex $(i,j)$. Consequently, by Theorem \ref{thm: cluster variables in Grkn and CMBkn}, there is a unique indecomposable rigid module $M_{v,m}^{(i)}$ corresponding to $L(Y_{i,m-2v+1}\cdots Y_{i,m-1})$.
Furthermore, by Theorem \ref{t:additve-cat} and the definition of the index, we obtain a distinguished triangle\[L_{i,m} \to  L_{i,m-2v} \to  M_{v,m}^{(i)}\to \tau L_{i,m}
\]
in the stable category $\underline{{\rm CM}}(B_{k,n})$.
}
% In terms of 
% $U_q(\widehat{\mathfrak{sl}_k})$-modules: at position $(i,m)$, we have  $L(Y_{i,m+1}Y_{i,m+3}\cdots Y_{i,\xi'(i)})$ and at position $(i,m-2v)$, we have 
% $L(Y_{i,m-2v+1} \cdots Y_{i,m-1}Y_{i,m+1}\cdots Y_{i,\xi'(i)})$ (Section~\ref{subsec:Generic kernels via Grassmannian cluster categories}). Using their ${\mathbf g}$-vectors, one finds that  the third module $M_{v,m}^{(i)}$ in the triangle corresponds to 
% $L(Y_{i,m-2v+1}\cdots Y_{m-1})$. 

As argued in the proof of Lemma~\ref{lem:plucker coordinate corresponding to Uqslkhat module at position i,m}, these translate to the following 
modules in $\underline{\rm CM}(B_{k,n})$: 
the module $L_{i,m}$ is the rank 1-module with $k$-subset $J_{i,m}$ from (\ref{eq:k-subsets-for-i_m}), the module 
$L_{i,m-2v}$ is the rank 1-module with $k$-subset $J_{i,m-2v}$ and 
the module 
$M_{v,m}^{(i)}$ is the rank 1-module with $k$-subset $I_{i,m}^{(v)}$ from (\ref{eq:the rank 1 module corresponding to a KR module}).

% The $U_q(\widehat{\mathfrak{sl}_k})$-module at position $(i,m)$ is $L(Y_{i,m+1}Y_{i,m+3} \cdots Y_{i,\xi(i)})$ and the $U_q(\widehat{\mathfrak{sl}_k})$-module at position $(i,m-2v)$ is $L(Y_{i,m-2v+1} \cdots Y_{i,m-1}Y_{i,m+1}\cdots Y_{i,\xi(i)})$. Therefore the $U_q(\widehat{\mathfrak{sl}_k})$-module corresponding to $M_{v,m}^{(i)}$ is $L(Y_{i,m-2v+1} \cdots Y_{i,m-1})$. 
% As argued in the proof of Lemma~\ref{lem:plucker coordinate corresponding to Uqslkhat module at position i,m}, the module $L(Y_{i,s})$ corresponds to the 
% Pl\"{u}cker coordinate with entries $[\frac{i-s}{2}, \frac{i-s}{2}+k] \setminus \{ \frac{i-s}{2}+k-i \}$ and the module $L(Y_{i,m-2v+1} \cdots Y_{i,m-1})$ corresponds to the Pl\"{u}cker coordinate $P_{J}$, where $J$ is given by (\ref{eq:the rank 1 module corresponding to a KR module}), 
% using \cite[Theorem 3.17]{CDFL20})
% .
% The rank one modules in ${\rm CM}(B_{k,n})$ are in bijection with Pl\"{u}cker coordinates in $\mathbb{C}[\Gr(k,n)]$, \cite[\S 5]{JKS16}.  Therefore $M_{v,m}^{(i)}$ is a rank one module whose profile has entries given by (\ref{eq:the rank 1 module corresponding to a KR module}).
\end{proof}

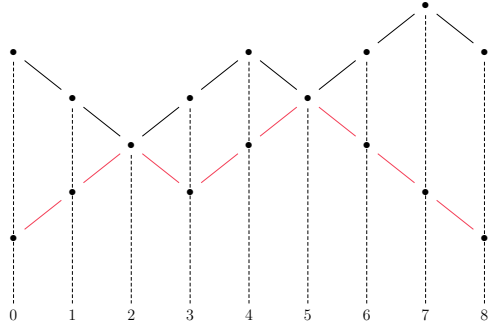
\begin{figure}
    \centering
    \adjustbox{scale=0.5}{
    \begin{tikzcd}
	&&&&&&& \bullet \\
	\bullet &&&& \bullet && \bullet && \bullet \\
	& \bullet && \bullet && \bullet \\
	&& \bullet && \bullet && \bullet \\
	& \bullet && \bullet &&&& \bullet \\
	\bullet &&&&&&&& \bullet \\
	\\
	0 & 1 & 2 & 3 & 4 & 5 & 6 & 7 & 8 
% 	\arrow["{x_1}", color={rgb,255:red,244;green,62;blue,86}, from=2-1, to=3-2]
	\arrow[no head, color={rgb,255:red,244;green,62;blue,86}, from=6-1, to=5-2]
	\arrow[no head, color={rgb,255:red,244;green,62;blue,86}, from=5-2, to=4-3]
	\arrow[no head, color={rgb,255:red,244;green,62;blue,86},from=4-3, to=5-4]
	\arrow[no head, color={rgb,255:red,244;green,62;blue,86}, from=5-4, to=4-5]
	\arrow[no head, color={rgb,255:red,244;green,62;blue,86}, from=4-5, to=3-6]
	\arrow[no head, color={rgb,255:red,244;green,62;blue,86}, from=3-6, to=4-7]
	\arrow[no head, color={rgb,255:red,244;green,62;blue,86}, from=4-7, to=5-8]
	\arrow[no head, color={rgb,255:red,244;green,62;blue,86}, from=5-8, to=6-9]
	\arrow[no head, from=2-1, to=3-2]
	\arrow[no head, from=3-2, to=4-3]
	\arrow[no head, from=4-3, to=3-4]
	\arrow[no head, from=3-4, to=2-5]
	\arrow[no head, from=2-5, to=3-6]
	\arrow[no head, from=3-6, to=2-7]
	\arrow[no head, from=2-7, to=1-8]
	\arrow[no head, from=1-8, to=2-9]
	\arrow[dashed, no head, from=2-1, to=8-1]
	\arrow[dashed, no head, from=3-2, to=8-2]
	\arrow[dashed, no head, from=4-3, to=8-3]
	\arrow[dashed, no head, from=3-4, to=8-4]
	\arrow[dashed, no head, from=2-5, to=8-5]
	\arrow[dashed, no head, from=3-6, to=8-6]
	\arrow[dashed, no head, from=2-7, to=8-7]
	\arrow[dashed, no head, from=1-8, to=8-8]
	\arrow[dashed, no head, from=2-9, to=8-9]
\end{tikzcd}
    }
    \caption{The top rim represents the rank-1 module with $4$-subset $1258$. The 
    $4$-subset of $\tau^{-1}$ of this module is $3678$, indicated by the second rim (in red). %The red rim below represents the module $3678$: We have that $\tau (3678) = 1258$.
    }
    \label{fig:tau3678 is 1258}
\end{figure}

One can determine the image of a rank one module $N$ in ${\rm CM}(B_{k,n})$ under $\tau^{-1}$ using the first syzygy of the module \cite[\S 2]{BB17}. Combinatorially, it amounts to subtracting the lattice diagram of $N$ from the one of its projective cover (see Section~\ref{subsec:grassmannian cluster categories} for the definition of lattice diagrams). In particular, if the $k$-subset of $N$ consists of two intervals, its syzygy is also a rank one module with a two-interval $k$-subset. In that situation, the $k$-subset of $\tau^{-1}(N)$ is obtained by extending the rim of $N$ beyond the two lowest points, 
see Figure~\ref{fig:tau3678 is 1258} for an example. 

Using this, one finds that the module $\tau M_{v,m}^{(i)}$ is also of rank one and that its $k$-subset is 
\begin{align}\label{eq:tauM case 1}
[\frac{1-i-m}{2}, \frac{i-m-1}{2}] \cup [\frac{i-m+2v+1}{2}, \frac{i-m+2v-1}{2}+k-i];
\end{align} 
if $\frac{1-i-m}{2}\ge 1$, respectively 
\begin{align} \label{eq:tauM case 2}
[1, \frac{i-m-1}{2}] \cup [\frac{i-m+2v+1}{2}, \frac{i-m+2v-1}{2}+k-i] \cup [n+\frac{1-i-m}{2}, n]. 
\end{align}
if $\frac{1-i-m}{2} \le 0$ (in the second case, the three intervals are in fact just two intervals as the first and the last one form a single interval including the elements $n$ and $1$). 

\begin{example}\label{ex:illustrate-M_vm} 
Consider the vertex $(i,m)=(3,-2)$ of 
$\Gamma_{-6-2}^-$ (i.e. $\ell=3$), let $k=4$ and $n=k+\ell+1=8$ and choose $v=2$.
Then the $U_q(\widehat{\mathfrak{sl}_4})$-module at position $(i,m)$ of $\Gamma^-_{-8}$ is $L(Y_{3,-1})$ and 
the module at position $(i,m-2v)$ is $L(Y_{3,-5}Y_{3,-3}Y_{3,-1})$. 
We have 
\begin{eqnarray*}
    J_{3,-2} & = & [2]\cup [4,6]\\
    J_{3,-6} & = & [2]\cup[6,8]\\
    I_{3,-2}^{(2)} & = & [3]\cup [6,8]
\end{eqnarray*}
and $\tau M_{v,m}^{(i)}$ is the module 
with $4$-subset $[1,2]\cup[5]\cup[8]$ as in 
(\ref{eq:tauM case 2}), see Figure \ref{fig:tau3678 is 1258} for the last two $4$-subsets. 

% The corresponding modules in ${\rm CM}(B_{4,8})$ are $L_{i,m}$ with $4$-subset $2456$ and $L_{i,m-2v}$ with $4$-subset $2678$ respectively. The module $M_{v,m}^{(i)}$ has $4$-subset $3678$ and $\tau M_{v,m}^{(i)}$ has $4$-subset $1258$, see Figure \ref{fig:tau3678 is 1258}. 
\end{example}

\begin{theorem} \label{thm:realization-generic-kernel}
Let $(i,m)$ be a vertex of $\Gamma^-$ and $v\geq 1$. Let $\ell\in \mathbb{N}$ such that $-2\ell\le m-2v-2$.
% Let $s=-2u\in \mathbb{Z}_{<0}$ such that $s \ll 0$.
Then $\Hom_{\underline{\rm CM}(B_{k,k+l+1})}(L,\tau M_{v,m}^{(i)})$ is the kernel of a generic $A$-module homomorphism from $I(i,m)$ to $I(i,m-2v)$. 
\end{theorem}

\begin{proof}
Denote by $s=-2\ell$.
According to \cite[Proposition 4.17]{HL16}, the potential of $(\Gamma_s^-,W_s)$ is rigid. 
Therefore $\mathcal{C}(s)\cong \underline{\rm CM}(B_{k,k+\ell+1})$ by \cite[Theorem 0.2]{CZ22}. Now the result follows from {Theorem} \ref{thm:generic-kernel-grass} and the discussion after Definition \ref{def:generic-kernel}.
\end{proof}

\begin{remark}\label{rem:criterion-compatibility}
Let $x_{v_1,m_1}^{(i)}$ and $x_{v_2,m_2}^{(j)}$ be the cluster variables of $K_0(\mathscr{C}_{\ell}^{\mathfrak{sl}_k})$ corresponding to $K_{v_1,m_1}^{(i)}$ and $K_{v_2,m_2}^{(j)}$ respectively. Combining Theorem~\ref{thm:realization-generic-kernel} with \cite[Theorem 3.1]{BB17} yields a combinatorial criterion for checking whether $x_{v_1,m_1}^{(i)}$ and $x_{v_2,m_2}^{(j)}$ belong to a common cluster, providing an alternative approach to that of \cite{oh2019simplicity, oh2025DenominatorsRmatrices}.  
By Theorem \ref{thm:realization-generic-kernel}, we have 
%$K_{v_1,m_1}^{(i)}\cong \Hom_{\underline{\rm CM}(B_{k,n})}(L,\tau M_{v_1,m_1}^{(i)})$ and similarly for $K_{v_2,m_2}^{(j)})$ and hence  
$x_{v_1,m_1}^{(i)}=\mathbf{X}^L_{M_{v_1,m_1}^{(i)}}$ and $x_{v_2,m_2}^{(j)}=\mathbf{X}^L_{M_{v_2,m_2}^{(j)}}$.
According to Corollary \ref{c:prod-cluster-monomial}, we clearly have that $x_{v_1,m_1}^{(i)}$ and $x_{v_2,m_2}^{(j)}$ belong to a cluster if and only if $\Ext^1_{{\rm CM}(B_{k,n})}(M_{v_1,m_1}^{(i)},M_{v_2,m_2}^{(j)})=0$, while the latter can be computed by \cite[Theorem 3.1]{BB17}. 
\end{remark}

%%%%%%%%
%
\section{Application to construction of indecomposable modules} \label{sec:application to construction of indecomposable modules}

Apart from the finite types (i.e. for $(k,n)\in \{(2,n)\}\cup \{(3,6),(3,7),(3,8)\}$), the Grassmannian cluster categories ${\rm CM}(B_{k,n})$ have infinitely 
many rigid indecomposable modules. We note that in the infinite types, there are rigid and non-rigid indecomposable modules {of} arbitrary ranks, so a complete characterisation seems very hard. However, our methods apply to any $(k,n)$. 
Also, it is known that any connected component in the Auslander-Reiten quiver is isomorphic to a tube (\cite[Proposition 2.1]{BBGL21}) and even more can be said in the two tame cases $(3,9)$ and $(4,8)$. The rank 2 rigid indecomposable modules in these two cases are determined in~\cite[\S 6,\S 7]{BBG20} and in~\cite[\S 4]{BBGL21}. 
For $(3,9)$,~\cite[\S 6.3]{BBGL21} determines the rank 3 rigid indecomposables.

In this section, we explicitly construct rigid and non-rigid indecomposable modules of Grassmannian cluster categories ${\rm CM}(B_{k,n})$ in the two tame cases. 
{We first use the mutation rule for tableaux to obtain cluster variables in $\CC[\Gr(k,n)]$, cf. subsection \ref{ss:mutation-tableaux} or \cite{CDFL20}.} Then we apply the results from Section~\ref{sec:correspondence between additive and monoid categorifications} to construct the corresponding rigid indecomposable modules in ${\rm CM}(B_{k,n})$. 

In order to construct non-rigid indecomposable modules in ${\rm CM}(B_{k,n})$, we first need to construct prime non-real elements in $\CC[\Gr(k,n)]$. {From these, we will then obtain non-rigid  modules 
in ${\rm CM}(B_{k,n})$, which are conjecturally  indecomposable.} 

The computations in this section can be found in \url{https://github.com/lijr07/from-tableaux-to-profiles}.

% Recall that a simple $U_q(\widehat{\mathfrak{g}})$-module $L(m)$ is called real if $L(m) \otimes L(m)$ is still simple \cite{Lec03}. A simple module $L(m)$ is called prime if $L(m) \not\cong L(m')\otimes L(m'')$ for any non-trivial $U_q(\widehat{\mathfrak{g}})$-modules $L(m')$, $L(m'')$, \cite{CP97}. A tableau $T$ is called real (resp. prime) if $L(M_T)$ is real (resp. prime), \cite{CDFL20}. A simple module $L(M)$ is real if and only if $\chi_q(L(M))\chi_q(L(M)) = \chi_q(L(M^2))$. A simple module $L(M)$ is prime if and only if $\chi_q(L(M)) \ne \chi_q(L(M'))\chi_q(L(M''))$ for any non-trivial $U_q(\widehat{\mathfrak{sl}_k})$-modules $L(M'), L(M'')$.

% \section{Reduction for the categories of modules of quantum affine algebras}

% First we introduce reduction method in quantum affine algebra setting. Then we translate the proof in $\Gr(3,9)$, $\Gr(4,8)$ cases in ${\rm CM}(B_{k,n})$ setting to $\mathcal{C}_{\ell}^{\mathfrak{sl}_k}$ and $\CC[\Gr(k,n)]$ setting. 

%\section{$f$-compatibility of dual canonical (triangular?) basis of cluster algebras}

% \subsection{One to one correspondence between rigid indecomposable modules and real prime tableaux in the cases of $\Gr(3,9)$ and $\Gr(4,8)$}

% \begin{theorem}
% Let $(k,n) \in \{(3,9), (4,8)\}$. Then there is a one to one correspondence between rigid indecomposable modules in ${\rm CM}(B_{k,n})$ and reachable real prime tableaux in $\SSYT(k,[n])$. 
% \end{theorem}

% \begin{proof}

% \end{proof}

\subsection{Rigid indecomposable modules in ${\rm CM}(B_{3,9})$ up to rank 4} \label{subsec:rigid indecomposable modules Gr39}

The \textit{rank of a tableau} is the number of {its columns.} 
%If $x$ is a cluster variable in $\CC[\Gr(k,n)]$ with corresponding to a tableau $\bT$, the \textit{rank of the cluster variable $x$} is defined to be the rank of $\bT$. 

In order to save space, we will only describe the tableaux up to {{\em promotion}, an operation on the set } of semistandard Young tableaux defined in terms of ``jeu de taquin'' sliding moves \cite{Sch63, Sch72, Sch77}. By~\cite{Ga80}, it can be described in terms of so-called Bender-Knuth involutions~\cite{BK72}. 
We recall these involutions: 

\begin{definition} 
Consider $\SSYT(k,[n])$, let $i\in [n]$. \\
(1) 
The $i$th Bender-Knuth involution 
\[{\rm BK}_i: \SSYT(k, [n]) \to \SSYT(k, [n]),\] 
is defined as follows: {for each column that does not contain both $i$ and $i+1$, replace $i$ by $i+1$ and $i+1$ by $i$.} Then one reorders every row such that the result is a semistandard Young tableau.
This process is clearly involutive. 

(2) For $\bT\in \SSYT(k,[n])$, the promotion 
${\rm pr}(\bT)$ is then defined as 
\begin{align*}
{\rm pr}(\bT) = {\rm BK}_1 \circ \cdots \circ {\rm BK}_{n-1}(\bT).
\end{align*} 
% the following procedure: for $i$ and $i+1$ which are not in the same column, we replace $i$ by $i+1$ and replace $i+1$ by $i$, and then reorder $i$, $i+1$ in each row such that the resulting tableau is semistandard. The promotion ${\rm pr}(\bT)$ of $\bT$ is defined by 
% \begin{align*}
% {\rm pr}(\bT) = {\rm BK}_1 \circ \cdots \circ {\rm BK}_{n-1}(\bT).
% \end{align*} 
\end{definition}
For example, when $n=8$, for $\bT= \scalemath{0.6}{
\begin{ytableau}
1 & 1  \\
2 & 4  \\
4 & 7 
\end{ytableau}}$, we have ${\rm pr}(\bT) = 
\scalemath{0.6}{\begin{ytableau}
2 & 2  \\
3 & 5  \\
5 & 8 
\end{ytableau} }$. 
% \begin{align*}
% & \bT= \scalemath{0.6}{
% \begin{ytableau}
% 1 & 1  \\
% 2 & 4  \\
% 4 & 7 
% \end{ytableau}}  
% \quad 
% {\rm pr}(\bT) = 
% \scalemath{0.6}{\begin{ytableau}
% 2 & 2  \\
% 3 & 5  \\
% 5 & 8 
% \end{ytableau} }
% \end{align*}

In $\CC[\Gr(3,9)]$, there are $168, 225$ and $288$ cluster variables with tableaux of rank $2,3$ and $4$. 
In~\cite{CDHHHL}, the authors explained how to obtain these. 
The corresponding rank $2$ and rank $3$ {rigid indecomposable modules} are determined in~\cite{BBG20, BBGL21}. So it remains to describe cluster variables whose tableaux have four columns, i.e. rank 4 modules. 
Up to tableaux promotion, the tableaux of rank $4$ in $\CC[\Gr(3,9)]$ are the following: 

\begin{align*}
& \scalemath{0.6}{
\begin{ytableau}
1 & 1 & 3 & 3 \\
2 & 2 & 6 & 7 \\
4 & 5 & 8 & 9
\end{ytableau}, \ \begin{ytableau}
1 & 1 & 3 & 4 \\
2 & 2 & 6 & 7 \\
4 & 5 & 8 & 9
\end{ytableau}, \ \begin{ytableau}
1 & 1 & 3 & 3 \\
2 & 4 & 5 & 6 \\
4 & 7 & 8 & 9
\end{ytableau}, \ \begin{ytableau}
1 & 1 & 2 & 3 \\
2 & 5 & 5 & 6 \\
4 & 7 & 8 & 9
\end{ytableau}, \ \begin{ytableau}
1 & 1 & 2 & 3 \\
2 & 5 & 6 & 6 \\
4 & 7 & 8 & 9
\end{ytableau}, \ \begin{ytableau}
1 & 1 & 2 & 3 \\
2 & 5 & 6 & 7 \\
4 & 7 & 8 & 9
\end{ytableau}, \ \begin{ytableau}
1 & 1 & 3 & 3 \\
2 & 5 & 5 & 6 \\
4 & 7 & 8 & 9
\end{ytableau}, }
\end{align*}
\begin{align*}
& \scalemath{0.6}{ \begin{ytableau}
1 & 1 & 3 & 3 \\
2 & 5 & 6 & 6 \\
4 & 7 & 8 & 9
\end{ytableau}, \  \begin{ytableau}
1 & 1 & 3 & 3 \\
2 & 5 & 6 & 7 \\
4 & 7 & 8 & 9
\end{ytableau}, \ \begin{ytableau}
1 & 1 & 3 & 4 \\
2 & 5 & 6 & 7 \\
4 & 7 & 8 & 9
\end{ytableau}, \ \begin{ytableau}
1 & 1 & 2 & 3 \\
2 & 5 & 6 & 7 \\
4 & 8 & 8 & 9
\end{ytableau}, \ \begin{ytableau}
1 & 1 & 2 & 3 \\
2 & 5 & 6 & 7 \\
4 & 8 & 9 & 9
\end{ytableau}, \ \begin{ytableau}
1 & 1 & 3 & 3 \\
2 & 5 & 6 & 7 \\
4 & 8 & 8 & 9
\end{ytableau}, \ \begin{ytableau}
1 & 1 & 3 & 4 \\
2 & 5 & 6 & 7 \\
4 & 8 & 8 & 9
\end{ytableau}, } 
\end{align*}
\begin{align*}
& \scalemath{0.6}{ \begin{ytableau}
1 & 1 & 3 & 3 \\
2 & 5 & 6 & 7 \\
4 & 8 & 9 & 9
\end{ytableau}, \ \begin{ytableau}
1 & 1 & 3 & 4 \\
2 & 5 & 6 & 7 \\
4 & 8 & 9 & 9
\end{ytableau}, \ \begin{ytableau}
1 & 2 & 3 & 3 \\
2 & 4 & 5 & 6 \\
4 & 7 & 8 & 9
\end{ytableau}, \ \begin{ytableau}
1 & 2 & 3 & 3 \\
2 & 5 & 5 & 6 \\
4 & 7 & 8 & 9
\end{ytableau}, \ \begin{ytableau}
1 & 2 & 3 & 3 \\
2 & 5 & 6 & 6 \\
4 & 7 & 8 & 9
\end{ytableau}, \ \begin{ytableau}
1 & 2 & 3 & 3 \\
2 & 5 & 6 & 7 \\
4 & 7 & 8 & 9
\end{ytableau}, \ \begin{ytableau}
1 & 2 & 3 & 3 \\
2 & 5 & 6 & 7 \\
4 & 8 & 8 & 9
\end{ytableau}, } 
\end{align*}
\begin{align*}
& \scalemath{0.6}{ \begin{ytableau}
1 & 2 & 3 & 4 \\
2 & 5 & 6 & 7 \\
4 & 8 & 8 & 9
\end{ytableau}, \ \begin{ytableau}
1 & 2 & 3 & 3 \\
2 & 5 & 6 & 7 \\
4 & 8 & 9 & 9
\end{ytableau}, \ \begin{ytableau}
1 & 2 & 3 & 4 \\
2 & 5 & 6 & 7 \\
4 & 8 & 9 & 9
\end{ytableau}, \ \begin{ytableau}
1 & 3 & 3 & 4 \\
2 & 5 & 6 & 7 \\
4 & 8 & 8 & 9
\end{ytableau}, \ \begin{ytableau}
1 & 3 & 3 & 4 \\
2 & 5 & 6 & 7 \\
4 & 8 & 9 & 9
\end{ytableau}, \ \begin{ytableau}
1 & 1 & 3 & 4 \\
2 & 2 & 6 & 7 \\
5 & 5 & 8 & 9
\end{ytableau}, \ \begin{ytableau}
1 & 1 & 3 & 4 \\
2 & 2 & 7 & 7 \\
5 & 6 & 8 & 9
\end{ytableau}, } 
\end{align*}
\begin{align*}
& \scalemath{0.6}{ \begin{ytableau}
1 & 1 & 3 & 4 \\
2 & 3 & 6 & 7 \\
5 & 6 & 8 & 9
\end{ytableau}, \ \begin{ytableau}
1 & 1 & 3 & 4 \\
2 & 3 & 7 & 7 \\
5 & 6 & 8 & 9
\end{ytableau}, \ \begin{ytableau}
1 & 1 & 4 & 4 \\
2 & 3 & 7 & 7 \\
5 & 6 & 8 & 9
\end{ytableau}, \ \begin{ytableau}
1 & 2 & 3 & 4 \\
2 & 5 & 6 & 7 \\
5 & 8 & 8 & 9
\end{ytableau}, \ \begin{ytableau}
1 & 2 & 3 & 4 \\
2 & 5 & 6 & 7 \\
5 & 8 & 9 & 9
\end{ytableau}, \ \begin{ytableau}
1 & 2 & 3 & 4 \\
2 & 6 & 6 & 7 \\
5 & 8 & 9 & 9
\end{ytableau}. }
\end{align*}

Using %the correspondence between monoidal categorifications and additive categorifications described in 
Section~\ref{sec:profiles}, we get the following rigid 
indecomposable modules {from the above tableaux:} 

\begin{align*}
& \scalemath{0.7}{ \ffrac{1 3 7}{1 3 6}{2 5 9}{2 4 8}, \ 
 \ffrac{1 4 7}{1 3 6}{2 5 9}{2 4 8}, \ 
 \ffrac{1 3 6}{3 5 9}{2 4 8}{1 4 7}, \ 
 \ffrac{1 3 6}{2 5 9}{2 5 8}{1 4 7}, \ 
 \ffrac{1 3 6}{2 6 9}{2 5 8}{1 4 7}, \ 
 \ffrac{1 3 7}{2 6 9}{2 5 8}{1 4 7}, \ 
 \ffrac{1 3 6}{3 5 9}{2 5 8}{1 4 7}, \ 
 \ffrac{1 3 6}{3 6 9}{2 5 8}{1 4 7}, \ 
 \ffrac{1 3 7}{3 6 9}{2 5 8}{1 4 7}, \ 
 \ffrac{1 4 7}{3 6 9}{2 5 8}{1 4 7}, \ 
 \ffrac{1 3 7}{2 6 9}{2 5 8}{1 4 8}, \ 
 \ffrac{1 3 7}{2 6 9}{2 5 9}{1 4 8}, \ 
 \ffrac{1 3 7}{3 6 9}{2 5 8}{1 4 8}, \ 
 \ffrac{1 4 7}{3 6 9}{2 5 8}{1 4 8}, \ 
 \ffrac{1 3 7}{3 6 9}{2 5 9}{1 4 8}, \ 
 \ffrac{1 4 7}{3 6 9}{2 5 9}{1 4 8}, \ 
 \ffrac{1 3 6}{3 5 9}{2 4 8}{2 4 7},}
 \end{align*}
\begin{align*}
& \scalemath{0.7}{
 \ffrac{1 3 6}{3 5 9}{2 5 8}{2 4 7}, \ 
 \ffrac{1 3 6}{3 6 9}{2 5 8}{2 4 7}, \ 
 \ffrac{1 3 7}{3 6 9}{2 5 8}{2 4 7}, \ 
 \ffrac{1 3 7}{3 6 9}{2 5 8}{2 4 8}, \ 
 \ffrac{1 4 7}{3 6 9}{2 5 8}{2 4 8}, \ 
 \ffrac{1 3 7}{3 6 9}{2 5 9}{2 4 8}, \ 
 \ffrac{1 4 7}{3 6 9}{2 5 9}{2 4 8}, \ 
 \ffrac{1 4 7}{3 6 9}{3 5 8}{2 4 8}, \ 
 \ffrac{1 4 7}{3 6 9}{3 5 9}{2 4 8}, \ 
 \ffrac{1 4 7}{1 3 6}{2 5 9}{2 5 8}, \ 
 \ffrac{1 4 7}{1 3 7}{2 6 9}{2 5 8}, \ 
 \ffrac{1 4 7}{1 3 6}{3 6 9}{2 5 8}, \ 
 \ffrac{1 4 7}{1 3 7}{3 6 9}{2 5 8}, \ 
 \ffrac{1 4 7}{1 4 7}{3 6 9}{2 5 8}, \ 
 \ffrac{1 4 7}{3 6 9}{2 5 8}{2 5 8}, \ 
 \ffrac{1 4 7}{3 6 9}{2 5 9}{2 5 8}, \ 
 \ffrac{1 4 7}{3 6 9}{2 6 9}{2 5 8}. }
\end{align*}

\begin{remark}\label{rem:cyclic-shifts}
The counterpart of the promotion operation on Young tableaux on profiles (or filtration factors) of modules in $\rm CM(B_{k,n})$ is to add a fixed number $a\in [n]$ to every entry of the profile. 

We say that two profiles $P_M=I_1\mid \cdots I_m$ 
%\begin{small}\begin{array}{c} I_1\\ \hline \vdots \\ \hline I_m \end{array}\end{small}$ 
and 
$P_N=J_1\mid \cdots J_m$ 
%$P_N=\begin{small}\begin{array}{c} J_1\\ \hline \vdots \\ \hline J_m \end{array}\end{small}$ 
are \textit{equal up to cyclic shift} if 
if there exists $a\in[n]$ such that $J_r=\{i+a\mid i\in I_r\}$ for all $r=1,\dots, m$ (reducing modulo $n$). 
\end{remark}
We conjecture that the 34 profile above form a complete list of representatives of rank $4$ rigid indecomposables, up to cyclic shift. 

\subsection{Rigid indecomposable modules in ${\rm CM}(B_{4,8})$ up to rank 4} \label{subsec:rigid indecomposable modules Gr48}
 
In $\CC[\Gr(4,8)]$, there are $120$, $174$ and $208$ cluster variables with tableaux of rank $2,3$ and $4$ respectively, \cite{CDHHHL}. 
The rank $2$ rigid indecomposable modules in ${\rm CM}(B_{4,8})$ are studied in detail in \cite{BBG20}. So it remains to describe cluster variables in $\CC[\Gr(4,8)]$ whose tableaux have rank $3$ and $4$. Up to promotion, the tableaux with three columns are the following: 

\begin{align*}
& \scalemath{0.6}{ \begin{ytableau}
1 & 1 & 3 \\
2 & 2 & 4 \\
3 & 4 & 7 \\
5 & 6 & 8
\end{ytableau}, \ \begin{ytableau}
1 & 1 & 3 \\
2 & 2 & 5 \\
3 & 4 & 7 \\
5 & 6 & 8
\end{ytableau}, \ \begin{ytableau}
1 & 1 & 2 \\
2 & 3 & 4 \\
3 & 6 & 6 \\
5 & 7 & 8
\end{ytableau}, \ \begin{ytableau}
1 & 1 & 2 \\
2 & 3 & 4 \\
3 & 6 & 7 \\
5 & 7 & 8
\end{ytableau}, \ \begin{ytableau}
1 & 1 & 2 \\
2 & 3 & 4 \\
3 & 6 & 7 \\
5 & 8 & 8
\end{ytableau}, \ \begin{ytableau}
1 & 1 & 2 \\
2 & 4 & 4 \\
3 & 5 & 6 \\
5 & 7 & 8
\end{ytableau}, \ \begin{ytableau}
1 & 1 & 2 \\
2 & 4 & 4 \\
3 & 6 & 6 \\
5 & 7 & 8
\end{ytableau}, \ \begin{ytableau}
1 & 1 & 2 \\
2 & 4 & 4 \\
3 & 6 & 7 \\
5 & 7 & 8
\end{ytableau}, } 
\end{align*}
\begin{align*}
& \scalemath{0.6}{ \begin{ytableau}
1 & 1 & 2 \\
2 & 4 & 5 \\
3 & 6 & 6 \\
5 & 7 & 8
\end{ytableau}, \ \begin{ytableau}
1 & 1 & 2 \\
2 & 4 & 5 \\
3 & 6 & 7 \\
5 & 7 & 8
\end{ytableau}, \ \begin{ytableau}
1 & 1 & 3 \\
2 & 4 & 5 \\
3 & 6 & 7 \\
5 & 7 & 8
\end{ytableau}, \ \begin{ytableau}
1 & 1 & 2 \\
2 & 4 & 4 \\
3 & 6 & 7 \\
5 & 8 & 8
\end{ytableau}, \ \begin{ytableau}
1 & 1 & 2 \\
2 & 4 & 5 \\
3 & 6 & 7 \\
5 & 8 & 8
\end{ytableau}, \ \begin{ytableau}
1 & 1 & 3 \\
2 & 4 & 4 \\
3 & 6 & 7 \\
5 & 8 & 8
\end{ytableau}, \ \begin{ytableau}
1 & 1 & 3 \\
2 & 4 & 5 \\
3 & 6 & 7 \\
5 & 8 & 8
\end{ytableau}, \ \begin{ytableau}
1 & 2 & 3 \\
2 & 4 & 4 \\
3 & 6 & 7 \\
5 & 8 & 8
\end{ytableau}, } 
\end{align*}
\begin{align*}
& \scalemath{0.6}{ \begin{ytableau}
1 & 2 & 3 \\
2 & 4 & 5 \\
3 & 6 & 7 \\
5 & 8 & 8
\end{ytableau}, \ \begin{ytableau}
1 & 1 & 3 \\
2 & 2 & 5 \\
3 & 4 & 7 \\
6 & 6 & 8
\end{ytableau}, \ \begin{ytableau}
1 & 2 & 3 \\
2 & 4 & 5 \\
3 & 6 & 7 \\
6 & 7 & 8
\end{ytableau}, \ \begin{ytableau}
1 & 2 & 3 \\
2 & 4 & 5 \\
3 & 6 & 7 \\
6 & 8 & 8
\end{ytableau}, \ \begin{ytableau}
1 & 2 & 3 \\
2 & 4 & 5 \\
3 & 7 & 7 \\
6 & 8 & 8
\end{ytableau}, \ \begin{ytableau}
1 & 1 & 3 \\
2 & 2 & 5 \\
4 & 4 & 7 \\
6 & 6 & 8
\end{ytableau}, \ \begin{ytableau}
1 & 1 & 3 \\
2 & 2 & 5 \\
4 & 5 & 7 \\
6 & 7 & 8
\end{ytableau}, \ \begin{ytableau}
1 & 1 & 3 \\
2 & 3 & 5 \\
4 & 5 & 7 \\
6 & 7 & 8
\end{ytableau}, \ \begin{ytableau}
1 & 2 & 3 \\
2 & 4 & 5 \\
4 & 6 & 7 \\
6 & 8 & 8
\end{ytableau}. }
\end{align*}

Using the correspondence between monoidal categorifications and additive categorifications described in Section \ref{sec:correspondence between additive and monoid categorifications}, we are able to construct the rigid indecomposable modules corresponding to the tableaux above. The profiles of the corresponding rank $3$ rigid indecomposable modules in ${\rm CM}(B_{4,8})$ are the following.
\begin{align*}
\scalemath{0.7}{ \thfrac{1347}{1246}{2358}, \ 
 \thfrac{1357}{1246}{2358}, \ 
 \thfrac{1246}{2368}{1357}, \ 
 \thfrac{1247}{2368}{1357}, \ 
 \thfrac{1247}{2368}{1358}, \ 
 \thfrac{1246}{2458}{1357}, \ 
 \thfrac{1246}{2468}{1357}, \ 
 \thfrac{1247}{2468}{1357}, \
 \thfrac{1256}{2468}{1357}, \ 
 \thfrac{1257}{2468}{1357}, \ 
 \thfrac{1357}{2468}{1357}, \ 
 \thfrac{1247}{2468}{1358}, \ 
 \thfrac{1257}{2468}{1358}, }
\end{align*}
\begin{align*}
\scalemath{0.7}{
\thfrac{1347}{2468}{1358}, \ 
 \thfrac{1357}{2468}{1358}, \ 
 \thfrac{1347}{2468}{2358}, \
\thfrac{1357}{2468}{2358}, \ 
 \thfrac{1357}{1246}{2368}, \ 
 \thfrac{1357}{2468}{2367}, \ 
 \thfrac{1357}{2468}{2368}, \ 
 \thfrac{1357}{2478}{2368}, \ 
 \thfrac{1357}{1246}{2468}, \ 
 \thfrac{1357}{1257}{2468}, \ 
 \thfrac{1357}{1357}{2468}, \ 
 \thfrac{1357}{2468}{2468}. }
\end{align*} 
{We conjecture that, up to cyclic shifts, these are all the rigid indecomposable modules of rank $3$ in ${\rm CM}(B_{4,8})$ (cf. Remark~\ref{rem:cyclic-shifts}).}

Up to promotion, the tableaux with four columns are the following:  
%Rank $4$ cluster variables in $\CC[\Gr(4,8)]$ are:
\begin{align*}
& \scalemath{0.6}{ \begin{ytableau}
1 & 1 & 1 & 3 \\
2 & 2 & 4 & 5 \\
3 & 4 & 6 & 7 \\
5 & 6 & 7 & 8
\end{ytableau}, \ \begin{ytableau}
1 & 1 & 1 & 3 \\
2 & 2 & 4 & 5 \\
3 & 4 & 6 & 7 \\
5 & 6 & 8 & 8
\end{ytableau}, \ \begin{ytableau}
1 & 1 & 1 & 3 \\
2 & 2 & 4 & 5 \\
3 & 4 & 7 & 7 \\
5 & 6 & 8 & 8
\end{ytableau}, \ \begin{ytableau}
1 & 1 & 2 & 3 \\
2 & 2 & 4 & 5 \\
3 & 4 & 6 & 7 \\
5 & 6 & 7 & 8
\end{ytableau}, \ \begin{ytableau}
1 & 1 & 2 & 3 \\
2 & 2 & 4 & 5 \\
3 & 4 & 6 & 7 \\
5 & 6 & 8 & 8
\end{ytableau}, \ \begin{ytableau}
1 & 1 & 2 & 3 \\
2 & 2 & 4 & 5 \\
3 & 4 & 7 & 7 \\
5 & 6 & 8 & 8
\end{ytableau}, \ \begin{ytableau}
1 & 1 & 3 & 3 \\
2 & 2 & 4 & 5 \\
3 & 4 & 6 & 7 \\
5 & 6 & 8 & 8
\end{ytableau},  \begin{ytableau}
1 & 1 & 3 & 3 \\
2 & 2 & 4 & 5 \\
3 & 4 & 7 & 7 \\
5 & 6 & 8 & 8
\end{ytableau}, } 
\end{align*}
\begin{align*}
& \scalemath{0.6}{ \begin{ytableau}
1 & 1 & 1 & 2 \\
2 & 3 & 4 & 4 \\
3 & 5 & 6 & 6 \\
5 & 7 & 8 & 8
\end{ytableau}, \ \begin{ytableau}
1 & 1 & 1 & 2 \\
2 & 3 & 4 & 4 \\
3 & 5 & 6 & 7 \\
5 & 7 & 8 & 8
\end{ytableau}, \ \begin{ytableau}
1 & 1 & 2 & 2 \\
2 & 3 & 4 & 4 \\
3 & 5 & 6 & 6 \\
5 & 7 & 7 & 8
\end{ytableau}, \ \begin{ytableau}
1 & 1 & 2 & 2 \\
2 & 3 & 4 & 4 \\
3 & 5 & 6 & 6 \\
5 & 7 & 8 & 8
\end{ytableau}, \ \begin{ytableau}
1 & 1 & 2 & 2 \\
2 & 3 & 4 & 4 \\
3 & 5 & 6 & 7 \\
5 & 7 & 8 & 8
\end{ytableau}, \ \begin{ytableau}
1 & 1 & 2 & 3 \\
2 & 3 & 4 & 4 \\ 
3 & 5 & 6 & 6 \\
5 & 7 & 8 & 8
\end{ytableau}, \ \begin{ytableau}
1 & 1 & 2 & 3 \\
2 & 3 & 4 & 4 \\
3 & 5 & 6 & 7 \\
5 & 7 & 8 & 8
\end{ytableau}, \ \begin{ytableau}
1 & 1 & 2 & 2 \\
2 & 3 & 4 & 4 \\
3 & 6 & 6 & 7 \\
5 & 7 & 8 & 8
\end{ytableau}, \ \begin{ytableau}
1 & 1 & 2 & 2 \\
2 & 3 & 4 & 5 \\
3 & 6 & 6 & 7 \\
5 & 7 & 8 & 8
\end{ytableau}, } 
\end{align*}
\begin{align*}
& \scalemath{0.6}{ \begin{ytableau}
1 & 1 & 2 & 3 \\
2 & 3 & 4 & 4 \\
3 & 6 & 6 & 7 \\
5 & 7 & 8 & 8
\end{ytableau}, \ \begin{ytableau}
1 & 1 & 2 & 3 \\
2 & 3 & 4 & 5 \\
3 & 6 & 6 & 7 \\
5 & 7 & 8 & 8
\end{ytableau}, \ \begin{ytableau}
1 & 1 & 2 & 2 \\
2 & 4 & 4 & 5 \\
3 & 6 & 6 & 7 \\
5 & 7 & 8 & 8
\end{ytableau}, \ \begin{ytableau}
1 & 1 & 2 & 3 \\
2 & 4 & 4 & 5 \\
3 & 6 & 6 & 7 \\
5 & 8 & 8 & 8
\end{ytableau}, \ \begin{ytableau}
1 & 1 & 2 & 3 \\
2 & 2 & 4 & 5 \\
3 & 4 & 7 & 7 \\
6 & 6 & 8 & 8
\end{ytableau}, \ \begin{ytableau}
1 & 1 & 2 & 3 \\
2 & 2 & 5 & 5 \\
3 & 4 & 7 & 7 \\
6 & 6 & 8 & 8
\end{ytableau}, \ \begin{ytableau}
1 & 1 & 3 & 3 \\
2 & 2 & 4 & 5 \\
3 & 4 & 7 & 7 \\
6 & 6 & 8 & 8
\end{ytableau}, \ \begin{ytableau}
1 & 1 & 1 & 3 \\
2 & 2 & 3 & 5 \\
4 & 4 & 5 & 7 \\
6 & 6 & 7 & 8
\end{ytableau}, \ \begin{ytableau}
1 & 1 & 2 & 3 \\
2 & 2 & 5 & 5 \\
4 & 4 & 7 & 7 \\
6 & 6 & 8 & 8
\end{ytableau}. }
\end{align*}

From these, we get the following profiles of rank $4$ rigid indecomposable modules in ${\rm CM}(B_{4,8})$: 
\begin{align*}
\scalemath{0.7}{
\ffrac{1357}{1246}{2468}{1357}, \ 
 \ffrac{1357}{1246}{2468}{1358}, \ 
 \ffrac{1357}{1247}{2468}{1358}, \ 
 \ffrac{1357}{1246}{2468}{2357}, \ 
 \ffrac{1357}{1246}{2468}{2358}, \ 
 \ffrac{1357}{1247}{2468}{2358}, \  
 \ffrac{1357}{1346}{2468}{2358}, \ 
\ffrac{1357}{1347}{2468}{2358}, \ 
\ffrac{1246}{2468}{1358}{1357}, \ 
 \ffrac{1247}{2468}{1358}{1357}, \ 
 \ffrac{1246}{2468}{2357}{1357}, \ 
 \ffrac{1246}{2468}{2358}{1357}, \ 
 \ffrac{1247}{2468}{2358}{1357}, }
 \end{align*}
 \begin{align*}
\scalemath{0.7}{
 \ffrac{1346}{2468}{2358}{1357}, \ 
 \ffrac{1347}{2468}{2358}{1357}, \ 
\ffrac{1247}{2468}{2368}{1357}, \ 
 \ffrac{1257}{2468}{2368}{1357}, \
\ffrac{1347}{2468}{2368}{1357}, \ 
 \ffrac{1357}{2468}{2368}{1357}, \ 
 \ffrac{1257}{2468}{2468}{1357}, \ 
 \ffrac{1357}{2468}{2468}{1358}, \ 
 \ffrac{1357}{1247}{2468}{2368}, \ 
 \ffrac{1357}{1257}{2468}{2368}, \ 
 \ffrac{1357}{1347}{2468}{2368}, \ 
 \ffrac{1357}{1357}{1246}{2468}, \ 
 \ffrac{1357}{1257}{2468}{2468}. }
\end{align*}
We conjecture that these are all rigid indecomposable modules of rank $4$ in ${\rm CM}(B_{4,8})$ (up to cyclic shifts). 
\subsection{Quasi-homomorphisms of cluster algebras and braid group actions} \label{subsec:Quasi-homomorphisms of cluster algebras and braid group actions} 

The concept of quasi-homomorphism of cluster 
algebras was introduced by Fraser, see~\cite[Definition 3.1]{Fra16}. It is a map between cluster algebras {sharing the same principal part of the initial exchange matrix but may differ in their coefficients.} 
Let $\mathcal{A}$ be a cluster algebra with labeled tree $\TT_n$. 
If ${\bf x}_t$ is the cluster of $t\in\TT_n$, we write 
$x_{i;t}$ for its variables. 
% and 
% $x_{i;t'}$ for the variables of the cluster ${\bf x}_{t'}$ of $\mathcal{A}'$. 
Similarly, we write 
$\widehat{y}_{i;t}$ 
%and $\widehat{y}_{i;t'}$ 
for the variables from (\ref{eq:definition of y hat}). 

Let $x$ and $y$ be elements of $\mathcal{A}$. We say that $x$ is proportional to $y$, written as $x \propto y$, if $x = a y$ for some Laurent monomial $a$ in frozen variables. Let $\mathcal{A}$ and $\mathcal{A}'$ be two cluster algebras with $n$ cluster variables in each seed, and with respective groups $\mathbb{P}$ and $\mathbb{P}'$ of Laurent monomials in frozen variables (as in Section~\ref{subsec:cluster algebras}).
% such that there is an isomorphism between the 
% %, and there is an isomorphism $t \mapsto t'$ 
% labeled tree $\TT_n$ of $\mathcal{A}$ and the labeled tree $\TT_n'$ of $\mathcal{A}'$ (i.e. the mutable parts of $\mathcal{A}$ and 
% $\mathcal{A}'$ are the same).
% Let $t'\in\TT'_n$ be the image of $t\in\TT_n$ under the isomorphism. 
An algebra homomorphism $f: \mathcal{A} \to \mathcal{A}'$ is called a \textit{quasi-homomorphism} if $f(\mathbb{P}) \subset \mathbb{P}'$, and if there exists a seed $\Sigma_t=({\bf x}_t, \widetilde{B}_t)$ for $\mathcal{A}$ and a seed $\Sigma_{t'}=({\bf x}_{t'}, \widetilde{B}_{t'})$ for $\mathcal{A}'$ such that
\begin{enumerate}[(i)]
\item $f(x_{i;t}) \propto x_{i;t'}$, $1 \le i \le n$,

\item $f(\widehat{y}_{i;t}) = \widehat{y}_{i; t'}$, $1 \le i \le n$, 

\item $B_t = B_{t'}$,
\end{enumerate} 
where $B_t$ and $B_{t'}$ are the upper $n \times n$ submatrices of $\widetilde{B}_t$ and $\widetilde{B}_{t'}$ respectively. A \textit{quasi-automorphism} is a quasi-homomorphism $f: \mathcal{A} \to \mathcal{A}$ such that there is a quasi-homomorphism $g: \mathcal{A} \to \mathcal{A}$ such that for any cluster variable $x$, $(g\circ f)(x) \propto x$. 

Fraser showed that quasi-homomorphisms send clusters to clusters up to factors in frozen variables, \cite{Fra20}.
Furthermore, he showed that 
%Quasi-cluster morphisms are categorified in the Appendix A.5 of \cite{KW23}. 
certain quasi-automorphisms on $\CC[\Gr(k,n)]$ satisfy braid relations and thus form a braid group. 
We now recall the definition of these 
quasi-automorphisms. See \cite[\S 5]{Fra20} for more details. 

Let $V$ be a $k$-dimensional 
complex vector space. An $n$-tuple of vectors $(v_1, \ldots, v_n)$ is called {\em consecutively generic} if every cyclically consecutive $k$-tuple of vectors is linearly independent, i.e., $\det( v_{i+1}, \ldots, v_{i+k} ) \ne 0$ for $i=1,\ldots,n$ where the indices are treated modulo $n$. 
Denote by $(V^n)^{\circ} \subset V^n$ the quasi-affine variety consisting of consecutively generic $n$-tuples.

Let $d={\rm gcd}(k,n)$. For $i \in [1,d-1]$, the map $\sigma_i: (V^n)^{\circ} \to (V^n)^{\circ}$ is defined as follows. 
{Every element $(v_1, \ldots, v_n)$ of $(V^n)^{\circ}$ gets divided into $\frac{n}{d}$ windows:
$[v_{1+jd}, \ldots, v_{(j+1)d}]$, $j \in [0, \frac{n}{d}-1]$. }
We note that it is enough to define $\sigma_i$ on one of the windows as it satisfies $d$-periodicity: If 
$\rho$ is the so-called twisted cyclic shift on $V^n$, given by  
$(v_1,\dots, v_n)\stackrel{\rho}{\mapsto} (v_2,\dots, v_n, (-1)^{k-1}v_1)$ then $\sigma_i\circ \rho^d=\rho^d\circ\sigma_i$.

The map $\sigma_i$ sends the first window (where $j=0$) to $[v_1, \ldots, v_{i-1}, v_{i+1}, w_1, v_{i+2}, \ldots, v_d]$, 
where 
\begin{align} \label{eq:braid group action formula}
    w_1 = \frac{\det(v_i, v_{i+2}, \ldots, v_{i+k})}{\det( v_{i+1}, v_{i+2}, \ldots, v_{i+k} )} v_{i+1} - v_{i}.
\end{align} 

The $\ell$th window is defined by the same recipe by $d$-periodically augmenting indices. 
The pullback $\sigma_i^*$ (with a slight abuse of notation, we also denote it as $\sigma_i$) is a quasi-automorphisms on $\CC[\Gr(k,n)]$. We will use this 
in particular when we deal with non-real 
elements in the dual canonical 
basis in Sections~\ref{subsec:non-rigid modules Gr39} and~\ref{subsec:non-rigid modules Gr48}.

\begin{example}\label{ex:sigmas-for-39}
If $(k,n)=(3,9)$, we have $d=3$. 
We write out the effect of the maps $\sigma_1,\sigma_2$ on the first window $[v_1,v_2,v_3]$ of any $9$-tuple 
$(v_1,\dots, v_9)\in (V^9)^\circ$ of consecutively generic vectors: 
\begin{eqnarray*}
\sigma_1: [v_1,v_2,v_3] & \mapsto & 
[v_2,\frac{\det(v_1,v_3,v_4)}{\det(v_2,v_3,v_4)}v_2-v_1,v_3], \\
\sigma_2: [v_1,v_2,v_3] & \mapsto & 
[v_1,v_3,\frac{\det(v_2,v_4,v_5)}{\det(v_3,v_4,v_5)}v_3-v_2].
\end{eqnarray*}
\end{example}

\begin{remark}
By Theorem 5.3 in~\cite{Fra20}, the 
maps $\sigma_i$ where $1\le i\le d-1$ and $d=\gcd(k,n)$ are 
algebra homomorphisms on $\CC[{\rm Gr}(k,n)]$ satisfying the braid relations and so they induce an action of the braid 
group ${\rm Br}_d$ on $\CC[{\rm Gr}(k,n)]$. 
This in turn induces a braid group 
action on the subcategory $\mathscr{C}_\ell=\mathscr{C}_\ell^{\mathfrak{sl}_k}$ 
of the 
category of finite dimensional representations of $U_q(\widehat{\mathfrak{sl}_k})$ 
(see Section~\ref{subsec:correspondence between modules and tableaux}) with $n=k+\ell+1$, using the isomorphism 
between its Grothendieck ring and the cluster algebra $\CC[{\rm Gr}(k,n,\sim)]$.

The braid group action Fraser defined coincides with the one from~\cite[Theorem 2.3]{KKOP21} on $\mathscr{C}_{\ell}$, see \cite[Theorem 5.5]{LP2024braid}. This action sends simple modules to 
simple modules, \cite[\S 1.1]{KQW22}. Combining Fraser's approach with the one of \cite{KKOP21}, up to factors arising from frozen variables, the braid group action should send prime modules to prime modules and real (resp. non-real) modules to real (resp. non-real) modules. 
\end{remark}

In the remaining sections, we study non-rigid indecomposable modules in ${\rm CM}(B_{3,9})$ and ${\rm CM}(B_{4,8})$ (Sections \ref{subsec:non-rigid modules Gr39} and \ref{subsec:non-rigid modules Gr48}), using non-real semistandard Young tableaux and Fraser's braid group action on them.

\subsection{Indecomposable non-rigid modules in ${\rm CM}(B_{3,9})$} \label{subsec:non-rigid modules Gr39}

There are three prime non-real tableaux in $\SSYT(3,[9])$ with $3$ columns:
\begin{align*}
\scalemath{0.6}{
\bT_1 = \begin{ytableau}
1 & 2 & 3 \\
4 & 5 & 6 \\
7 & 8 & 9
\end{ytableau}, \quad 
\bT_2 = \begin{ytableau}
1 & 2 & 5 \\
3 & 4 & 8 \\
6 & 7 & 9
\end{ytableau}, \quad 
\bT_3 = \begin{ytableau}
1 & 3 & 4 \\
2 & 6 & 7 \\
5 & 8 & 9
\end{ytableau}, }
\end{align*}
see Section 8 in \cite{CDFL20}. 
We take the initial seed of $\CC[\Gr(3,9)]$ from Section~\ref{subsec:non-rigid-module-example-Gr39}, with $3$-subsets:
\begin{align} \label{eq:initial seed Gr39}
\begin{split}
& 124, 125, 126, 127, 128, 134, 145, 156, 167, 178, \\
& 123, 234, 345, 456, 567, 678, 789, 129, 189,
\end{split}
\end{align}
the second row contains the frozen variables. 
In terms of this seed,
the extended $\mathbf{g}$-vectors of $\bT_1, \bT_2, \bT_3$ are 
\begin{align*}
& {\bf \hat{g}}_{1} = (0, -1, -1, 0, 1, -1, 0, 1, 1, 0, 0, 1, 1, 0, 0, 0, 0, 1, 0),  \\
& {\bf \hat{g}}_2 = (-1, -1, 1, 1, 0, 0, 1, 0, -1, 0, 0, 1, 0, 0, 1, 0, 0, 0, 1), \\
& {\bf \hat{g}}_3 = (0, 1, 0, -1, 0, 0, -1, -1, 1, 1, 0, 0, 1, 1, 0, 0, 0, 1, 0),
\end{align*}
respectively. Let $\mathbf{g}_i={\rm pr}_n(\hat{\mathbf{g}}_i)$ for $i=1,2,3$.

These three extended $\mathbf{g}$-vectors are non-real. As argued in 
Section~\ref{subsec:non-rigid-module-example-Gr39}, modules in ${\rm CM}(B_{3,9})$ with these extended $\mathbf{g}$-vectors are non-rigid. 

There is an open dense subset $\mathcal{O}_{\mathbf{g}_1}$ of $\Hom_{\underline{\rm CM}(B_{3,9})}(T^{{\mathbf{g}_1}_-},T^{{\mathbf{g}_1}_+})$ such that for any $f\in \mathcal{O}_{\mathbf{g}_1}$, $\operatorname{cone}(f)$ has the same $F$-polynomial (cf. Section \ref{ss:function-e}). We refer to each such module $\operatorname{cone}(f)$ as a generic module of ${\mathbf{g}}_1$. We conjecture that {a generic 
module of ${\mathbf{g}}_1$} has profile $\scalemath{0.6}{\thfrac{369}{258}{147}}$. Similarly, modules in ${\rm CM}(B_{3,9})$ corresponding to ${\bf {g}}_2$ and ${\bf {g}}_3$ are also non-rigid. We conjecture that {generic modules of $\mathbf{g}_2$ and $\mathbf{g}_3$} have profiles $\scalemath{0.6}{\thfrac{258}{147}{369}}$ and $\scalemath{0.6}{\thfrac{147}{369}{258}}$ respectively. We expect that these modules are indecomposable. %\textcolor{red}{can we show that they are indecomposable?}

We now use the braid group action to find new candidates for non-rigid indecomposable modules.
The braid group for $\CC[\Gr(3,9)]$ is ${\rm Br}_3$, it has two generators $\sigma_1, \sigma_2$, see 
formula (\ref{eq:braid group action formula}) and Example~\ref{ex:sigmas-for-39}. Using the correspondence between dual canonical basis elements and tableaux \cite{CDFL20} and the braid group action (\ref{eq:braid group action formula}) in Section \ref{subsec:Quasi-homomorphisms of cluster algebras and braid group actions}, 
one obtains: 
\begin{align*}
\sigma_2(\bT_1) = \scalemath{0.6}{ \begin{ytableau}
1 & 1 & 2 & 2 & 3 & 6 \\
3 & 4 & 4 & 5 & 5 & 8 \\
6 & 7 & 7 & 8 & 9 & 9
\end{ytableau} }\,, \ \sigma_1(\bT_2) = \scalemath{0.6}{ \begin{ytableau}
1 & 1 & 2 & 4 & 4 & 5 \\
2 & 3 & 3 & 7 & 7 & 8 \\
5 & 6 & 6 & 8 & 9 & 9
\end{ytableau}}\,, \ \sigma_2\sigma_1^2(\bT_2) = \scalemath{0.6}{ \begin{ytableau}
1 & 1 & 2 & 3 & 3 & 4 \\
2 & 5 & 5 & 6 & 6 & 7 \\
4 & 7 & 8 & 8 & 9 & 9
\end{ytableau}}.
\end{align*}
These tableaux are all prime non-real. They have extended $\mathbf{g}$-vectors:
\begin{align*}
& {\bf \hat{g}}( \sigma_2(\bT_1) ) = (-1, -2, -1, 2, 1, -1, 1, 2, 0, -1, 0, 2, 1, 0, 0, 1, 0, 1, 1), \\
& {\bf \hat{g}}( \sigma_1(\bT_2) ) = (-2, 1, 2, 0, -1, 1, 0, -2, -1, 2, 0, 1, 0, 2, 1, 0, 0, 1, 1), \\
& {\bf \hat{g}}( \sigma_2\sigma_1^2(\bT_2) ) = (1, 0, -2, -1, 1, -1, -2, 1, 2, 1, 0, 1, 2, 1, 0, 0, 0, 2, 0),
\end{align*}
respectively, with respect to the initial seed (\ref{eq:initial seed Gr39}). 
{Conjecturally,} the associated modules in ${\rm CM}(B_{3,9})$ have the following profiles: 
\begin{align*}
\begin{array}{c}
   369 \\
   \hline
   258 \\ 
   \hline
   258 \\ 
   \hline
   147 \\ 
   \hline
   147 \\ 
  \hline
   369  
\end{array}, \quad \begin{array}{c} 
   258 \\ 
   \hline
   147 \\ 
   \hline
   147 \\ 
  \hline
   369  \\
   \hline
   369 \\
   \hline
   258 
\end{array}, \quad \begin{array}{c} 
   147 \\ 
  \hline
   369  \\
   \hline
   369 \\
   \hline
   258 \\
   \hline
   258 \\ 
   \hline
   147  
\end{array}. 
%\scalemath{0.7}{ \sfrac{369}{258}{258}{147}{147}{369}, \ \sfrac{258}{147}{147}{369}{369}{258}, \ \sfrac{147}{369}{369}{258}{258}{147}. }
\end{align*}
respectively. 
We expect that they are indecomposable and non-rigid. 

% \begin{remark}
% There are indecomposable non-rigid modules in ${\rm CM}(B_{3,9})$ of arbitrary high rank. 
% \end{remark}

%%
%
\subsection{Indecomposable non-rigid modules in ${\rm CM}(B_{4,8})$}  \label{subsec:non-rigid modules Gr48}
There are two prime non-real tableaux in $\SSYT(4,[8])$ with $2$ columns:
\begin{align*} 
\bT_1 = \scalemath{0.6}{ \begin{ytableau}
1 & 2 \\
3 & 4 \\
5 & 6 \\
7 & 8
\end{ytableau} }, \quad 
\bT_2 = \scalemath{0.6}{  \begin{ytableau}
1 & 3 \\
2 & 5 \\
4 & 7 \\
6 & 8
\end{ytableau} }. 
\end{align*}
Take the following initial seed of $\CC[\Gr(4,8)]$ with the following $4$-subsets: 
\begin{align} \label{eq:initial seed Gr48}
\begin{split}
& 1235, 1236, 1237, 1245, 1256, 1267, 1345, 1456, 1567, \\
& 1234, 2345, 3456, 4567, 5678, 1238, 1278, 1678.
\end{split}
\end{align}
In terms of this seed, the extended $\mathbf{g}$-vectors of $\bT_1, \bT_2$ are
\begin{align*}
& {\bf \hat{g}}_1 = (0, -1, 0, -1, 0, 1, 0, 1, 0, 0, 1, 0, 0, 0, 1, 0, 0),  \\
& {\bf \hat{g}}_2 = (-1, 1, 0, 1, 0, -1, 0, -1, 1, 0, 0, 1, 0, 0, 0, 1, 0),
\end{align*}
respectively. 

As shown in Section \ref{subsec:non-rigid-module-example-Gr48}, the modules in ${\rm CM}(B_{4,8})$ corresponding to ${\bf g}_1$ are non-rigid. We conjecture that a generic module of $\mathbf{g}_1$ has profile $\scalemath{0.6}{\frac{1357}{2468}}$. Similarly, we can show that the modules in ${\rm CM}(B_{4,8})$ corresponding to ${\bf g}_2$ are also non-rigid. We conjecture that a generic module of $\mathbf{g}_2$ has profile $\scalemath{0.6}{\frac{2468}{1357}}$. 
Explicit 1-parameter families for indecomposable rank 2 modules with this profile were given in~\cite{BBL23}. 

The braid group for $\CC[\Gr(4,8)]$ is ${\rm Br}_4$ and it has three generators $\sigma_1$, $\sigma_2$, $\sigma_3$. Using the formula (\ref{eq:braid group action formula}) and the correspondence between dual canonical basis elements and tableaux \cite{CDFL20}, one obtains:  
\begin{align*}
& \sigma_1(\bT_1) = \scalemath{0.6}{ \begin{ytableau}
1 & 1 & 2 & 4 \\
2 & 3 & 3 & 6 \\
4 & 5 & 5 & 7 \\
6 & 7 & 8 & 8
\end{ytableau} }, \ \sigma_2(\bT_2) = \scalemath{0.6}{ \begin{ytableau}
1 & 1 & 2 & 3 \\
2 & 4 & 4 & 5 \\
3 & 6 & 6 & 7 \\
5 & 7 & 8 & 8
\end{ytableau}}, \ \sigma_1^2(\bT_1) = \scalemath{0.6}{ \begin{ytableau}
1 & 1 & 1 & 2 & 3 & 4 \\
2 & 2 & 3 & 3 & 5 & 6 \\
4 & 4 & 5 & 5 & 7 & 7 \\
6 & 6 & 7 & 8 & 8 & 8
\end{ytableau}},  
\end{align*}
\begin{align*} 
& \sigma_2^2(\bT_2) = \scalemath{0.6}{ \begin{ytableau}
1 & 1 & 1 & 2 & 2 & 3 \\
2 & 3 & 4 & 4 & 4 & 5 \\
3 & 5 & 6 & 6 & 6 & 7 \\
5 & 7 & 7 & 8 & 8 & 8
\end{ytableau}},  \ \sigma_1 \sigma_2^2(\bT_2) = \scalemath{0.6}{ \begin{ytableau}
1 & 1 & 1 & 2 & 2 & 4 \\
2 & 3 & 3 & 3 & 4 & 6 \\
4 & 5 & 5 & 5 & 6 & 7 \\
6 & 7 & 7 & 8 & 8 & 8
\end{ytableau}},  \ \sigma_3 \sigma_2 \sigma_1^2 \sigma_2^2(\bT_2) = \scalemath{0.6}{ \begin{ytableau}
1 & 1 & 1 & 2 & 3 & 3 \\
2 & 2 & 4 & 4 & 5 & 5 \\
3 & 4 & 6 & 6 & 7 & 7 \\
5 & 6 & 7 & 8 & 8 & 8
\end{ytableau}}.
\end{align*}
These tableaux are all prime non-real. They have extended $\mathbf{g}$-vectors:
\begin{align*}
& {\bf \hat{g}}( \sigma_1(\bT_1) ) = (0, -1, 0, -1, 0, 1, 0, 1, 0, 0, 1, 0, 0, 0, 1, 0, 0), \\
& {\bf \hat{g}}( \sigma_2(\bT_2) ) =  (0, -1, 0, -1, 0, 1, 0, 1, 0, 0, 1, 0, 0, 0, 1, 0, 0), 
\end{align*}
\begin{align*} 
& {\bf \hat{g}}( \sigma_1^2(\bT_1) ) =  (-2, 0, 1, 0, 2, -1, 1, -1, 0, 0, 1, 1, 1, 0, 1, 1, 1), \\
& {\bf \hat{g}}( \sigma_2^2(\bT_2) ) =  (1, -1, -1, -1, -2, 2, -1, 2, 1, 0, 2, 1, 0, 0, 2, 1, 0), 
\end{align*}
\begin{align*} 
& {\bf \hat{g}}( \sigma_1 \sigma_2^2(\bT_2) ) =  (-1, -2, 1, -2, 2, 1, 1, 1, -1, 0, 2, 0, 1, 0, 2, 0, 1), \\
& {\bf \hat{g}}( \sigma_3 \sigma_2 \sigma_1^2 \sigma_2^2(\bT_2) ) =  (0, 1, -1, 1, -2, 0, -1, 0, 2, 0, 1, 2, 0, 0, 1, 2, 0),
\end{align*}
respectively, with respect to the initial seed (\ref{eq:initial seed Gr48}). 
The associated modules in ${\rm CM}(B_{4,8})$ have profiles: 
\begin{align*}
\begin{array}{c}
   2468 \\
   \hline
   1357 \\ 
   \hline
   1357 \\ 
   \hline           
   2468 
\end{array}, \quad \begin{array}{c} 
   1357 \\ 
   \hline
   2468 \\ 
   \hline
   2468 \\   
  \hline
   1357 
\end{array}, \quad \begin{array}{c} 
   2468 \\ 
  \hline
   1357  \\
   \hline
   1357 \\
   \hline     
   1357 \\
   \hline
   2468 \\ 
   \hline
   2468  
\end{array}, \ \quad \begin{array}{c} 
   1357 \\ 
  \hline
   2468  \\
   \hline
   2468 \\
   \hline     
   2468 \\
   \hline
   1357 \\ 
   \hline
   1357  
\end{array}, \ \quad \begin{array}{c} 
   2468 \\ 
  \hline
   2468  \\
   \hline
   1357 \\
   \hline     
   1357 \\
   \hline
   1357 \\ 
   \hline
   2468  
\end{array}, \ \quad \begin{array}{c} 
   1357 \\ 
  \hline
   1357  \\
   \hline
   2468 \\
   \hline     
   2468 \\
   \hline
   2468 \\ 
   \hline
   1357  
\end{array},
\end{align*}
respectively. We conjecture that they are indecomposable and non-rigid.

\bibliographystyle{acm}
\bibliography{biblio-CJK}

\end{document}